\documentclass[reqno, 11pt, a4paper]{amsart}

\usepackage{latexsym,ifthen,xspace}
\usepackage{amsmath,amssymb,amsthm}
\usepackage{enumerate, calc}
\usepackage{bbm}
\usepackage[colorlinks=true, pdfstartview=FitV, linkcolor=blue,
  citecolor=blue, urlcolor=blue, pagebackref=false]{hyperref}
\usepackage[text={33pc,605pt},centering]{geometry}    
\usepackage{orcidlink}
\usepackage{graphicx}

\newtheorem{theorem}{Theorem}[section]
\newtheorem{lemma}[theorem]{Lemma}
\newtheorem{prop}[theorem]{Proposition}
\newtheorem{assumption}[theorem]{Assumption}
\newtheorem{cor}[theorem]{Corollary}

\theoremstyle{definition}
\newtheorem{definition}[theorem]{Definition}
\newtheorem{example}[theorem]{Example}

\theoremstyle{remark}
\newtheorem{remark}[theorem]{Remark}

\numberwithin{equation}{section}

\DeclareMathAlphabet{\mathsl}{OT1}{cmss}{m}{sl}
\SetMathAlphabet{\mathsl}{bold}{OT1}{cmss}{bx}{sl}

\newcommand{\me}{\ensuremath{\mathrm{e}}}
\newcommand{\md}{\ensuremath{\mathrm{d}}}

\newcommand{\scpr}[3]{%
  \ensuremath{%
    \big\langle
      #1, #2
    \big\rangle_{\raisebox{-0ex}{$\scriptstyle \ell^{\raisebox{.1ex}{$\scriptscriptstyle 2$}} (#3)$}}
  }
}

\DeclareMathOperator{\mean}{\mathbb{E}}
\DeclareMathOperator{\Mean}{\mathrm{E}}
\DeclareMathOperator{\prob}{\mathbb{P}}
\DeclareMathOperator{\Prob}{\mathrm{P}}
\DeclareMathOperator{\var}{\mathbb{V}ar}
\DeclareMathOperator{\cov}{\mathbb{C}ov}

\DeclareMathOperator{\supp}{\mathrm{supp}}
\DeclareMathOperator{\dom}{\mathrm{dom}}

\newcommand{\ldef}{\ensuremath{\mathrel{\mathop:}=}}
\newcommand{\rdef}{\ensuremath{=\mathrel{\mathop:}}}

\newcommand{\indicator}{\mathbbm{1}}

\makeatletter
\newcommand\initCounter[1]{%
  \@ifundefined{c@constcnt-#1}%
    {\newcounter{constcnt-#1}}%
    {}%
}

\newcommand\const[2][c]{%
  \initCounter{#1}%
  \@ifundefined{const-#1-#2}%
    {%
      \stepcounter{constcnt-#1}%
      \expandafter\xdef\csname const-#1-#2\endcsname{\arabic{constcnt-#1}}%
    }%
    {}%
  #1_{\csname const-#1-#2\endcsname}%
}
\makeatother

\begin{document}

\title[Gradient estimates of the heat kernel in the RCM]{Gradient estimates of the heat kernel for random walks among time-dependent random conductances}

%
%

\author[J.-D.~Deuschel]{Jean-Dominique Deuschel}
\address{Technische Universit\"at Berlin}
\curraddr{Strasse des 17. Juni 136, 10623 Berlin, Germany}
\email{deuschel@math.tu-berlin.de}
\thanks{}

\author[T.~Kumagai]{Takashi Kumagai \orcidlink{0000-0001-7515-1055}}
\address{Waseda University}
\curraddr{3-4-1 Okubo, Shinjuku-ku, Tokyo 169-8555, Japan}
\email{t-kumagai@waseda.jp}
\thanks{}

\author[M.~Slowik]{Martin Slowik \orcidlink{0000-0001-5373-5754}}
\address{University of Mannheim}
\curraddr{Mathematical Institute, B6, 26, 68159 Mannheim}
\email{slowik@math.uni-mannheim.de}
\thanks{}


\date{\today}

\dedicatory{}

\begin{abstract}
  In this paper we consider a time-continuous random walk in $\mathbb{Z}^{d}$ in a dynamical random environment with symmetric jump rates to nearest neighbours.  We assume that these random conductances are stationary and ergodic and, moreover, that they are bounded from below but unbounded from above with finite first moment.  We derive sharp on-diagonal estimates for the annealed first and second discrete space derivative of the heat kernel which then yield local limit theorems for the corresponding kernels.  Assuming weak algebraic off-diagonal estimates, we then extend these results to the annealed Green function and its first and second derivative.  Our proof which extends the result of \cite{DD05} to unbounded conductances with first moment only, is an adaptation of the recent entropy method of \cite{BD-CKY15}.
\end{abstract}

\maketitle

\tableofcontents

\section{Introduction}
It is well known that the heat kernel associated with a uniform elliptic operator in divergence form possesses favourable regularity properties, particularly H\"older continuity in both time and space variables -- the so-called De Giorgi-Nash-Moser theory, cf.~\cite{De57, Na58, Mo61, Mo64}. These results can be extended to the discrete setting for the heat kernels of random walks with symmetric jump rates, uniformly in the diffusive scaling parameter, cf.~\cite{De99}.  However, unless the parameters are sufficiently smooth, the situation is markedly different for the gradient of the heat kernel.  In particular, in the random media setting, one cannot expect good estimates in the quenched regime, i.e., for almost all realisations of the environment, but only in the annealed regime, i.e., averaged over the law of the environment, cf.~\cite{DD05, GNO15, MO15, CGO17}.  The analytical methods used in \cite{DD05} and \cite{GNO15, MO15, CGO17} rely heavily on the uniform ellipticity assumption, rendering them ill-suited for degenerate situations such as unbounded coefficients or random walks on random percolation clusters.  Under nice mixing conditions on the environment, including supercritical percolation models, estimates for the gradient of the Green's function have recently been obtained even in the quenched setting, cf.~\cite{AD18, DG21, B-RCD25}, using a powerful technique of quantitative homogenisation, cf.~\cite{AKM19}.  In \cite{BD-CKY15}, a different approach is proposed based on a more robust entropy method for a discrete-time random walk on the supercritical percolation cluster.  They obtained very sharp annealed results for the discrete gradient of the corresponding heat kernel, which they utilised to prove a Liouville principle for sublinear harmonic functions.  Our objective in the present paper is to adapt the entropy method to the  variable speed random walk  with dynamical random conductances, bounded from below but unbounded from above.  Under the minimal assumption of time-space stationarity of the environment and a finite first moment of the conductance, we then obtain sharp, scaling-invariant annealed on-diagonal estimates for the first and second discrete derivative of the heat kernel, which we subsequently use to prove a local limit theorem for the annealed heat kernel and for the discrete first derivative.  Moreover, assuming further off-diagonal estimates, we can extend these annealed results to the Green function and the corresponding first and second discrete derivatives.

\subsection{The model}
Consider the $d$-dimensional Euclidean lattice, $(\mathbb{Z}^{d}, E_{d})$, for $d \geq 1$, where $E_{d}$ denotes the set of all oriented nearest neighbour bonds, i.e.\ $E_{d} \ldef \{ (x, y): x,y \in \mathbb{Z}^{d}, |x - y| = 1\}$, where $|\cdot|$ is the usual graph-distance on $\mathbb{Z}^{d}$.  We also write $x \sim y$ if $(x,  y) \in E_{d}$.  Moreover, by $B(x, n) \ldef \{y \in \mathbb{Z}^{d} : |x - y| \leq n\}$ we denote the closed ball with centre $x$ and radius $n$.  The graph $(\mathbb{Z}^{d}, E_{d})$ is endowed with a family $\omega = \{\omega_{t}(e) : e \in E_{d}, \, t\in \mathbb{R} \} \in \Omega \ldef (0, \infty)^{\mathbb{R} \times E_{d}}$ of time-dependent, positive weights.  We equip $\Omega$ with the Borel-$\sigma$-algebra, $\mathcal{F}$.

To simplify notation, for $x, y \in \mathbb{Z}^{d}$ and $t \in \mathbb{R}$ we set $\omega_{t}(x,y) = 0$ for $(x,y) \notin E_{d}$.  We introduce a \emph{time-space} shift $\tau_{s,z}$ by $(s, z) \in \mathbb{R} \times \mathbb{Z}^{d}$ through
\begin{align}
  \label{eq:shift:space-time}
  (\tau_{s, z\,} \omega)_{t}(x,y) \;\ldef\; \omega_{t+s}(x+z,y+z),
  \qquad \forall\, t \in \mathbb{R},\; (x, y) \in E_{d}.
\end{align}
Assume that, for every $A \in \mathcal{F}$, the mapping $(\omega,t,x) \mapsto \indicator_{A}(\tau_{t, x} \omega)$ is jointly measurable with respect to the $\sigma$-algebra $\mathcal{F} \otimes \mathcal{B}(\mathbb{R}) \otimes \mathcal{P}(\mathbb{Z}^{d})$.

Further, consider a probability measure, $\prob$, on the measurable space $(\Omega, \mathcal{F})$, and we write $\mean$ to denote the corresponding expectation with respect to $\prob$.
\begin{assumption}\label{ass:law}
  Assume that $\prob$ satisfies the following conditions:
  \begin{enumerate}[(i)]
  \item
    $\prob$ is stationary and ergodic with respect to time-space shifts.
    
  \item
    The time dependent random conductances are $\prob$-a.s.\ symmetric, that is,
    \begin{align*}
      \prob\bigl[\omega_{t}(x,y) = \omega_{t}(y,x) \quad \forall\, t\in \mathbb{R} \bigr]
      \;=\;
      1.
    \end{align*}

  \item
    $\mean\bigl[\omega_{t}(e)\bigr] < \infty$ for all $e \in E_{d}$ and $t \in \mathbb{R}$.
    
  \item
    For every $A \in \mathcal{F}$ the mapping $(\omega, t, x) \mapsto \mathbbm{1}_{A}(\tau_{t, x} \omega)$ is jointly measurable with respect to the $\sigma$-algebra $\mathcal{F} \otimes \mathcal{B}(\mathbb{R}) \otimes \mathcal{P}(\mathbb{Z}^{d})$.
  \end{enumerate}
\end{assumption}
\begin{remark}
  In fact the time-space ergodicity assumption is not necessary for most of our results on the annealed estimates, cf.\ \cite{DD05}.  However, it is important when we are dealing with limit theorems in Section~\ref{sec:alclt}.  Assumption~\ref{ass:law}-(ii) is essential, because it allows us to identify the dual process, and it is used in several crucial steps.  Assumption~\ref{ass:law}-(iii) is quite natural under general ergodic assumption and it also guarantees non-explosion of the corresponding process, cf. \cite[Lemma 4.1]{ACDS18}, otherwise explosion could occur in finite time, cf.\ \cite[Remark 6.6]{BD10}.  In particular, the joint measurability, as states in Assumption~\ref{ass:law}-(iv), along with the stationarity of $\prob$ with respect to time-space shifts implies, by Fubini's theorem, that, $\prob$-a.s., the conductances are locally integrable in time, that is,
  \begin{align}
    \label{eq:local:integrability}
    \prob\biggl[ \int_{I} \omega_{s}(e)\, \md s < \infty \biggr]
    \;=\;
    1,
    \qquad\text{for all finite interval } I \subset \mathbb{R},\,
    e \in E_{d}.
  \end{align}
\end{remark}
For any fixed realisation $\omega \in \Omega$, we consider the time-inhomogeneous Markov process, $X \equiv \bigl((X_{t} : t \geq s), \Prob_{s,x}^{\omega} : (s,x) \in \mathbb{R} \times \mathbb{Z}^{d}\bigr)$ in the random environment $\omega$, where $\Prob_{s,x}^{\omega}$ denotes the law of $X$ on $\mathcal{D}(\mathbb{R}, \mathbb{Z}^{d})$, the space of $\mathbb{Z}^{d}$-valued c\`{a}dl\`{a}g functions on $\mathbb{R}$, starting at time $s$ in the vertex $x$, i.e.\
\begin{align*}
  \Prob_{s,x}^{\omega}\bigl[X_{s} = x\bigr] \;=\; 1
  \qquad \prob\text{-a.s.}
\end{align*}
The time-dependent generator, $\mathcal{L}_{t}^{\omega}$, acts on bounded functions $f\colon \mathbb{Z}^{d} \to \mathbb{R}$ as
\begin{align}
  (\mathcal{L}_{t}^{\omega} f)(x)
  &\;\ldef\;
  \sum_{y: (x,y) \in E_{d}}\mspace{-9mu} \omega_{t}(x,y)\, 
  \bigl(f(y) \,-\, f(x) \bigr).
\end{align}
Clearly, under $\Prob_{s,x}^{\omega}$, the distribution of the Markov process $X$ depends explicitly on $s \in \mathbb{R}$, $x \in \mathbb{Z}^{d}$ and $\omega \in \Omega$.  
For any $s \in \mathbb{R}$, we denote by $(P_{s,t}^{\omega} : t \geq s)$ the Markov semi-group associated to the Markov process $X$, i.e.\ $(P_{s,t}^{\omega}f)(x) = \Mean_{s,x}^{\omega}[f(X_{t})]$ for any bounded function $f\colon \mathbb{Z}^{d} \to \mathbb{R}$, $s \leq t$ and $x \in \mathbb{Z}^{d}$.  Moreover, for any $x, y \in \mathbb{Z}^{d}$ and $s \leq t$, we denote by 
\begin{align*}
  p_{s,t}^{\omega}(x,y) \;\ldef\; \Prob_{s,x}^{\omega}\!\bigl[X_{t} = y\bigr]
  \qquad \text{and} \qquad
  \bar{p}_{s,t}(x,y) \;\ldef\; \mean\bigl[p_{s,t}^{\omega}(x,y)\bigr].
\end{align*}
the (quenched) transition density with respect to the counting measure (heat kernel of the so-called VSRW) and the annealed (averaged) transition density, respectively.  As a consequence of \eqref{eq:shift:space-time} we have
\begin{align}\label{eq:shift:hk}
  p_{s,t}^{\tau_{h, z} \omega}(x,y)
  \;=\;
  p_{s+h, t+h}^{\omega}(x+z, y+z)
  \qquad \forall\, h \in \mathbb{R},\, z \in \mathbb{Z}^{d}. 
\end{align}
In particular, in view of Assumption~\ref{ass:law}-(i), it holds
\begin{align*}
  \bar{p}_{s,t}(x,y)
  \;=\;
  \bar{p}_{0,t-s}(0,y-x)
  \;=\;
  \bar{p}_{s-t,0}(x-y,0).
\end{align*}
Further, define $\tilde{\omega}_{t}(e) \ldef \omega_{-t}(e)$ for any $t \in \mathbb{R}$ and $e \in E_{d}$.  Then, in view of Lemma~\ref{lem:time_inverse}, $p_{s,t}^{\tilde{\omega}}(x,y) = p_{-t,-s}^{\omega}(y,x)$ for all $x,y \in \mathbb{Z}^{d}$ and $s,t \in \mathbb{R}$ with $s \leq t$.  It is also easily check that, for any $s \leq t \leq u$,
\begin{align*}
  p_{s,u}^{\tilde{\omega}}(x,y)
  \;=\;
  \sum_{z \in \mathbb{Z}^{d}} p_{s,t}^{\tilde{\omega}}(x,z)\, p_{t,u}^{\tilde{\omega}}(z,y)
  \qquad \forall\, x, y \in \mathbb{Z}^{d}. 
\end{align*}
Thus, by setting $(P_{s,t}^{\tilde{\omega}}f)(x) \ldef \sum_{y \in \mathbb{Z}^{d}} p_{s,t}^{\tilde{\omega}}(x,y) f(y)$ for any bounded function $f\colon \mathbb{Z}^{d} \to \mathbb{R}$, $s \leq t$ and $x \in \mathbb{Z}^{d}$, $(P_{s,t}^{\tilde{\omega}}: t \geq s)$ is a Markov semi-group.  For any fixed $\omega \in \Omega$, we write $\tilde{X} \equiv \bigl((\tilde{X}_{t} : t \geq s), \Prob_{s,x}^{\tilde{\omega}}: (s,x) \in \mathbb{R} \times \mathbb{Z}^{d}\bigr)$ to denote the associated time-inhomogeneous Markov process in the random environment $\tilde{\omega}$.  Again, as a consequence of \eqref{eq:shift:space-time},
\begin{align*}
  p_{s,t}^{\tau_{-h,z} \tilde{\omega}}(x,y)
  \;=\;
  p_{s+h,t+h}^{\tilde{\omega}}(x + z,y + z)
  \qquad \forall\, h \in \mathbb{R},\, z \in \mathbb{Z}^{d}. 
\end{align*}
Finally, let $(P^{\omega}_{s,t})^{*}$ be the adjoint of $P^{\omega}_{s,t}$ in $\ell^{2}(\mathbb{Z}^{d})$.  Then, $P_{s,t}^{\tilde{\omega}} = (P_{-t,-s}^{\omega})^{*}$.

In the next lemma we collect some important properties of the stochastic process $X$, $\tilde{X}$ and the associated \emph{environmental processes} -- also known as \emph{processes as seen from the point of view of the particle}. Since the time-dependent weights, $\omega_{t}(e)$, are not bounded from above, the processes $X$ and $\tilde{X}$ may explode in finite time, cf.\ \cite[Remark~6.6.]{BD10} for an example. However, Assumption~\ref{ass:law}-(iii) is enough to exclude this.
\begin{lemma}\label{lem:non-explosion}
  Suppose that Assumption~\ref{ass:law} is satisfied.  Then,
  \begin{enumerate}[(i)]
  \item
    for $\prob$-a.e.\ $\omega$, $\Prob_{0, 0}^{\omega}$-a.s.\ $X$ and $\tilde{X}$ do not explode, that is, there are only finitely many jumps in finite time.  In particular, for any $s \leq t$
    \begin{align}\label{eq:p:sum}
      \sum_{y \in \mathbb{Z}^{d}} p_{s,t}^{\omega}(x,y)
      \;=\;
      1
      \;=\;
      \sum_{y \in \mathbb{Z}^{d}} p_{s,t}^{\omega}(y,x)
      \qquad \forall\, x \in \mathbb{Z}^{d}
      \qquad \prob\text{-a.s.}
    \end{align}
    
  \item
    the environmental processes $(\tau_{t,X_{t}}\omega \,:\, t \geq 0)$ and $(\tau_{-t, \tilde{X}_{t}} \omega : t \geq 0)$ are Markov, and $\prob$ is a stationary distribution for both processes. If, additionally, $\prob\bigl[\omega_{t}(e) > 0\bigr] = 1$ for every $e \in E_{d}$ and $t \in \mathbb{R}$ then $\prob$ is ergodic.
    %
  \end{enumerate}
\end{lemma}
Although the proof of the above lemma is fairly standard and follows from \cite[Lemma 4.1]{ACDS18}, Lemma~\ref{lem:time_inverse} and by similar arguments as in the proof of \cite[Lemma~2.4]{ADS15}, we will provide a self-consistent proof in Section~\ref{sec:alclt} for the convenience of the reader.
%
%
%
%
%
%

\subsection{Main results}
As a further consequence of Assumption~\ref{ass:law} the averaged mean displacement of the stochastic process, $(X_{t} : t \geq 0)$, behaves diffusively.  For this purpose, consider the diffusively rescaled process $X^{(n)} = (X_{t}^{(n)} : t \geq 0)$, that is,
\begin{align*}
  X_{t}^{(n)} \;\ldef\; \frac{1}{n} X_{t n^{2}},
  \qquad t \geq 0, n \in \mathbb{N},
\end{align*}
and write $\prob_{0,0}^{*}[G] \ldef \mean\bigl[\Prob_{0,0}^{\omega}[G]\bigr]$ for any $G \in \mathcal{B}(\mathcal{D}([0, \infty),\mathbb{R}^{d}))$ to denote the \emph{annealed probability measure}.
\begin{prop}[Mean displacement]\label{prop:mean:displacement}
  Suppose that Assumption~\ref{ass:law} holds.  Then, there exists $\const[C]{mean:displacement} < \infty$ such that the following holds for all $T > 0$,
  \begin{align}\label{eq:mean:displacement}
    \mean\biggl[
      \Mean_{0,0}^{\omega}\biggl[\sup_{0 \leq t \leq T}|X_{t}|\biggr] 
    \biggr]
    \;\leq\;
    \const[C]{mean:displacement} \sqrt{T}.
  \end{align}
  In particular, for any $t > 0$, the family of scaled displacements $\bigl\{X_{t}^{(n)} : n \in \mathbb{N}\bigr\}$ is tight under $\prob_{0, 0}^{*}$.
\end{prop}
Our further results rely on the following additional assumption.
\begin{assumption}\label{ass:integrability}
  There exists a non-random constant $\const[C]{lower:elliptic} > 0$ such that 
  \begin{align}
    \prob[\omega_{t}(x,y) \geq \const[C]{lower:elliptic}] \;=\; 1
    \qquad \text{for any } (x,y) \in E_{d}.
  \end{align}
\end{assumption}
Our main objective is to prove a spatial derivative estimate of the annealed heat kernel.  For a given function $f\colon \mathbb{Z}^{d} \to \mathbb{R}$ the discrete spatial derivative in direction $(0, x) \in E_{d}$ is defined by
\begin{align*}
  \nabla_{\!x} f(y)
  \;\ldef\;
  f(x+y)-f(y),
  \qquad  \forall\, y \in \mathbb{Z}^{d}.
\end{align*}
Notice that $\nabla_{\!x}$ is a bounded linear operator on $\ell^p(\mathbb{Z}^{d})$ for any $1 \leq p \leq \infty$.  For functions $f\colon\mathbb{Z}^{d} \times \mathbb{Z}^{d} \to \mathbb{R}$ we write $\nabla_{\!x}^1$ and $\nabla_{\!x}^{2}$, respectively, to denote that the discrete spatial derivative acts either on the first or the second variable of $f$.  For any $(0,x), (0,x') \in E_{d}$ and $y, y' \in \mathbb{Z}^{d}$, higher spatial derivatives are obtained iteratively via the formula
\begin{align}\label{def:second:deriv}
  \nabla^1_{\!x}\nabla^1_{\!x'}f(y,y')
  &\;\ldef\;
  \nabla^1_{\!x}f(x'+y,y') - \nabla^1_{\!x}f(y,y')
  \nonumber\\[.5ex]  
  &\;=\;
  f(x+x'+y,y') - f(x+y,y') - f(x'+y,y') + f(y,y')
\end{align}
(an analogue expression holds for $\nabla^{2}_{\!x} \nabla^{2}_{\!x'}f$), whereas the mixed derivative is given by
\begin{align}\label{def:mixed:deriv}
  \nabla^{1}_{\!x} \nabla^{2}_{\!x'} f(y,y')
  &\;\ldef\;
  \nabla^{1}_{\!x}f(y, x' + y') - \nabla^{1}_{\!x}f(y,y')
  \nonumber\\[.5ex]
  &\;=\;
  f(x+y,x'+y') - f(y,x'+y') - f(x+y,y') + f(y,y').
\end{align}
Note that $\nabla^{1}_{\!x} \nabla^{2}_{\!x'} f = \nabla^{2}_{\!x'} \nabla^{1}_{\!x} f$.
\begin{theorem}[Gradient estimates]\label{thm:gradient:est}
  Suppose that Assumption~\ref{ass:law} and \ref{ass:integrability} hold true.
  \begin{enumerate}[(i)]
    \item There exists $\const[C]{gradient:point} < \infty$ such that for any $(0,x) \in E_{d}$, $y,y' \in \mathbb{Z}^{d}$ and $t \geq 1$,
    \begin{align}\label{eq:gradient:on-diag}
      \mean\Bigl[
        \bigl| \nabla^1_{\!x}\, p_{0, t}^{\omega}(y, y') \bigr|^{2}
      \Bigr]^{\!1/2}
      \;\leq\;
      \const[C]{gradient:point}\, t^{-(d+1)/2}.
    \end{align}
    Moreover, the same estimate holds if $\nabla^1_{\!x}\, p_{0, t}^{\omega}(y, y')$ is replaced by $\nabla^{2}_{\!x}\, p_{0, t}^{\omega}(y, y')$.
    
    \item There exists $\const[C]{gradient:sum} < \infty$ such that, for any $p \in [1,2]$, $(0,x)\in E_{d}$, $y, y' \in \mathbb{Z}^{d}$ and $t \geq 1$
    \begin{align}\label{eq:gradient:sum:on-diag}
      \biggl(
        \sum_{y' \in \mathbb{Z}^{d}}
        \mean\Bigl[ 
          \bigl| \nabla^1_{\!x}\, p_{0, t}^{\omega}(y, y') \bigr|^p 
        \Bigr]
      \biggr)^{\!\!1/p}
      \;\leq\;
      \const[C]{gradient:sum}\, t^{-(1+d(1-1/p))/2}.
    \end{align}
    Moreover, the same estimate holds if $\nabla^1_{\!x}\, p_{0, t}^{\omega}(y, y')$ is replaced by $\nabla^{2}_{\!x}\, p_{0, t}^{\omega}(y, y')$.
    
    \item For any $(0,x), (0,x')\in E_{d}$, $y, y' \in \mathbb{Z}^{d}$, and $t \geq 1$,
    \begin{align}
      \bigl|
        \nabla^{2}_{\!x} \nabla^{2}_{\!x'}\, \bar{p}_{0,t}(y,y')
      \bigr|
      &\;\leq\;
      \const[C]{gradient:sum}\, t^{-(d+2)/2},
      \label{eq:hk:deriv:22}
      \\[.5ex]
      \mean\Bigl[
        \Bigl|\nabla^{1}_{-x}\nabla^{2}_{x'}\, p^{\omega}_{0,t}(y,y')\Bigr|
      \Bigr]
      &\;\leq\;
      \const[C]{gradient:sum}\, t^{-(d+2)/2}.
      \label{eq:hk:deriv:12}
    \end{align}
  \end{enumerate}
\end{theorem}
\begin{remark}
  In \cite[Equation~(1.6)]{DD05}, for time-dependent conductances that are \emph{uniformly elliptic}, i.e.\ there exists $c_{0} \in [1, \infty)$ such that $\prob\bigl[c_{0}^{-1} \leq \omega_{t}(e) \leq c_{0}\; \forall\, t \in \mathbb{R}\bigr] = 1$ for any edge $e \in E_{d}$, the following time-derivative estimate
  \begin{align*}
    \Bigl| \frac{\partial}{\partial t}\, \bar{p}_{0,t}(0,x) \Bigr|
    \;\leq\;
    c_{1}\, t^{-d/2-1},
    \qquad \forall\, x \in \mathbb{Z}^{d},\, t \geq 1
  \end{align*}
  for some $c_{1} \in (0, \infty)$, is derived from a spatial-derivative estimate analogous to that given in \eqref{eq:gradient:sum:on-diag}. Actually, in \cite[Theorem~1.1]{DD05}, an even stronger estimate has been proven. However, when the conductances are unbounded from above, the proof provided therein fails, as it fundamentally relies on the uniform ellipticity of the conductances.
\end{remark}
In order to obtain the estimate of time derivative for the annealed heat kernel, we need the following additional assumption.
\begin{assumption}\label{ass:time:derivative}
  There exists a non-random constant $\const[C]{lower:elliptic:time} \in [1, \infty)$ such that
  \begin{align}
    \prob\bigl[
      \const[C]{lower:elliptic:time}^{-1}
      \leq \omega_{t}(e) / \omega_{0}(e) \leq
      \const[C]{lower:elliptic:time}
    \bigr]
    \;=\;
    1,
    \qquad \forall\, e \in E_{d},\, t \in \mathbb{R}. 
  \end{align}
\end{assumption}
\begin{cor}\label{cor:time:derivative}
  Suppose that Assumption~\ref{ass:law}, \ref{ass:integrability} and \ref{ass:time:derivative} hold true.  Then, there exists $\const[C]{time:elliptic} < \infty$ such that, for any $y, y' \in \mathbb{Z}^{d}$ and $t \geq 1$,
  \begin{align}\label{eq:time:derivative}
    \Bigl| \frac{\partial}{\partial t}\, \bar{p}_{0, t}(y, y') \Bigr|
    \;\leq\;
    \const[C]{time:elliptic}\, t^{-d/2-1}.
  \end{align}
\end{cor}
\begin{remark}
  We note that for time-independent conductances, such time derivative estimates hold very generally and the proof simply uses the spectral resolution -- see, for instance \cite[Section 5.2]{GT12}.  However, for time-dependent conductances we do not know how to prove it without Assumption~\ref{ass:time:derivative}.  It is an interesting open problem how to establish time derivative estimates for more general time-dependent random conductances.
\end{remark}
As a next result we are interested in an \emph{averaged} or (\emph{annealed}) CLT (ACLT) as well as in the convergences of the corresponding annealed heat kernel and its first derivative.
\begin{prop}[Annealed CLT]\label{prop:annealed:CLT}
  Suppose that Assumption~\ref{ass:law} and \ref{ass:integrability} holds.  Then, for each $k \in \mathbb{N}$ and any $t_{1}, \ldots, t_{k}$ satisfying $0 \leq t_{1} < t_{2} < \cdots < t_{k} < \infty$, the distribution of $(X_{t_{1}}^{(n)}, \ldots, X_{t_{k}}^{(n)})$ under $\prob_{0, 0}^{*}$ converges weakly, as $n \to \infty$, to the distribution of $(W_{t_{1}}, \ldots, W_{t_{k}})$, where $(W_{t} : t \geq 0)$ is a Brownian motion on $\mathbb{R}^{d}$ with deterministic and non-degenerate covariance matrix $\Sigma^{2}$.
\end{prop}
\begin{remark}
  We note that a quenched FCLT has been obtained under Assumption~\ref{ass:law} and some additional moment bounds on both $\omega_{t}(e)$ and $\omega_{t}(e)^{-1}$ in \cite[Theorem 1.7]{ACDS18} when $d \geq 2$, and in \cite[Theorem~1.4]{DS16} when $d = 1$.  The latter result on the quenched invariance principle in $d = 1$ has been improved if Assumption~\ref{ass:law} holds true and $\mean\bigl[\omega_{t}(e)^{-1}\bigr]< \infty$, see \cite[Theorem~1.2]{Bi19}.  In case the dynamical conductances are bounded from above the additional assumption on $\omega_{t}(e)^{-1}$ has been recently relaxed, cf.\ \cite{BR18} and \cite{BP23}.  Clearly, quenched FCLT implies annealed FCLT.
\end{remark}
\begin{theorem}[Annealed local CLT]\label{thm:annealed:LCLT}
  Suppose that Assumption~\ref{ass:law} and \ref{ass:integrability}
  holds.  Further, for any $t > 0$, let $k_{t}$ be the density of the distribution $\mathcal{N}(0, t\Sigma^{2})$.  Then, 
  \begin{enumerate}[(i)]
    \item for any $t > 0$ and $K \subset \mathbb{R}^{d}$ compact,
    \begin{align}\label{eq:alclt}
      \lim_{n \to \infty} \sup_{y \in K}
      \bigl|
        n^{d} \bar{p}_{0, tn^{2}}(0, [y n]) - k_{t}^{\Sigma}(y)
      \bigr|
      \;=\;
      0.
    \end{align}
    If, additionally, Assumption~\ref{ass:time:derivative} holds, then for any $T_{0} > 0$ and $K \subset \mathbb{R}^{d}$ compact,
    \begin{align}\label{eq:alclt:time:constant}
      \lim_{n \to \infty} \sup_{y \in K} \sup_{t \geq T_{0}}
      \bigl|
        n^{d} \bar{p}_{0, tn^{2}}(0, [y n]) - k_{t}^{\Sigma}(y)
      \bigr|
      \;=\;
      0.
    \end{align}    
    \item for any $t > 0$, $K \subset \mathbb{R}^{d}$ compact and each unit vector $e_{i} \sim 0$, $i \in \{1, \ldots, d\}$, 
    \begin{align}\label{eq:alclt:gradient}
      \lim_{n \to \infty} \sup_{y \in K}
      \bigl|
        n^{d+1} \nabla_{e_{i}}^{2}\bar{p}_{0, tn^{2}}(0, [y n])
        - \bigl(\partial_{i}k_{t}^{\Sigma}\bigr)(y)
      \bigr|
      \;=\;
      0.
    \end{align}
  \end{enumerate}
\end{theorem}
Finally, as a further consequence of Theorem~\ref{thm:gradient:est}-(iii), the properly scaled second discrete derivatives of the heat kernel satisfies the following uniform upper bounds. For any $n \in \mathbb{N}$, $(0,x), (0,x') \in E_{d}$, $y, y' \in \mathbb{R}^{d}$, and $t > 0$ we have that
\begin{align}
  \bigl|
    n^{d+2} \nabla^{2}_{\!x} \nabla^{2}_{\!x'}\, \bar{p}_{0,tn^{2}}([yn], [y'n])
  \bigr|
  &\;\leq\;
  \const[C]{gradient:sum}\, (1 \vee t)^{-(d+2)/2},
  \\[1ex]
  \mean\Bigl[
    \Bigl|
      n^{d+2} \nabla^{1}_{-x}\nabla^{2}_{x'}\, p^{\omega}_{0,tn^{2}}([yn], [y'n])
    \Bigr|
  \Bigr]
  &\;\leq\;
  \const[C]{gradient:sum}\, (1 \vee t)^{-(d+2)/2}.
\end{align}
In particular, we have that the families
\begin{align*}
  \bigl\{
    n^{d+2} \nabla_{-x}^{1} \nabla_{x'}^{2}p_{0, tn^{2}}^{\omega}(0,[yn])
    : n \in \mathbb{N}
  \bigr\}
  \quad \text{and} \quad
  \bigl\{
    n^{d+2}\nabla_{x}^{2} \nabla_{x'}^{2}p_{0, tn^{2}}^{\omega}(0, [yn])
    : n \in \mathbb{N}
  \bigr\}
\end{align*}
are both tight under $\prob_{0, 0}^{*}$.
\medskip

It will be interesting to consider continuous versions of the results in this section.  Note that apriori it is not known whether $\nabla_{\!x}$ exists or not.
\begin{remark}
  Although under our assumptions, local CLT holds for the annealed heat kernel, the corresponding quenched result could fail.  Indeed, in \cite[Proposition~1.5]{DF20}, a random walk in layered conductances on $\mathbb{Z} \times \mathbb{Z}^{d-1}$ is considered with \emph{time-homogeneous} conductances
  \begin{align*}
    \omega(x, x \pm e_{1})
    \;=\;
    \begin{cases}
      Z(x_{2}), &i = 1,
      \\
      1, &i \in \{2, \ldots, d\},
    \end{cases}
    \qquad \forall\, x = (x_{1}, x_{2}) \in \mathbb{Z} \times \mathbb{Z}^{d-1},
  \end{align*}
  where $\{Z(x_{2}) : x_{2} \in \mathbb{Z}^{d-1}\}$ are i.i.d.\ random variables on some probability space $(\Omega, \mathcal{F}, \prob)$ taking values in $[0, \infty)$ such that
  \begin{align*}
    \prob\bigl[ Z(x_{2}) > L \bigr] 
    \;=\;
    L^{-\alpha}, 
    \qquad L \geq 1.
  \end{align*}
  If $d \geq 4$ and $\alpha \in (1, (d-1)/2)$ then it is show that, for any $R > 0$ and $t > 0$,
  \begin{align*}
    \lim_{n \to \infty} \inf_{|y| \leq R \sqrt{t}}
    \frac{n^{d} p_{0, t n^{2}}^{\omega}(0, [y n])}{k_{t}^{\Sigma}(y)}
    \;=\;
    0
    \qquad \prob\text{-a.s.}
  \end{align*}
  Further, notice that in this case $\mean\bigl[\omega(x, x+e_{1})\bigr] < \infty$ but $\mean\bigl[ \omega(x, x+e_{1})^{1/\alpha}\bigr] = \infty$, thus $\omega(x, x+e_{1}) \notin L^{p}(\prob)$ for $p > 2/(d-1)$, which is the sharp bound for quenched local CLT to hold, cf. \cite{BS22}.
\end{remark} 

\subsection{Structure of the paper}
In Section~\ref{subsec:entropy} we first establish estimates on the entropy which serves as the main building block to derive in Section~\ref{subsec:gradient:hk} and \ref{subsec:gradient:hk+green} space-derivatives on both the heat kernel and the Green kernel.  In Section~\ref{sec:alclt} we then prove Proposition~\ref{prop:annealed:CLT} and Theorem~\ref{thm:annealed:LCLT}.  In Appendix~\ref{appendix:semigroup} we verify properties of the transition semi-group of the time-inhomogeneous Markov process $X$, whereas Appendix~\ref{appendix:technical} just contains a technical lemma needed in the proofs.

\section{Estimates for the space-derivative of the heat kernel}
\label{section:gradient:estimates}

\subsection{Entropy estimates}
\label{subsec:entropy}
For any (sub-) probability measure $\mu$ on $\mathbb{Z}^{d}$ we define the entropy, $H(\mu)$, of $\mu$ by
\begin{align}
  H(\mu) \;\ldef\; \sum_{x \in \mathbb{Z}^{d}} \phi(\mu(x)),
\end{align}
where $\phi(0) = 0$ and $\phi(t) \ldef -t \ln t$ for any $t>0$.  Clearly, $H(\mu)$ is well-defined with values in $[0,\infty]$.  In particular, for any $t \geq 0$ and $\prob$-a.e.\ $\omega$, we set
\begin{align}
  H_{t}^{\omega}
  \;\ldef\;
  H^{\omega}(X_{t})
  \;\ldef
  H\bigl(\Prob_{0, 0}^{\omega} \circ X_{t}^{-1}\bigr)
  \;=\;
  H\bigl(\,p_{0,t}^{\omega}(0,\cdot)\, \bigr)
  .
\end{align}
Further, for any $0 < s < t < u$ and $\prob$-a.e.\ $\omega$, the conditional entropy is defined by
\begin{align}
  \label{def:conditional:entropy}
  \begin{split}
    H^{\omega}(X_{s} \mid X_{t})
    &\;\ldef\;
    \Mean_{0, 0}^{\omega}\bigl[
      H\bigl(\Prob_{0, 0}^{\omega}[X_{s} = \cdot \mid X_{t}]\bigr)
    \bigr]
    \\[.5ex]
    H^{\omega}(X_{s} \mid X_{t}, X_{u})
    &\;\ldef\;
    \Mean_{0, 0}^{\omega}\bigl[
      H\bigl(\Prob_{0, 0}^{\omega}[X_{s} = \cdot \mid X_{t}, X_{u}]\bigr)
    \bigr]
  \end{split}
\end{align}
Likewise, we define the entropy $H^{\tilde{\omega}}_{t} \ldef H^{\tilde{\omega}}(\tilde{X}_{t})$ and the conditional entropy $H^{\tilde{\omega}}(\tilde{X}_{s} \mid \tilde{X}_{t})$ and $H^{\tilde{\omega}}(\tilde{X}_{s} \mid \tilde{X}_{t}, \tilde{X}_{u})$. Notice that an alternative representation of $H^{\omega}(X_{s} \mid X_{t})$ is given in \eqref{eq:conditional:entropy}.

In order to address the question of $\prob$-a.s.\ finiteness of both $H_{t}^{\omega}$ and $H_{t}^{\tilde{\omega}}$ for any $t \geq 0$, in the next lemma we first prove a quenched result that holds true without Assumption~\ref{ass:law}.
\begin{lemma}\label{lem:Bass-Nash}
  There exists $\const[C]{entropy:ub} \equiv \const[C]{entropy:ub}(d) \in (0, \infty)$ such that for $\prob$-a.e.\ $\omega$ and every $t \geq 0$
  \begin{align}\label{eq:entropy:general:upper:bound}
    H_{t}^{\omega}
    \;\leq\;
    \const[C]{entropy:ub} + d \ln \bigl( 1 + \Mean_{0,0}^{\omega}\bigl[|X_{t}|\bigr] \bigr).
  \end{align}
\end{lemma}
\begin{proof}
  The estimate is based on a method that was originally introduced by Nash in \cite[p.~938]{Na58} and later improved by Bass in \cite{Ba02}.  Note that a discrete version of the Bass-Nash argument is discussed in \cite[Lemma~6.14]{Ba17} which implies immediately \eqref{eq:entropy:general:upper:bound}.  Nevertheless, for the convenience of the reader we will provide here a detailed proof.

  Let us fix some realisation $\omega \in \Omega$. First, for any $n \in \mathbb{N}$ and $t \geq 0$ set
  \begin{align*}
    M_{n}^{\omega}(t)
    \;\ldef\;
    \sum_{y \in B(0, n)} |y|\, p_{0, t}^{\omega}(0, y)
    \qquad \text{and} \qquad
    H_{n}^{\omega}(t)
    \;\ldef\;
    \sum_{y \in B(0, n)} \phi(p_{0, t}^{\omega}(0, y)).
  \end{align*}
  Second, notice that $\min_{s \geq 0} (\phi(s) + \lambda s) = - \me^{-\lambda-1}$ for any $\lambda \in \mathbb{R}$. Moreover, there exists $c \equiv c(d) < \infty$ such that $\sum_{y \in \mathbb{Z}^{d}} \me^{-a|y|_{2}} \leq c a^{-d}$, where $| \cdot |_{2}$ denotes the usual Euclidean norm on $\mathbb{Z}^{d}$. Thus, for any $a \in (0, 1]$ and $b \in [0, \infty)$ we get that, for any $t \geq 0$,
  \begin{align*}
    -H_{n}^{\omega}(t) + a M_{n}^{\omega}(t) + b
    &\;\geq\;
    \sum_{y \in B(0, n)} \bigl(
      \phi(p_{0, t}^{\omega}(0, y)) + (a |y| + b)\, p_{0, t}^{\omega}(0, y)
    \bigr)
    \\
    &\;\geq\;
    - \me^{-b-1} \sum_{y \in \mathbb{Z}^{d}} \me^{-a |y|}
    \\
    &\;\geq\;
    -c\, \me^{-b-1} a^{-d}.
  \end{align*}
  By setting $a \ldef (1 + M_{n}(t))^{-1}$ and $b \ldef -d \ln a = d \ln (1 + M_{n}(t))$, we obtain from Lemma~\ref{lem:non-explosion}-(i) that, $\prob$-a.s.\ $a \in (0, 1]$ and $b \in [0, \infty)$. Thus, for $\prob$-a.e.\ $\omega$ and any $t \geq 0$
  \begin{align*}
    -H_{n}^{\omega}(t) + \frac{M_{n}^{\omega}(t)}{1 + M_{n}^{\omega}(t)}
    + d \ln\bigl(1 + M_{n}^{\omega}(t) \bigr)
    \;\geq\;
    -c\, \me^{-1}
  \end{align*}
  Since $M_{n}^{\omega}(t) / (1 + M_{n}^{\omega}(t)) \leq 1$, by setting $\const[C]{entropy:ub} \ldef 1 + c\,\me^{-1}$ and rearranging the terms in the above expression, we get that $H_{n}^{\omega}(t) \leq \const[C]{entropy:ub} + d \ln(1 + M_{n}^{\omega}(t))$ for $\prob$-a.e.\ $\omega$ and every $t \geq 0$. The final assertion now follows by monotone convergence.
\end{proof}
Thus, together with Proposition~\ref{prop:mean:displacement}, it follows that
\begin{align}
  \bar{H}_{t}
  \;\ldef\;
  \mean\bigl[H_{t}^{\omega}\bigr]
  \;\leq\;
  \const[C]{entropy:ub} + d \ln \bigl(1 + \const[C]{mean:displacement} \sqrt{t} \bigr)
  \qquad \Longrightarrow \qquad
  H_{t}^{\omega} \;<\; \infty
  \quad \prob\text{-a.s.}
\end{align}
In particular, as a further consequence of \eqref{eq:shift:space-time}, we have that $\mean[H_{t}^{\tilde{\omega}}] = \mean[H_{t}^{\omega}]$ and, hence, $H_{t}^{\tilde{\omega}}< \infty$ $\prob$-a.s.

We now give fundamental properties of $H_{t}^{\omega}$, $H_{t}^{\tilde{\omega}}$ and $\bar{H}_{t}$.
\begin{lemma}\label{lem:entropy:convex}
  Suppose Assumption~\ref{ass:law} holds true.  Then,
  \begin{enumerate}[(i)]
  \item
    for each $s, t \geq 0$ and $\prob$-a.e.\ $\omega$, it holds that $H_{t}^{\omega} \leq H_{t+s}^{\omega}$ and $H_{t}^{\tilde{\omega}} \leq H_{t+s}^{\tilde{\omega}}$.
    
  \item
    for each fixed $s>0$, $t \mapsto \bar{H}_{t+s}-\bar{H}_{t}$ is decreasing for $t \geq 0$. In other words, the function $t \mapsto \bar{H}_{t}$ is concave.
  \end{enumerate}
\end{lemma}
\begin{proof}
  (i) By using the Chapman-Kolmogorov equation, we can rewrite the entropy $H_{s+t}^{\omega}$, for $\prob$-a.e.\ $\omega$ and any $s, t \geq 0$, as
  \begin{align*}
    H_{t+s}^{\omega}
    \;=\;
    \sum_{x \in \mathbb{Z}^{d}} \phi\bigl( p_{0,t+s}^{\omega}(0,x) \bigr)
    \;=\;
    \sum_{x \in \mathbb{Z}^{d}}
    \phi\biggl(
      \sum_{y \in \mathbb{Z}^{d}} p_{0,t}^{\omega}(0,y)\, p_{t,s+t}^{\omega}(y,x)
    \biggr).
  \end{align*}
  Since, by Lemma~\ref{lem:non-explosion}-(i), $\sum_{y \in \mathbb{Z}^{d}} p_{t,s+t}^{\omega}(y,x) = 1$ for $\prob$-a.e.\ $\omega$, an application of Jensen's inequality yields
  \begin{align*}
    H_{t+s}^{\omega}
    \;\geq\;
    \sum_{x \in \mathbb{Z}^{d}}
    \sum_{y \in \mathbb{Z}^{d}} \phi
    \bigl( p_{0,t}^{\omega}(0,y) \bigr)\, p_{t,s+t}^{\omega}(y,x)
    \;=\;
    \sum_{y \in \mathbb{Z}^{d}} \phi \bigl( p_{0,t}^{\omega}(0,y) \bigr)
    \;=\;
    H_{t}^{\omega}.
  \end{align*}
  By the same argument the assertion also follows for $H_{t}^{\tilde{\omega}}$.
  
  (ii) The proof relies on the argument given in \cite[Corollary~10]{BD-CKY15}.  First, notice that, for $\prob$-a.e.\ $\omega$ and any $t \geq 0$ and $s > 0$, the conditional entropy $H^{\omega}$ can also be expressed as
  \begin{align}
    \label{eq:conditional:entropy}
    H^{\omega}\bigl(X_{s} \,|\, X_{s+t}\bigr)
    &\;=\;
    \sum_{y \in \mathbb{Z}^{d}} \Prob_{0,0}^{\omega}\bigl[X_{s+t} = y\bigr]
    \sum_{x \in \mathbb{Z}^{d}}
    \phi\bigl(\Prob_{0,0}^{\omega}[X_{s} = x \mid X_{s+t} = y]\bigr)
    \nonumber\\[.5ex]
    &\;=\;
    \sum_{x,y \in \mathbb{Z}^{d}}  
    \phi\bigl(\Prob_{0,0}^{\omega}\bigl[X_{s} = x, X_{s+t} = y\bigr]\bigr)
    \,-\,
    \sum_{y \in \mathbb{Z}^{d}} 
    \phi\bigl(\Prob_{0,0}^{\omega}\bigl[X_{s+t} = y\bigr]\bigr).
  \end{align}
  By means of the Markov property, the first term can be computed further as 
  \begin{align*}
    &\sum_{x,y \in \mathbb{Z}^{d}}
    \phi\bigl(
      \Prob_{0,0}^{\omega}\bigl[X_{s} = x\bigr]
      \Prob_{s,x}^{\omega}\bigl[X_{s+t} = y\bigr]
    \bigr)
    \\[.5ex]
    &\mspace{36mu}=\;
    \sum_{x \in \mathbb{Z}^{d}} \phi\bigl(\Prob_{0,0}^{\omega}\bigl[X_{s} = x\bigr]\bigr)
    \,+\,
    \sum_{x,y \in \mathbb{Z}^{d}} \Prob_{0,0}^{\omega}\bigl[X_{s} = x\bigr]
    \phi\bigl(\Prob_{s,x}^{\omega}\bigl[X_{s+t} = y\bigr]\bigr)
    \\[.5ex]
    &\mspace{36mu}=\;
    H_{s}^{\omega}
    \,+\, 
    \Mean_{0,0}^{\omega}\bigl[H\bigl(p^\omega_{s,s+t}(X_{s},\,\cdot\,)\bigr)\bigr].
  \end{align*}
  By taking the expectation with respect to $\prob$ and using the fact that
  \begin{align}\label{eq:noeoj2g}
    \mean\bigl[
      \Mean_{0,0}\bigl[H\bigl(p_{s,s+t}^{\omega}(X_{s}, \cdot )\bigr) \bigr]
    \bigr]
    \overset{\eqref{eq:shift:hk}}{\;=\;}
    \mean\bigl[\Mean_{0,0}\bigl[ H_{t}^{\omega} \circ \tau_{s, X_{s}}\bigr] \bigr]
    \;=\;
    \bar{H}_{t}, 
  \end{align}
  we therefore obtain
  \begin{align*}
    \mean\bigl[H^{\omega}(X_{s} \,|\, X_{s+t})\bigr]
    \;=\;
    \bar{H}_{s} \,+\, \bar{H}_{t} \,-\, \bar{H}_{s+t}.
  \end{align*}
  Further, we claim that
  \begin{align}\label{eq:enosgb}
    H^{\omega}\bigl(X_{s} \,|\, X_{s+t}\bigr)
    \;=\;
    H^{\omega}\bigl(X_{s} \,|\, X_{s+t}, X_{s+t'}\bigr)
    \;\leq\;
    H^{\omega}\bigl(X_{s} \,|\, X_{s+t'}\bigr)
    \qquad \text{for } t \leq t'. 
  \end{align}
  While the equality is trivial, the inequality relies on the general fact, cf.\ \cite[Equation~(2.60) and (2.92)]{CT06}, that for any discrete random variables $X$, $Y$ and $Z$ (defined on the same probability space) $H(X \,|\, Y,Z) \leq H(X \,|\, Z)$. Indeed, by Jensen's inequality, we have that 
  \begin{align*}
    H(X \mid Y,Z)
    &\;=\;
    \sum_{x,y,z} \Prob[X=x, Y=y, Z=z] \log
    \biggl(\frac{1}{\Prob[X=x \mid Y=y, Z=z]}\biggr)
    \\
    &\;\leq\;
    \sum_{x,z} \Prob[X=x, Z=z] \log
    \biggl(
      \sum_{y} \frac{\Prob[Y=y \mid X=x, Z=z]}{\Prob[X=x \mid Y=y, Z=z]}
    \biggr)
    \\
    &\;=\;
    \sum_{x,z} \Prob[X=x, Z=z] \log
    \biggl(
      \sum_{y} \frac{\Prob[Y=y, Z=z]}{\Prob[X=x, Z=z]}
    \biggr)
    \;=\;  
    H(X \,|\, Z),
  \end{align*}
  where the summation is only over all those elements that are contained in the distribution of $\Prob \circ\, (X,Y,Z)^{-1}$ and $\Prob \circ\, (X,Z)^{-1}$, respectively.  Thus, by combining \eqref{eq:noeoj2g} with \eqref{eq:enosgb}, we obtain, for any $t \leq t'$,
  \begin{align*}
    \bar{H}_{s} \,+\, \bar{H}_{t} \,-\, \bar{H}_{s+t}
    \;\leq\;
    \bar{H}_{s} \,+\, \bar{H}_{t'} \,-\, \bar{H}_{s+t'},
  \end{align*}
  which completes the proof.
\end{proof}
%
%
\begin{lemma}\label{lem:Nashclas}
  Suppose Assumption~\ref{ass:law} and \ref{ass:integrability} hold true.  Then, there exists a non-random constant $\const[C]{heatkernel:ub} \in (0, \infty)$ such that, for all $s < t$,
  \begin{align}\label{eq:nash}
    \max_{x,y \in \mathbb{Z}^{d}}\, p_{s,t}^{\omega}(x,y)
    \;\leq\;
    \const[C]{heatkernel:ub} \bigl(1 \vee (t-s)\bigr)^{-d/2}.
  \end{align}
  Moreover, the same estimate holds if $p_{s,t}^{\omega}(x,y)$ is replaced by $p_{s,t}^{\tilde{\omega}}(x,y)$.
\end{lemma}
\begin{proof}
  Let $B_{n} \ldef B(0,n)$ and set $\ell_{B_{n}}(\mathbb{Z}^{d}) \ldef \{f\colon \mathbb{Z}^{d}\to \mathbb{R} \,|\, f(x)=0$ for $x \notin B_{n}\}$.  Further, for any $f, g\colon \mathbb{Z}^{d} \to \mathbb{R}$ and $t \geq 0$ define the time-dependent Dirichlet form
  \begin{align*}
    \mathcal{E}_{t}(f, g)
    \;\ldef\;
    \sum_{(x,y) \in E_{d}}\mspace{-9mu} \omega_{t}(x,y)\, \bigl(f(y) - f(x) \bigr)\bigl(g(x) - g(y)\bigr),
  \end{align*}
  whenever this sum converges absolutely.  Then, for any $f \in \ell_{B_{n}}(\mathbb{Z}^{d})$,
  \begin{align*}
    \mathcal{E}_{t}^{\omega}(f,f)
    \;=\;
    \sum_{(x,y) \in E_{d}}\mspace{-9mu} \omega_{t}(x,y)\, \bigl(f(y) - f(x) \bigr)^{2}
    \;=\;
    \scpr{f}{-\mathcal{L}_{t}^{\omega} f}{\mathbb{Z}^{d}},
  \end{align*}
  where the second equality holds true due to the fact that the sums are each over finite sets.  Moreover, it is well-known, see e.g.\ \cite[Lemma~3.13]{Ba17}, that on the Euclidean lattice $(\mathbb{Z}^{d}, E_{d})$ the Nash inequality holds true, that is, for any $f \in \ell^{1}(\mathbb{Z}^{d}) \cap \ell^{2}(\mathbb{Z}^{d})$
  \begin{align}
    \label{eq:Nash:ineq:Z^{d}}
    \|f\|_{\ell^{2}(\mathbb{Z}^{d})}^{2+4/d}
    \;\leq\;
    C_{\mathrm{Nash}} \sum_{(x,y) \in E_{d}} \bigl(f(x) - f(y)\bigr)^{2}\, \|f\|_{\ell^{1}(\mathbb{Z}^{d})}^{4/d}.
  \end{align}
  On the other hand, in view of Assumption~\ref{ass:law}-(ii) and~(iv) combined with Assumption~\ref{ass:integrability}, an application of Fubini's theorem yields
  \begin{align*}
    0
    \;=\;
    \int_{0}^{t}
    \mean\bigl[
      \mathbbm{1}_{\omega_{s}(e) < C_{2}}
    \bigr]\,
    \md s
    \;=\;
    \mean\biggl[
      \int_{0}^{t} \mathbbm{1}_{\omega_{s}(e) < C_{2}}\, \md s
    \biggr],
    \qquad \forall\, t > 0,\, e \in E_{d}.
  \end{align*}
  Hence, $\prob$-a.s.\ $\omega_{s}(e) \geq \const[C]{lower:elliptic}$ for a.e.\ $s \in [0, t]$ and any $e \in  E_{d}$. In particular, it follows from \eqref{eq:Nash:ineq:Z^{d}}, that $\prob$-a.s., for any $n \in \mathbb{N}$ and a.e.\ $s \in [0, t]$, 
  \begin{align}
    \|f\|_{\ell^{2}(\mathbb{Z}^{d})}^{2+4/d}
    \;\leq\;
    \frac{C_{\mathrm{Nash}}}{\const[C]{lower:elliptic}}\,
    \mathcal{E}_{s}^{\omega}(f,f)\, \|f\|_{\ell^1(\mathbb{Z}^{d})}^{4/d},
    \qquad \forall\, f \in \ell_{B_{n}}(\mathbb{Z}^{d}).
  \end{align}
  Consider now the time-inhomogeneous random walk, $X$, killed upon exiting the ball $B_{n}$. Its heat kernel can be written as
  \begin{align*}
    p_{s, t}^{\omega, B_{n}}(x, y)
    \;\ldef\;
    \Prob_{s,x}^{\omega}\bigl[X_{t} = y, T_{n} > t\bigr],
  \end{align*}
  where $T_{n} \ldef \inf\{t \geq s: X_{t} \not\in B(0,n)\}$ denotes the exit time of the ball $B(0,n)$. Using \eqref{eq:beqweak} for $P_{s,t}^{\omega, B_{n}}f$ which is finitely supported, we have for $\prob$-a.e.\ $\omega$, a.e.\ $s \in (0,t)$ and $f \in \ell_{B_{n}}(\mathbb{Z}^{d})$, 
  \begin{align*}
    \frac{\md}{\md s}\, \|P_{s,t}^{\omega, B_{n}}f\|_{\ell^{2}(\mathbb{Z}^{d})}^{2}
    \;=\;
    2
    \langle
    -\mathcal{L}_{s}^{\omega} P_{s,t}^{\omega, B_{n}} f, P_{s,t}^{\omega, B_{n}} f
    \rangle_{\ell^{2}(\mathbb{Z}^{d})}
    \;=\;
    2\, \mathcal{E}_{s}^{\omega}(P_{s,t}^{\omega, B_{n}}f, P_{s,t}^{\omega, B_{n}} f). 
  \end{align*}
  We note that the last equality (the Gauss-Green formula) uses exchanges of the order of sums; in this case the sums are finite, so the equality holds without any problem.  For any $f \in \ell_{B_{n}}(\mathbb{Z}^{d})$ with $\|f\|_{\ell^1(\mathbb{Z}^{d})} = 1$ and $t > 0$ we set $f_{s} \ldef P_{s,t}^{\omega,B_{n}}f$.  Then, we obtain, for $\prob$-a.e.\ $\omega$ and a.e.\ $s \in (0,t)$,
  \begin{align*}
    u'(s)
    \;\ldef\;
    \frac{\md}{\md s} \|f_{s}\|_{\ell^{2}(\mathbb{Z}^{d})}^{2}
    \;=\;
    2\, \mathcal{E}_{s}^{\omega}(f_{s},f_{s})
    \;\geq\;
    \frac{2 \const[C]{lower:elliptic}}{C_{\mathrm{Nash}}}\, \|f_{s}\|_{\ell^{2}(\mathbb{Z}^{d})}^{2+4/d}
    \;\rdef\;
    c_{1}\, u(s)^{1+2/d}.
  \end{align*}
  Set $v(s) \ldef u(s)^{-2/d}$, then we obtain $v'(s) \leq -2c_{1}/d$. Since $\lim_{s \uparrow t} v(s) = \|f\|_{\ell^{2}(\mathbb{Z}^{d})}^{-4/d} > 0$, it follows
  \begin{align*}
    -v(s)
    \;\leq\;
    \int_{s}^{t} v'(r)\, \md r
    \;\leq\;
    -\frac{2c_{1}}{d}\, (t-s)
    \qquad\Longleftrightarrow\qquad
    u(s) \;\leq\; c_{2}\, (t-s)^{-d/2}.
  \end{align*}
  By introducing the notation to denote the norm of an operator $A$ from $\ell^{p}(\mathbb{Z}^{d})$ to $\ell^{q}(\mathbb{Z}^{d})$ for some $p, q \in [1, \infty]$ by
  \begin{align*}
    \|A\|_{p \to q}
    \;\ldef\;
    \sup\bigl\{
      \|A f\|_{\ell^{q}(\mathbb{Z}^{d})} : \|f\|_{\ell^{p}(\mathbb{Z}^{d})} \leq 1
    \bigr\},
  \end{align*}
  the above estimate yields $\|P_{s,t}^{\omega, B_{n}}\|_{1 \to 2} \leq c_{3} (t-s)^{-d/4}$.  On the other hand, by the differential forward equation in weak sense, we find that, for $\prob$-a.e.\ $\omega$ and a.e.\ $t \in (s, \infty)$,
  \begin{align*}
    \frac{\md}{\md t} \scpr{f}{\bigl(P_{s,t}^{\omega, B_{n}}\bigr)^{\!*}g}{\mathbb{Z}^{d}}
    &\;=\;
    \frac{\md}{\md t} \scpr{P_{s,t}^{\omega, B_{n}}f}{g}{\mathbb{Z}^{d}}
    \;=\;
    \scpr{P_{s,t}^{\omega} (\mathcal{L}_{t}^{\omega} f)}{g}{\mathbb{Z}^{d}}
    \\[.5ex]
    &\;=\;
    \scpr{\mathcal{L}_{t}^{\omega}f}{\bigl(P_{s,t}^{\omega, B_{n}}\bigr)^{\!*} g}{\mathbb{Z}^{d}}
  \end{align*}
  for any $f,g \in \ell_{B_{n}}(\mathbb{Z}^{d})$.  In particular, for $\prob$-a.e.\ $\omega$ and a.e.\ $t \in (s,\infty)$, it holds that
  \begin{align*}
    \frac{\md}{\md t}\, \|\bigl(P_{s,t}^{\omega}\bigr)^{\!*} f\|_{\ell^{2}(\mathbb{Z}^{d})}^{2}
    &\;=\;
    2\,
    \scpr{\mathcal{L}_{t}^{\omega} \bigl(P_{s,t}^{\omega, B_{n}}\bigr)^{\!*} f}
    {\bigl(P_{s,t}^{\omega, B_{n}}\bigr)^{\!*} f}{\mathbb{Z}^{d}}
    \\[.5ex]
    &\;=\;
    -2\,
    \mathcal{E}_{t}^{\omega}\bigl(
      \bigl(P_{s,t}^{\omega, B_{n}}\bigr)^{\!*} f, 
      \bigl(P_{s,t}^{\omega, B_{n}}\bigr)^{\!*} f
    \bigr),
    \qquad \forall\, f \in \ell_{B_{n}}(\mathbb{Z}^{d}).
  \end{align*}
  Thus, by a computation similar to the one that we did above, we obtain that $\|\bigl(P_{s,t}^{\omega, B_{n}}\bigr)^{\!*}\|_{1 \to 2} \leq c_{3} (t-s)^{-d/4}$.  Recall that, by the Chapman-Kolmogorov equation it holds that $P_{s,t}^{\omega, B_{n}} = P_{s,(s+t)/2}^{\omega, B_{n}} \circ P_{(s+t)/2,t}^{\omega, B_{n}}$ and $\bigl(P_{s,t}^{\omega, B_{n}}\bigr)^{\!*} = \bigl(P_{(s+t)/2,t}^{\omega, B_{n}}\bigr)^{\!*} \circ \bigl(P_{s,(s+t)/2}^{\omega, B_{n}}\bigr)^{\!*}$.  Thus, by using the Cauchy-Schwarz equation, we find that
  \begin{align*}
    \left.
      \begin{array}{c}
        \|P_{s,t}^{\omega, B_{n}}\|_{1 \to \infty} \\[1ex]
        \|\bigl(P_{s,t}^{\omega, B_{n}}\bigr)^{\!*}\|_{1 \to \infty}
      \end{array}
    \right\}
    \;\leq\;
    \|\bigl(P_{s,(s+t)/2}^{\omega, B_{n}}\bigr)^{\!*}\|_{1 \to 2}\,
    \|P_{(s+t)/2,t}^{\omega, B_{n}}\|_{1 \to 2}.
  \end{align*}
  Moreover, by Proposition~\ref{prop:mean:displacement},
  \begin{align*}
    \Prob_{s,x}^{\omega}\Bigl[ \lim_{n \to \infty} T_{n}=\infty \Bigr] \;=\; 1.
  \end{align*}
  Indeed, if not we have $\Prob_{s,x}^{\omega}\bigl[\lim_{n \to \infty} T_{n}\leq M \bigr] > \varepsilon$ for some large $M > 0$ and small $\varepsilon > 0$, which contradicts 
  Proposition~\ref{prop:mean:displacement}.  Hence, by the monotone convergence theorem, we have, for each $x,y \in \mathbb{Z}^{d}$,
  \begin{align*}
    p_{s,t}^{\omega}(x,y)
    \;=\;
    \Prob_{s,x}^{\omega}\bigl[X_{t} = y\bigr]
    \;=\;
    \lim_{n \to \infty}\Prob_{s,x}^{\omega}\bigl[X_{t} = y, T_{n} \geq t\bigr]
    \;=\;
    \lim_{n \to \infty}p_{s,t}^{\omega,  B_{n}}(x,y).
  \end{align*}
  In particular, by the dominated convergence theorem,
  \begin{align*}
    \bigl|\scpr{f}{P_{s,t}^{\omega}g}{\mathbb{Z}^{d}}\bigr|
    \;=\;
    \lim_{n \to \infty} \bigl|\scpr{f}{P_{s,t}^{\omega, B_{n}}g}{\mathbb{Z}^{d}}\bigr|
    \qquad \forall\, f,g \in \ell^1(\mathbb{Z}^{d}).
  \end{align*}
  Thus, we conclude
  \begin{align*}
    \|\bigl(P_{s,t}^{\omega}\bigr)^{*}\|_{1 \to \infty}
    \;=\;
    \|P_{s,t}^{\omega}\|_{1 \to \infty}
    \;\leq\;
    \limsup_{n \to \infty} \|P_{s,t}^{\omega,  B_{n}}\|_{1 \to \infty}
    \;\leq\;
    \frac{c_{3}^{2}}{(t-s)^{d/2}},
  \end{align*}
  and the assertion follows.
\end{proof}
\begin{prop}\label{prop:entropy:estimate}
  Suppose that Assumption~\ref{ass:law} and \ref{ass:integrability} hold.  Then, for any $t_{0} > 0$, there exists $\const[C]{entropy:difference} \equiv \const[C]{entropy:difference}(t_{0}) > 0$ such that for all $s \geq 0$ and $t \geq t_{0}$, it holds that
  \begin{align}\label{eq:entropy:difference}
    \bar{H}_{t+s}-\bar{H}_{t}
    \;\leq\;
    \const[C]{entropy:difference}\, \frac{s}{t}.
  \end{align}
\end{prop}
\begin{proof}
  First, by Lemma~\ref{lem:Nashclas} it holds that $\bar{H}_{t} \geq \frac{d}{2} \ln t - c$.  Also, by Lemma~\ref{lem:Bass-Nash} and Proposition~\ref{prop:mean:displacement}, using Jensen's inequality we have $\bar{H}_{t} \leq \frac{d}{2} \ln (1+t) + c$ for some $c < \infty$.  These two estimates imply that, for any $s,t > 0$, 
  \begin{align}\label{eq:entropy:difference:apriori}
    \bar{H}_{t+s} - \bar{H}_{t}
    \;\leq\;
    \frac{d}{2}\, \ln\Bigl(1 + \frac{1+s}{t}\Bigr) + c.
  \end{align}
  Let us now fix some $t_{0} > 0$.  Then, for any $t_{0} \leq t \leq s$, we immediately deduce from \eqref{eq:entropy:difference:apriori} that
  \begin{align*}
    \bar{H}_{t+s} - \bar{H}_{t}
    \;\leq\;
    \frac{s}{t}\,
    \biggl(
      \frac{t}{s} \cdot \frac{d}{2}
      \ln\Bigl(1 + \frac{1}{t_{0}} + \frac{s}{t}\Bigr) + c
    \biggr)
    \;\leq\;
    \frac{s}{t}\, \biggl(\frac{d}{2} \Bigl(1+\frac{1}{t_{0}}\Bigr) + c \biggr).
  \end{align*}
  So it remains to prove the estimate in case $0 \leq s < t$ and $t \ge t_{0}$.  For this purpose, choose $k \in \mathbb{N}$ such that $s (k-1) \leq t/2 < ks$.  Then, in view of Lemma~\ref{lem:entropy:convex}-(ii),
  \begin{align*}
    \bar{H}_{t+s} - \bar{H}_{t}
    \;\leq\;
    \frac{1}{k}\,
    \sum_{i=1}^{k} \bigl( \bar{H}_{t/2 + i s} - \bar{H}_{t/2 + (i-1)s} \bigr)
    \;=\;
    \frac{1}{k}\, \bigl(\bar{H}_{t/2 + k s} - \bar{H}_{t/2}\bigr).
  \end{align*}
  But,
  \begin{align*}
    \frac{1}{k}\, \bigl(\bar{H}_{t/2 + k s} - \bar{H}_{t/2}\bigr)
    \overset{\!\!\!\eqref{eq:entropy:difference:apriori}\!\!\!}{\;\leq\;}
    \frac{2 s}{t}\,
    \biggl(\frac{d}{2} \ln \Bigl( 1 + \frac{2(1+k s)}{t} \Bigr) + c \biggr)
    \;\leq\;
    \frac{2s}{t}\, \biggl(\frac{d}{2}\Bigl(3+\frac{2}{t_{0}}\Bigr) + c \biggr).
  \end{align*}
  By choosing $\const[C]{entropy:difference}(t_{0}) \ldef d(3+2/t_{0}) + 2c$, the assertion follows.
\end{proof}
For any probability measure $\mu, \nu$ on $\mathbb{Z}^{d}$ we further define the $\Delta$-distance,
$\Delta(\mu, \nu)$, between $\mu$ and $\nu$ by
\begin{align}
  \Delta(\mu,\nu)^{2}
  \;\ldef\;
  \sum_{\substack{x \in \mathbb{Z}^{d} \\ \mu(x) + \nu(x) > 0}}\mspace{-6mu}
  \frac{(\mu(x)-\nu(x))^{2}}{\mu(x)+\nu(x)}.
\end{align}
We also write
\begin{align*}
  \Delta_{s,t}^{\!\omega}(0,x)
  &\;\ldef\;
  \Delta
  \bigl(\,
    p_{0,s+t}^{\omega}(0,\cdot),\, 
    p_{s,s+t}^{\omega}(x,\cdot)\,
  \bigr)
  \qquad \forall\, x \in \mathbb{Z}^{d},\, s,t \geq 0.
\end{align*}
Likewise, we define $\Delta_{s,t}^{\tilde{\omega}}$.
\begin{prop}\label{BD-CKY15:theorem8}
  Suppose that Assumption~\ref{ass:law} holds true.
  \begin{enumerate}[(i)]
  \item
    Then, for every $t > 0$, 
    \begin{align}\label{eq:BD-CKY15:theorem8}
      \mean\bigl[\Mean_{0,0}^{\omega}\bigl[\Delta_{s,t}^{\omega}(0, X_{s})^{2} \bigr] \bigr]
      \;\leq\;
      2\, \bigl( \bar{H}_{s+t} \,-\, \bar{H}_{t} \bigr).
    \end{align}
  \item
    If Assumption~\ref{ass:integrability} is additionally satisfied, then there exists $\const[C]{mean:entropy:difference} \equiv \const[C]{mean:entropy:difference}(t_{0}) > 0$ such that the following holds for all $t \geq t_{0}$,
    \begin{align}\label{eq:mean:entropy:est}
      \mean\bigl[ \Mean_{0,0}\bigl[ \Delta_{s,t}^{\omega}(0, X_{s})^{2} \bigr] \bigr]
      \;\leq\; 
      2 \bigl(
        \bar{H}_{s + t} \,-\, \bar{H}_{t}
      \bigr)
      \;\leq\;
      \const[C]{mean:entropy:difference}\, \frac{s}{t}.
    \end{align}
  \end{enumerate}
  Further, the same estimates hold if $\Mean_{0,0}^{\omega}\bigl[\Delta_{s,t}^{\omega}(0, X_{s})^{2}\bigr]$ is replaced by $\Mean_{0,0}^{\tilde{\omega}}\bigl[\Delta_{s,t}^{\tilde{\omega}}(0, \tilde{X}_{s})^{2}\bigr]$.
\end{prop}
\begin{proof}
  (i) The proof is based on arguments originally given in \cite[Lemma~7]{BD-CKY15}.  First, notice that for any $a > 0$
  \begin{align*}
    a \ln a
    \;\geq\;
    \frac{1}{2}\, \frac{(a-1)^{2}}{a+1} + a-1.
  \end{align*}
  By choosing $a = p_{s,s+t}^{\omega}(X_{s},y) / p_{0,s+t}^{\omega}(0,y)$, we obtain
  \begin{align*}
    \Mean_{0,0}^{\omega}\bigl[\Delta_{s,t}^{\omega}(0, X_{s})^{2}\bigr]
    &\leq\;
    2 \sum_{y \in \mathbb{Z}^{d}}
    \Mean_{0,0}^{\omega}\Bigl[
     p_{s,s+t}^{\omega}(X_{s},y) 
      \ln\bigl(
        p_{s,s+t}^{\omega}(X_{s},y)/p_{0,s+t}^{\omega}(0,y)    
      \bigr)
    \Bigr]
    \\
    &\mspace{72mu}-\,
    2 \sum_{y \in \mathbb{Z}^{d}}
    \Bigl(
      \Mean_{0,0}^{\omega}\bigl[ p_{0,s+t}^{\omega}(0,y) \bigr]
      \,-\,
      \Mean_{0,0}^{\omega}\bigl[ p_{s,s+t}^{\omega}(X_{s},y) \bigr]
    \Bigr)
    \\[.5ex]
    &=\;
    2\, \sum_{y \in \mathbb{Z}^{d}}
    \Bigl(  
      \phi\bigl( p_{0,s+t}^{\omega}(0,y) \bigr)
      \,-\,
      \Mean_{0,0}^{\omega}\Bigl[
        \phi\bigl( p_{s,s+t}^{\omega}(X_{s},y) \bigr)
      \Bigr]
    \Bigr),
  \end{align*}
  where we used that the Chapman-Kolmogorov equation implies 
  \begin{align*}
    &\Mean_{0,0}^{\omega}\Bigl[ 
      p_{s,s+t}^{\omega}(X_{s},y) \,\ln p_{0,s+t}^{\omega}(0,y)
    \Bigr]
    \\[.5ex]
    &\mspace{36mu}=\;
    \sum_{x' \in \mathbb{Z}^{d}} p_{0,s}^{\omega}(0,x')\, p_{s,s+t}^{\omega}(x',y)\,
    \ln p_{0,s+t}^{\omega}(0,y)
    \;=\;
    - \phi\bigl(p_{0,s+t}^{\omega}(0,y)\bigr).
  \end{align*}
  By taking expectation with respect to ${\prob}$ and using \eqref{eq:noeoj2g}, the assertion follows.

  (ii) This follows from Proposition~\ref{BD-CKY15:theorem8}-(i) combined with \eqref{eq:entropy:difference}.  
\end{proof}

\subsection{Discrete space-derivatives of the heat kernel}
\label{subsec:gradient:hk}
The main objective in this subsection is to prove Theorem~\ref{thm:gradient:est}.
\begin{proof}[Proof of Theorem~\ref{thm:gradient:est}]
  (i) In view of \eqref{eq:shift:hk} and Assumption~\ref{ass:law}-(i), it suffices to prove \eqref{eq:gradient:on-diag} for $y=0$.  The proof comprises two steps.

  \textit{Step 1.}
  We first show that, for any $x, y \in \mathbb{Z}^{d}$, $t \geq 1$ and $s > 0$,
  \begin{align}\label{eq:gradient:claim}
    \mean\Bigl[
      p_{0,s}^{\omega}(0,x)\,
      \bigl| p_{s,s+t}^{\omega}(x,y) - p_{0,s+t}^{\omega}(0,y) \bigr|^{2}
    \Bigr]
    \;\leq\;
    c\, \frac{s}{t^{d+1}}.
  \end{align}
  By using the Chapman-Kolmogorov equation and the Cauchy-Schwarz inequality we find that
  \begin{align}
    &\bigl|
      p^{\omega}_{s,s+t}(x, y) \,-\, p_{0, s+t}^{\omega}(0, y) 
    \bigr| 
    \nonumber\\
    &\mspace{36mu}\leq\;
    \sum_{z \in \mathbb{Z}^{d}}
    \bigl|
      p^{\omega}_{s, s+t/2}(x, z) \,-\, p_{0, s+t/2}^{\omega}(0,z) 
    \bigr|\,
    p_{s+t/2,s+t}^{\omega}(z,y)
    \label{eq:gradient:chapman}\\[.5ex]
    &\mspace{36mu}\leq\;
    \Delta\bigl(
      p_{0,s+t/2}^{\omega}(0, \cdot),\, p_{s,s+t/2}^{\omega}(x, \cdot)
    \bigr)
    \nonumber\\
    &\mspace{72mu}\times\;
    \Biggl(
      \sum_{z \in \mathbb{Z}^{d}} 
      \bigl(p_{0,s+t/2}^{\omega}(0,z) + p_{s,s+t/2}^{\omega}(x,z)\bigr)\,
      p_{s+t/2,s+t}^{\omega}(z,y)^{2}
    \Biggr)^{\!1/2}
    \nonumber\\[.5ex]
    &\mspace{36mu}\leq\;
    \Delta\bigl(
      p_{0,s+t/2}^{\omega}(0, \cdot),\, p_{s,s+t/2}^{\omega}(x, \cdot)
    \bigr)
    \nonumber\\
    &\mspace{72mu}\times\;
    \Bigl(
      \bigl(
        p_{0, s+t}^{\omega}(0, y) \,+\, p^{\omega}_{s,s+t}(x, y)
      \bigr)\,
      \max_{z \in \mathbb{Z}^{d}} p_{s+t/2,s+t}^{\omega}(z,y)
    \Bigr)^{\!1/2}.
    \label{eq:gradient:chapman:final} 
  \end{align}
  Thus, in view of \eqref{eq:nash}, we obtain
  \begin{align}\label{eq:gradient:on-off-estimate}
    &\mean\Bigl[
      p_{0,s}^{\omega}(0,x)\,
      \bigl|p_{s,s+t}^{\omega}(x, y) - p_{0,s+t}^{\omega}(0, y)\bigr|^{2}
    \Bigr]
    \nonumber\\[.5ex]
    &\mspace{36mu}\leq\;
    c\, \Bigl(\frac{t}{2}\Bigr)^{\!-d/2}  
    \mean\Bigl[
      \Mean_{0,0}^{\omega}\Bigl[\Delta_{s,t/2}^{\omega}(0, X_{s})^{2}\Bigr]
      \bigl(p_{0,s+t}^{\omega}(0,y) + p_{s,s+t}^{\omega}(x,y)\bigr)
    \Bigr],
  \end{align}
  and \eqref{eq:gradient:claim} follows by combining \eqref{eq:gradient:on-off-estimate} with \eqref{eq:nash} and \eqref{eq:mean:entropy:est}.

  \textit{Step 2.}
  By using the stationarity of $\prob$ with respect to time shifts, \eqref{eq:shift:hk} and Minkowski's inequality, we obtain that, for any $s > 0$ and $(0,x) \in E_{d}$,
  \begin{align}
    \label{eq:split:gradient:hk}
    &
    \mean\Bigl[
      \big| \nabla^1_{\!x}\, p_{0, t}^{\omega}(0, y') \big|^{2}
    \Bigr]^{1/2}
    \nonumber\\[.5ex] 
    &\mspace{32mu}\overset{\!\!\eqref{eq:shift:hk}\!\!}{\;=\;}
    \mean\Bigl[
      \big| p_{s,s+t}^{\omega}(x, y') - p_{s,s+t}^{\omega}(0,y')\big|^{2}
    \Bigr]^{1/2}
    \nonumber\\[.5ex] 
    &\mspace{36mu}\leq\;
    \mean\Bigl[
      \big|p_{s,s+t}^{\omega}(x,y') - p_{0,s+t}^{\omega}(0,y')\big|^{2}
    \Bigr]^{1/2}
    \!+\,
    \mean\Bigl[
      \bigl|p_{s,s+t}^{\omega}(0,y') - p_{0,s+t}^{\omega}(0,y')\bigr|^{2}
    \Bigr]^{1/2}.
  \end{align}
  Set $K \ldef (1 \wedge \const[C]{lower:elliptic})/2$.  Then,
  \begin{align}
    \label{eq:gradient:lim:split:1}
    &\mean\Bigl[
      \bigl|p_{s,s+t}^{\omega}(x,y') - p_{0,s+t}^{\omega}(0,y')\bigr|^{2}
    \Bigr]
    \nonumber\\[.5ex]
    &\mspace{36mu}=\;
    \mean\biggl[
      \biggl(
        \frac{s^{-1} p_{0,s}^{\omega}(0,x) \vee K}
        {s^{-1} p_{0,s}^{\omega}(0,x) \vee K}
      \biggr)\,
      \bigl|p_{s,s+t}^{\omega}(x,y') - p_{0,s+t}^{\omega}(0,y')\bigr|^{2}
    \biggr]
    \nonumber\\[.5ex]
    &\mspace{36mu}\leq\;
    \frac{1}{K}\,  
    \mean\Bigl[
      s^{-1} p_{0,s}^{\omega}(0,x)\,
      \bigl|p_{s,s+t}^{\omega}(x,y') - p_{0,s+t}^{\omega}(0,y')\bigr|^{2}
    \Bigr]
    \nonumber\\[-1.1ex]
    \\[-1.9ex]
    &\mspace{62mu}+\,
    \frac{1}{K}\,
    \mean\Bigl[
      \bigl|s^{-1} p_{0,s}^{\omega}(0,x) \vee K - s^{-1} p_{0,s}^{\omega}(0,x)\bigr|
    \Bigr],
    \nonumber
  \end{align}
  where we used in the second term that $\bigl|p_{s,s+t}^{\omega}(x,y') - p_{0,s+t}^{\omega}(0,y')\bigr|^{2} \leq 1$.  Hence,
  \begin{align}
    &\mean\Bigl[
      \bigl|p_{s,s+t}^{\omega}(x,y') - p_{0,s+t}^{\omega}(0,y')\bigr|^{2}
    \Bigr]
    \nonumber\\[.5ex]  
    &\mspace{31mu}\overset{\!\!\!\eqref{eq:gradient:claim}\!\!\!}{\;\leq\;}
    \frac{1}{K}\,  
    \biggl(
      \frac{c}{t^{d+1}}
      +
      \mean\Bigl[
        \bigl|s^{-1} p_{0,s}^{\omega}(0,x) \vee K - s^{-1} p_{0,s}^{\omega}(0,x)\bigr|
      \Bigr]
    \biggr).
    \label{eq:gradient:split}  
  \end{align}
  Note that, due to Assumption~\ref{ass:integrability}, $\lim_{s \downarrow 0} s^{-1} p_{0,s}^{\omega}(0,x) = \omega_{0}(0,x) \geq \const[C]{lower:elliptic} > K$ for $\prob$-a.e.\ $\omega$.  Since $|s^{-1} p_{0,s}^{\omega}(0,x) \vee K - s^{-1} p_{0,s}^{\omega}(0,x)| \leq K$ for all $s \geq 0$, Lebesgue's dominated convergence theorem implies that
  \begin{align}
    \label{eq:gradient:lim:1}
    \lim_{s \searrow 0}
    \mean\Bigl[
      \bigl|s^{-1} p_{0,s}^{\omega}(0,x) \vee K - s^{-1} p_{0,s}^{\omega}(0,x)\bigr|
    \Bigr]
    \;=\;
    0.
  \end{align}
  Likewise, we obtain that
  \begin{align}
    &\limsup_{s \searrow 0}
    \mean\Bigl[
      \bigl|p_{s,s+t}^{\omega}(0,y') - p_{0,s+t}^{\omega}(0,y')\bigr|^{2}
    \Bigr]
    \nonumber\\
    &\mspace{36mu}=\;
    \limsup_{s \searrow 0}  
    \mean\biggl[
      \biggl(
        \frac{p_{0,s}^{\omega}(0,0) \vee K}
        {p_{0,s}^{\omega}(0,0) \vee K}
      \biggr)\,
      \bigl|p_{s,s+t}^{\omega}(0,y') - p_{0,s+t}^{\omega}(0,y')\bigr|^{2}
    \biggr]
    \nonumber\\[.5ex]
    &\mspace{31mu}\overset{\!\!\!\eqref{eq:gradient:claim}\!\!\!}{\;\leq\;}
    \limsup_{s \searrow 0}  
    \frac{1}{K}\,  
    \biggl(
      \frac{c s}{t^{d+1}}
      +
      \mean\Bigl[
        \bigl|p_{0,s}^{\omega}(0,0) \vee K - p_{0,s}^{\omega}(0,0)\bigr|
      \Bigr]
    \biggr)
    \;=\;
    0,
    \label{eq:gradient:lim:2}  
  \end{align}
  where we used in the last step again Lebesgue's dominated convergence theorem, since $|p_{0,s}^{\omega}(0,0) \vee K - p_{0,s}^{\omega}(0,0)| \leq K$ for all $s \geq 0$, and $\lim_{s \downarrow 0} p_{0,s}^{\omega}(0,0) = 1 > K$.  Thus, \eqref{eq:gradient:on-diag} follows from \eqref{eq:gradient:split}, \eqref{eq:gradient:lim:1} and \eqref{eq:gradient:lim:2}.  

  In order to show the last claim note that
  \begin{align*}
    \mean\bigl[|\nabla^1_x p^\omega_{0,t}(y',y)|^{2}\bigr]
    \;=\;
    \mean\bigl[|\nabla^1_x p^{\omega}_{-t,0}(y',y)|^{2}\bigr]
    \;=\;
    \mean\bigl[|\nabla^{2}_x p^{\tilde{\omega}}_{0,t}(y,y')|^{2}\bigr].
  \end{align*}
  Thus, the assertion follows by the same arguments used in \emph{Step 1} and \emph{Step 2}.
  \smallskip
  
  (ii) Again, in view of \eqref{eq:shift:hk} and Assumption~\ref{ass:law}-(i), it suffices to prove \eqref{eq:gradient:sum:on-diag} for $y=0$.  By interpolation, the assertion \eqref{eq:gradient:sum:on-diag} follows, once we have that shown that, for any $t \geq 1$,
  \begin{align*}
    \biggl(
      \sum_{y' \in \mathbb{Z}^{d}}
      \mean\Bigl[ 
        \bigl| \nabla^1_{\!x}\, p_{0, t}^{\omega}(y, y') \bigr|^p 
      \Bigr]
    \biggr)^{\!\!1/p}
    \;\leq\;
    \begin{cases}
      \const[C]{gradient:sum}\, t^{-1/2}, &p=1,
      \\[.5ex]
      \const[C]{gradient:sum}\, t^{-(d+2)/4}, &p=2.
    \end{cases}
  \end{align*}
  For $p = 2$ the proof of \eqref{eq:gradient:sum:on-diag} goes literally along the lines of the proof of Theorem~\ref{thm:gradient:est}-(i), except from the fact that \eqref{eq:gradient:claim} has to be replaced by
  \begin{align}\label{eq:gradient:sum:weighted}
    \sum_{y \in \mathbb{Z}^{d}}
    \mean\Bigl[
      p_{0,s}^{\omega}(0,x)\,
      \bigl| p_{s,s+t}^{\omega}(x,y) - p_{0,s+t}^{\omega}(0,y) \bigr|^{2}
    \Bigr]
    \;\leq\;
    c\, \frac{s}{t^{d/2+1}},
    \quad \forall\, x \in \mathbb{Z}^{d},\; s>0,
  \end{align}
  which follows immediately by summing \eqref{eq:gradient:on-off-estimate} over all $y \in \mathbb{Z}^{d}$ combined with \eqref{eq:mean:entropy:est}.  For $p=1$ we obtain instead of \eqref{eq:gradient:claim} the following estimate
  \begin{align*}
    \sum_{y \in \mathbb{Z}^{d}}
      \mean\Bigl[
       p_{0,s}^{\omega}(0,x)\,
      \bigl| p_{s,s+t}^{\omega}(x,y) - p_{0,s+t}^{\omega}(0,y) \bigr|
    \Bigr]
    \;\leq\;
    c\, \sqrt{\frac{s}{t}},
    \qquad \forall\, x \in \mathbb{Z}^{d},\; s>0,
  \end{align*}
  by first summing \eqref{eq:gradient:chapman} over all $y \in \mathbb{Z}^{d}$, and then applying the Cauchy-Schwarz inequality and \eqref{eq:mean:entropy:est}.  By a minor modification of the second step in the proof of Theorem~\ref{thm:gradient:est}-(i), the desired claim \eqref{eq:gradient:sum:on-diag} for $p=1$ follows.  Finally, we use that $\mean[|\nabla_x^1p_{0,t}^{\omega}(y',y)|^p] = \mean[|\nabla_x^{2}p_{0,t}^{\tilde{\omega}}(y,y')|^p]$ to prove the last result.
  \smallskip
  
  (iii) The proof relies on the same argument as in the \cite[Proof of (1.5b)]{DD05}.  First, by using the stationarity of $\prob$ with respect to space-shifts and \eqref{eq:shift:hk} we find that
  \begin{align}\label{eq:noebi12-1}
    \mean\Bigl[\nabla^{2}_{\!x}\nabla^{2}_{\!x'}\, p^{\omega}_{0,t}(y,y')\Bigr]
    \;=\;
    \mean\Bigl[\nabla^{1}_{\!-x}\nabla^{2}_{\!x'}\, p^{\omega}_{0,t}(y,y')\Bigr].
  \end{align}
  Hence, it is enough to prove \eqref{eq:hk:deriv:12}.  By using the Chapman-Kolmogorov equation, we get that
  \begin{align}\label{eq:property:mixed:deriv}
    \nabla^{1}_{\!-x}\nabla^{2}_{\!x'}\, p^{\omega}_{0,t}(y,y')
    \;=\;
    \sum_{z \in \mathbb{Z}^{d}} \Big(\nabla^1_{\!-x}\, p^{\omega}_{0,t/2}(y,z)\Big)
    \Bigl(\nabla^{2}_{\!x'}\, p^{\omega}_{t/2,t}(z,y') \Bigr).
  \end{align}
  By applying the Cauchy-Schwarz inequality, we obtain that
  \begin{align*}
    &\mean\Bigl[
      \bigl|\nabla^{1}_{\!-x}\nabla^{2}_{\!x'}\, p^{\omega}_{0,t}(y,y') \bigr|
    \Bigr]
    \\[.5ex]
    &\mspace{36mu}\leq\;
    \biggl(
      \sum_z
      \mean\Bigl[\bigl| \nabla^1_{\!-x}\, p_{0,t/2}^{\omega}(y, z)\bigr|^{2}\Bigr]
    \biggr)^{\!\!1/2}
    \biggl(
      \sum_z
      \mean\Bigl[
        \bigl| \nabla_{\!x'}^{2}\, p_{t/2,t}^{\omega}(z,y')\bigr|^{2}
      \Bigr]
    \biggr)^{\!\!1/2},
  \end{align*}
  and \eqref{eq:hk:deriv:12} follows from \eqref{eq:gradient:sum:on-diag}.
\end{proof}
\begin{proof}[Proof of Corollary~\ref{cor:time:derivative}]
  Since, $\lim_{s \downarrow 0} s^{-1} p_{0,s}^{\omega}(0,x) = \omega_{0}(0,x)$ for $\prob$-a.e.\ $\omega$ and the heat kernel is continuous in time in both arguments, cf.~Proposition~\ref{prop:backward}, an application of Fatou's lemma yields that, for any $t > 0$ and $\{0, x\} \in E_{d}$,
  \begin{align}\label{eq:gradient:sum:weighted:1}
    &\sum_{z \in \mathbb{Z}^{d}}
    \mean\Bigl[
      \omega_{0}(0,x)\,
      \bigl| \nabla_{\!x}^{1} p_{0,t}^{\omega}(0, z) \bigr|^{2}
    \Bigr]
    \nonumber\\[.5ex]
    &\mspace{36mu}\leq\;
    \liminf_{s \downarrow 0}
    \sum_{z \in \mathbb{Z}^{d}}
    \mean\Bigl[
      s^{-1} p_{0,s}^{\omega}(0, x)\,
      \bigl| p_{s,s+t}^{\omega}(x, z) - p_{0,s+t}^{\omega}(0, z) \bigr|^{2}
    \Bigr]
    \overset{\!\!\!\eqref{eq:gradient:sum:weighted}\!\!\!}{\;\leq\;}
    \frac{c}{t^{d/2+1}},
  \end{align}
  By noting that, for any $z \in \mathbb{Z}^{d}$,
  \begin{align*}
    &\mean\Bigl[
      \omega_{t}(0,x)\,
      \bigl| \nabla_{\!x}^{2} p_{0,t}^{\omega}(z, 0) \bigr|^{2}
    \Bigr]
    \;=\;
    \mean\Bigl[
      \tilde{\omega}_{t}(0,x)\,
      \bigl| \nabla_{\!x}^{1} p_{0,t}^{\tilde{\omega}}(0, z) \bigr|^{2}
    \Bigr]
  \end{align*}
  and repeating the proof of \eqref{eq:gradient:sum:weighted} with $\omega$ replaced by $\tilde{\omega}$, we also obtain that
  \begin{align}\label{eq:gradient:sum:weighted:2}
    \sum_{z \in \mathbb{Z}^{d}}
    \mean\Bigl[
      \omega_{t}(0,x)\,
      \bigl| \nabla_{\!x}^{2} p_{0,t}^{\omega}(z, 0) \bigr|^{2}
    \Bigr]
    \;\leq\;
    \frac{c}{t^{d/2+1}}.
  \end{align}
  Let us now address the proof of \eqref{eq:time:derivative}.  First, by using the Chapman-Kolmogorov equation, the stationarity of $\prob$ with respect to time-space shifts as given in Assumption~\ref{ass:law}-(i) and \eqref{eq:shift:hk}, we obtain, for any $y, y' \in \mathbb{Z}^{d}$, $t \in (0, \infty)$ and $s \in (0, t)$,
  \begin{align*}
    \bar{p}_{0, t}(y, y')
    \;=\;
    \sum_{z \in \mathbb{Z}^{d}}
    \mean\bigl[ p_{0, s}^{\omega}(y-z, 0)\, p_{s, t}^{\omega}(0, y'-z) \bigr].
  \end{align*}
  Recall that the forward equation, see Proposition~\ref{prop:forward}, implies that, for $\prob$-a.e.\ $\omega$, $\partial_{t} p_{s, t}^{\omega}(0, z) = (\mathcal{L}^{\omega} p_{s, t}(0, \cdot))(z)$ for any $z \in \mathbb{Z}^{d}$ and a.e.\ $t \in (s, \infty)$. Since, by Assumption~\ref{ass:time:derivative}, $|\mathcal{L}_{t}^{\omega} p_{s, t}^{\omega}(0, \cdot)(y)| \leq \const[C]{lower:elliptic:time} \mu_{0}^{\omega}(y)$ for any $y \in \mathbb{Z}^{d}$, $\prob$-a.e.\ $\omega$, and any $t \in (s, \infty)$, with $\mean[\mu_{0}^{\omega}(y)] < \infty$, we can interchange the partial derivative and the expected value. This yields, for almost every $t \in (0, \infty)$,
  \begin{align*}
    \frac{\partial}{\partial t}\, \bar{p}_{0, t}(y, y')
    &\;=\;
    \sum_{z \in \mathbb{Z}^{d}}
    \mean\bigl[
      p_{0, s}^{\omega}(y-z, 0)\,
      \bigl(
        \mathcal{L}_{t}^{\omega} p_{s, t}^{\omega}(0, \cdot )
      \bigr)(y'-z)
    \bigr]
    \\[.5ex]
    &\;=\;
    -\frac{1}{2} \sum_{z \in \mathbb{Z}^{d}} \sum_{x \sim 0}
    \mean\bigl[
      \omega_{t}(y', y'+x) \nabla_{\!x}^{1} p_{0, s}^{\omega}(y, z) \nabla_{\!x}^{2} p_{s, t}^{\omega}(z, y')
    \bigr].
  \end{align*}
  Thus, by choosing $s = t/2$, using the Cauchy-Schwarz inequality and that
  \begin{align*}
    \mean\bigl[
      \omega_{t}(y', y'+x) \bigl| \nabla_{\!x}^{1} p_{0, t/2}^{\omega}(y, z) \bigr|^{2}
    \bigr]
    &\;=\;
    \mean\bigl[
      \omega_{t}(0, x) \bigl| \nabla_{\!x}^{1} p_{0, t/2}^{\omega}(y-y', z) \bigr|^{2}
    \bigr]
    \\
    &\;\leq\;
    \const[C]{lower:elliptic:time}\,
    \mean\bigl[
      \omega_{0}(0, x) \bigl| \nabla_{\!x}^{1} p_{0, t/2}^{\omega}(y-y', z) \bigr|^{2}
    \bigr]
    \intertext{and}
    \mean\bigl[
      \omega_{t}(y', y'+x) \bigl| \nabla_{\!x}^{2} p_{t/2, t}^{\omega}(z, y') \bigr|^{2}
    \bigr]
    &\;=\;
    \mean\bigl[
      \omega_{t}(0, x) \bigl| \nabla_{\!x}^{2} p_{t/2, t}^{\omega}(z, 0) \bigr|^{2}
    \bigr],
  \end{align*}
  it follows from \eqref{eq:gradient:sum:weighted:1} and \eqref{eq:gradient:sum:weighted:2} that
  \begin{align*}
    \Bigl|
      \frac{\partial}{\partial t}\, \bar{p}_{0, t}(y, y')
    \Bigr|
    &\;\leq\;
    \frac{1}{2} \sum_{z \in \mathbb{Z}^{d}} \sum_{x \sim 0}
    \mean\bigl[
      \omega_{t}(y', y'+x) |\nabla_{\!x}^{1} p_{0, s}^{\omega}(y, z)| |\nabla_{\!x}^{2} p_{s, t}^{\omega}(z, y')|
    \bigr]
    \;\leq\; 
    \frac{c d \sqrt{\const[C]{lower:elliptic:time}}}{t^{d/2 + 1}},
  \end{align*}
  and by setting $\const[C]{time:elliptic} \ldef c d \sqrt{\const[C]{lower:elliptic:time}}$ the assertion follows. 
\end{proof}

\subsection{Off-diagonal upper bound and Green kernel estimates}
\label{subsec:gradient:hk+green}
In this subsection, we aim to establish estimates for the quenched and annealed Green's function (or potential kernel density) defined by
\begin{align}
  \mathbb{Z}^{d} \times \mathbb{Z}^{d} \ni (y, y')
  \;\longmapsto\;
  \int_{0}^{\infty} p_{0, t}^{\omega}(y, y')\, \md t
  \in [0, \infty],
\end{align}
along with its first and second derivatives. Such estimates are crucial for analysing the time-space covariances of the height of an interface separating two pure thermodynamic phases, commonly referred to in the literature as the \emph{Landau-Ginzburg $\nabla \phi$-model}, cf.~\cite{GOS01, Fu05, DD05, AT21}.

While Lemma~\ref{lem:Nashclas} immediately implies that, for any $d \geq 3$, the Green's function is finite, such on-diagonal estimate of the quenched heat kernel is not sufficient to derive finer estimates.  For this purpose, we introduce the following deterministic kernel.
\begin{definition}\label{def:offdiag}
  For some $\alpha \in (0, \infty)$ define, for any $t > 0$ and $y \in \mathbb{Z}^{d}$,
  \begin{align}\label{def:kernel}
    k_{\alpha}(t, y)
    \;\ldef\;
    \begin{cases}
      |y|^{-d/2 - \alpha/2},
      &\text{if }t \leq |y|
      \\[1ex]
      (1 \vee t)^{-d/2}\,
      \biggl(1 \vee \dfrac{|y|}{\sqrt{t}}\biggr)^{\!\!-\alpha},
      \quad
      &\text{if }t \geq |y|.
    \end{cases}
  \end{align}
\end{definition}
\begin{remark}\label{rem:kernel}
  Observe that, for any $\alpha \in (0, \infty)$ and $t \geq |y|^{2}$, $k_{\alpha}(t, y) = (1 \vee t)^{-d/2}$. 
\end{remark}
In the following lemma we show that the above defined kernel is stable under convolution in the following sense.  Since the proof is rather technical we postpone it to Appendix~\ref{appendix:technical}.  
\begin{lemma}\label{lem:kernel:convolution}
  Suppose that $\alpha > d$.  Then, there exists $\const[C]{kernel:convolution} \equiv \const[C]{kernel:convolution}(d, \alpha) < \infty$ such that, for any $y \in \mathbb{Z}^{d}$ and $t > 0$,
  \begin{align}\label{eq:kernel:convolution}
    \sum_{z \in \mathbb{Z}^{d}} k_{\alpha}(t/2, z)\, k_{\alpha}(t/2, y-z)
    \;\leq\;
    \const[C]{kernel:convolution}\, k_{\alpha}(t, y). 
  \end{align}
\end{lemma}
Now, let us impose the following quenched assumption on the off-diagonal behaviour of the heat kernel.
\begin{assumption}\label{ass:offdiag}
  For some $\alpha, \beta \in (0, \infty)$ assume that there exist non-random constants $c_{1}, c_{2} \in (0, \infty)$ and a positive random variable $N(\omega)$ such that the following holds true: for $\prob$-a.e.~$\omega$ and for all $s \geq 0$, $t > 0$ and $y, y' \in \mathbb{Z}^{d}$ such that $|y-y'| \vee \sqrt{t} \geq N_{s, y}(\omega) \ldef N(\omega) \circ \tau_{s,y}$, 
  \begin{align}\label{eq:ass:hk:offdiag}
    p_{s, s+t}^{\omega}(y, y')
    \;\leq\;
    c_{1}\, k_{\alpha}(t, y'- y).
  \end{align}
  Moreover,
  \begin{align}\label{eq:prob:N_y}
    \prob\bigl[N_{s, y}(\omega) \geq L\bigr]
    \;\leq\;
    c_{2}\, L^{-\beta},
    \quad \forall\, L > 0.
  \end{align}
\end{assumption}
\begin{remark}
  (i) The assumptions \eqref{eq:ass:hk:offdiag} and \eqref{eq:prob:N_y} that we impose on the off-diagonal behaviour of the heat kernel on large time-space scales, as well as the tail behaviour of the distribution of the random variables involved, are formulated in a weakest possible sense. Nevertheless, these assumptions are sufficient to deduce the correct order of decay for the corresponding Green's function and its first and second derivatives, cf.~Proposition~\ref{prop:Green:function:estimates}.
    
  (ii) For time-independent weights, the heat kernel upper bound \eqref{eq:ass:hk:offdiag} can typically be upgraded to a Gaussian-type bound, which is expressed as
  \begin{align}\label{eq:heat:kernel:Gaussian:ub}
    p_{s, s+t}^{\omega}(y, y')
    \;\leq\;
    \begin{cases}
      c_{2} t^{-d/2} \exp\biggl(-c_{3} \dfrac{d^{\omega}(y, y')^{2}}{t}\biggr),
      \quad &d^{\omega}(y, y') \leq c_{1} t
      \\[2ex]
      c_{2} t^{d/2}
      \exp\biggl(-c_{4} d^{\omega}(y, y')\biggl(1 \vee \log\dfrac{d^{\omega}(y, y')}{t}\biggr)\biggr),
      \quad &d^{\omega}(y, y') \geq c_{5} t
    \end{cases}
  \end{align}
  for some constants $c_{1}, \ldots, c_{5} \in (0, \infty)$. The distance function, $d^{\omega}$, appearing therein is defined, for any distinct points $y, y' \in \mathbb{Z}^{d}$, by
  \begin{align}\label{eq:def:chem:distance}
    d^{\omega}(y, y')
    \;\ldef\;
    \inf_{\gamma}\Biggl\{
      \sum_{i=0}^{l_{\gamma}-1} \biggl(1 \vee \frac{1}{\omega(z_{i}, z_{i+1})}\biggr)^{\!\!1/2}
    \Biggr\}
  \end{align}
  where the infimum is taken over all nearest-neighbour paths, $\gamma = (z_{0}, \cdots, z_{l_{\gamma}})$, connecting $y$ to $y'$. It is well known that such Gaussian bounds hold, for instance, for the VSRW on (i.i.d.\ and dependent) supercritical percolation clusters (see \cite[Lemma~1.1, Theorem~1]{Ba04}, \cite[Corollary~1.14]{Sa17}), for VSRW with i.i.d.\ conductances that are uniformly bounded away from zero (see \cite[Theorem~2.19]{BD10}), and for VSRW with ergodic conductances that satisfy certain integrability conditions (see \cite[Theorem~3.2]{ADS19}) up to a multiplicative factor $(1 + |y-y'|/\sqrt{t})^{\gamma}$, $\gamma \in (0, \infty)$.
  
  For sufficiently mixing environments as considered in e.g.~\cite{Ba04, Sa17}, the tail behaviour \eqref{eq:prob:N_y} can also be upgraded to a stretched exponential decay. Moreover, for such environments it has been shown that the distance $d^{\omega}$ is comparable with the usual graph distance on $\mathbb{Z}^{d}$, see e.g.\cite[Theorem~1]{AP96}, \cite[Theorem~1.3]{DRS14}~.
\end{remark}
%
%
%
\begin{example}
  (i) Consider $\omega \equiv \big\{\omega_{t}(e) : e \in E_{d},\, t \in \mathbb{R} \big\}$ that satisfy Assumption~\ref{ass:law}.  Additionally, assume that there exist i.i.d. $\bar{\omega} \equiv \big\{\bar{\omega}(e) : e \in E_{d} \big\}$  which is bounded from below (i.e.\ $\prob[\bar{\omega}(e) \geq \const[C]{lower:elliptic}] = 1$ for all $e \in E_{d}$) and $0 < c_{1} < c_{2} < \infty$ such that
  \begin{align}\label{eq:fneobi1c}
    \prob\bigl[
      c_{1}\, \bar{\omega}(e) \leq \omega_{t}(e) \leq c_{2}\, \bar{\omega}(e)\;
      \forall\, t \in \mathbb{R}
    \bigr]
    \;=\;
    1
    \qquad \forall e \in E_{d}.
  \end{align}
  Note that \eqref{eq:fneobi1c} can be deduced from Assumption~\ref{ass:time:derivative} if the maps $t \mapsto \omega_{t}(e)$ are separable for every $e \in E_{d}$.
  Then, for any $s \geq 0$ and $y \in \mathbb{Z}^{d}$, there exists a random variable whose distribution has sub-exponential tails such that, for $\prob$-almost every $\omega$ and any $t > 0$, $y' \in \mathbb{Z}^{d}$ satisfying $|y - y'| \vee \sqrt{t} \geq N_{s, y}(\omega)$, the corresponding heat kernel $p^{\omega}_{s, s+t}(y, y)$ satisfies a Gaussian upper bound as stated in \eqref{eq:heat:kernel:Gaussian:ub}, with $d^{\omega}$ replaced by the usual graph metric on $\mathbb{Z}^{d}$.  Hence, Assumption~\ref{ass:offdiag} holds, and the results in this section are applicable.  Indeed, by setting for any $\psi \in \ell^{\infty}(\mathbb{Z}^{d})$
  \begin{align*}
    \Lambda^{\omega}(\psi)
    \;\ldef\;
    \sup_{y \in \mathbb{Z}^{d}} \sum_{x \sim 0} \omega_{t}(y, y+x)\,
    \bigl(\cosh(\nabla_{\!x} \psi(y)) - 1\bigr),
  \end{align*}
  and using that, $\prob$-a.s., $\Lambda^{\omega}(\psi) \leq c_{2} \Lambda^{\bar{\omega}}(\psi)$, we obtain from \cite[Proposition~5]{Da93} the following Davies-Gaffney type estimate
  \begin{align*}
    p_{s, s+t}(y, y')
    \;\leq\;
    \inf_{\psi \in \ell^{\infty}(\mathbb{Z}^{d})}
    \exp\bigl( \psi(y') - \psi(y) + c_{2} \Lambda^{\bar{\omega}}(\psi) t \bigr)
  \end{align*}
  Hence, \cite[Theorem 2.3]{BD10} is applicable. Furthermore, the exit time estimates in \cite[Proposition 2.18]{BD10} are based on the strong Markov property of the process, which also holds true in our example. Consequently, we can prove (4.6)--(4.8) of \cite{BD10} (which is stronger than Assumption \ref{ass:offdiag}) by following the same arguments as in the proofs in \cite[Sections 2 and 4]{BD10}, with appropriate adjustments to the constants in the inequalities.

  (ii) Note that the case of time-dependent weights, $\omega$, are uniformly elliptic, that is, there exists $0 < c_{1} < c_{2} < \infty$ such that $\prob\bigl[c_{1} \leq \omega_{t}(e) \leq c_{2}\; \forall\, t \in \mathbb{R}\bigr] = 1$ for all $e \in E_{d}$, is a special case of \eqref{eq:fneobi1c}. Indeed, as shown in e.g.~\cite[Proposition~4.2]{DD05} or \cite[Proposition~B.3]{GOS01}, the heat kernel exhibits Gaussian upper bounds with $N_{s, y} = 1$ for all $s \geq 0$, $y \in \mathbb{Z}^{d}$.
\end{example}
\begin{lemma}\label{lem:hk:offdiag:p}
  Suppose that Assumption~\ref{ass:law}, \ref{ass:integrability} hold and let $p \in [1, \infty)$.  Further, let Assumption~\ref{ass:offdiag} be satisfied for some $\alpha > 0$ and $\beta \geq p \alpha/2$.  Then, there exists $\const[C]{mean:p:heatkernel} < \infty$ such that for any $y, y' \in \mathbb{Z}^{d}$ and $t > 0$
  \begin{align}\label{eq:hk:offdiag:p}
    \mean\bigl[p_{0,t}^{\omega}(y, y')^p\bigr]^{1/p}
    \;\leq\:
    \const[C]{mean:p:heatkernel}\, k_{\alpha}(t, y'-y).
  \end{align}
\end{lemma}
\begin{proof}
  In view of \eqref{eq:shift:hk}, Assumption~\ref{ass:law}-(i) and Remark~\ref{rem:kernel}, it suffices to prove \eqref{eq:hk:offdiag:p} for $y = 0$ and $y' \in \mathbb{Z}^{d}$ such that $|y'| \geq 1$ and $0 < t \leq |y'|^{2}$. For this purpose, consider the event $A_{0}(\omega) \ldef \{N_{0}(\omega) \leq |y'| \vee \sqrt{t}\}$.  Then, by Minkowski's inequality, we get that
  \begin{align*}
    \mean\bigl[p_{0,t}^{\omega}(0, y')^p\bigr]^{1/p}
    &\overset{\!\!\eqref{eq:nash}\!\!}{\;\leq\;}
    \mean\bigl[p_{0, t}^{\omega}(0, y')^{p} \indicator_{A_{0}(\omega)}\bigr]^{1/p}
    + c\, (1 \vee t)^{-d/2} \prob\bigl[A_{0}^{\mathrm{c}}(\omega)\bigr]^{1/p}
    \\[.5ex]
    &\overset{\!\!\eqref{eq:ass:hk:offdiag}\!\!}{\;\leq\;}
    c\, k_{\alpha}(t, y') + c\, (1 \vee t)^{-d/2} |y'|^{-\beta/p}.
  \end{align*}
  Since $\beta \geq p\alpha / 2$, the assertion follows.
\end{proof}
\begin{prop}\label{prop:gradient:offdiag}
  Suppose that Assumption \ref{ass:law} and \ref{ass:integrability} hold.  Additionally, let Assumption~\ref{ass:offdiag} be satisfied for some $\alpha \geq 2$ and $\beta > d + \alpha/2$.  Then, there exists $\const[C]{mean:gradient} < \infty$ such that for any $(0,x) \in E_{d}$, $y, y' \in \mathbb{Z}^{d}$ and $t > 0$,
  \begin{align}\label{eq:gradient:offdiag}
    \mean\Bigl[
      \bigl| \nabla^1_{\!x}\, p_{0, t}^{\omega}(y, y')\bigr|^{2}
    \Bigr]^{1/2}
    \;\leq\;
    \const[C]{mean:gradient}\, (1 \vee t)^{-1/2} k_{\alpha/2}(t, y' - y). 
  \end{align}
\end{prop}
\begin{proof}
  Again, in view of \eqref{eq:shift:hk} and Assumption~\ref{ass:law}-(i), it suffices to prove \eqref{eq:gradient:offdiag} for $y=0$ and $y' \in \mathbb{Z}^{d}$. Further, observe that for any $y' \in \mathbb{Z}^{d}$ and $t \geq 1$ such that $|y'|^{2} \leq t$, the assertion \eqref{eq:gradient:offdiag} follows directly from Theorem~\ref{thm:gradient:est}-(i). This is due to the fact that, as noted in Remark~\ref{rem:kernel}, $k_{\alpha}(t, y') = t^{-d/2}$. Moreover, for $y' = 0$ and $t \in (0,1]$ the assertion is trivial. Thus, it remains to prove \eqref{eq:gradient:offdiag} for any $|y'| \geq 1$ and $0 < t \leq |y'|^{2}$. For this purpose, set $A_z(\omega) \ldef \{N_z(\omega) \leq |z-y'| \vee \sqrt{t}\}$, $z \in \{0,x\}$.  We proceed by considering the cases $0 < t \leq |y'|$ and $|y'| \leq t \leq |y'|^{2}$ separately. 
  
  First, consider the regime $0 < t < |y'|$.  Clearly, for $|y'| = 1$ the assertion follows by bounding the heat kernel from above by one.  Therefore, it remains to consider the case that $|y'| > 1$.  By Minkowski's inequality, we get for any $(0, x) \in E_{d}$
  \begin{align*}
    &\mean\Bigl[
      \bigl| \nabla_{\!x}^1\, p_{0, t}^{\omega}(0, y') \bigr|^{2}
    \Bigr]^{1/2}
    \nonumber\\[.5ex]
    &\mspace{32mu}\overset{\!\!\eqref{eq:nash}\!\!}{\;\leq\;}
    \mean\Bigl[
      \bigl|
        \nabla_{\!x}^{1} p_{0, t}^{\omega}(0, y') \bigr|^{2}\,
      \indicator_{A_{0}(\omega) \cap A_{x}(\omega)}
    \Bigr]^{1/2}
    +\,
    c\, (1 \vee t)^{-d/2}\,
    \prob\bigl[A_{0}^{\mathrm{c}}(\omega) \cup A_{x}^{\mathrm{c}}(\omega)\bigr]^{1/2}
    \nonumber\\[.5ex]
    &\mspace{32mu}\overset{\!\!\eqref{eq:ass:hk:offdiag}\!\!}{\;\leq\;}
    c\, |y'|^{-d/2-\alpha/2}
    +\,
    c\, (1 \vee t)^{-d/2}\, |y'|^{-\beta/2},
  \end{align*}
  where we used in the last step that $\sqrt{t} \vee |y' - x| \geq (\sqrt{t} \vee |y'|)/2$.  Indeed, since $y' \in \mathbb{Z}^{d}$ with $|y'| > 1$ and $(0, x) \in E_{d}$ we have that $|y'-x| \geq 1 = |x|$ which implies that
  \begin{align*}
    \sqrt{t} \vee |y'|
    \;\leq\;
    \sqrt{t} \vee \bigl(|y' - x| + |x|\bigr)
    \;\leq\;
    2\bigl(\sqrt{t} \vee |y' - x|\bigr).
  \end{align*}
  Since $\alpha \geq 2$ and $\beta \geq d + \alpha/2$ we therefore obtain the desired bound, namely
  \begin{align*}
    \mean\Bigl[
      \bigl| \nabla_{\!x}^1\, p_{0,t}^{\omega}(0, y') \bigr|^{2}
    \Bigr]^{1/2}
    \;\leq\;
    c (1 \vee t)^{-1/2}\, k_{\alpha/2}(t, y').
  \end{align*}
  
  Next, consider the case that $|y'| \leq t \leq |y'|^{2}$.  By using similar arguments as in the proof of Theorem~\ref{thm:gradient:est}-(i) (\textit{Step 2}), we first obtain that, for any $s > 0$ and $(0, x) \in E_{d}$,
  \begin{align*}
    &\mean\Bigl[
      \bigl| \nabla_{\!x}^1\, p_{0,t}^{\omega}(0,y')\bigr|^{2}
    \Bigr]^{1/2}
    \\
    &\mspace{32mu}\overset{\!\!\eqref{eq:split:gradient:hk}\!\!}{\;\leq\;}
    \mean\Bigl[
      \big|p_{s,s+t}^{\omega}(x,y') - p_{0,s+t}^{\omega}(0,y')\big|^{2}
    \Bigr]^{1/2}
    \!+\,
    \mean\Bigl[
      \bigl|p_{s,s+t}^{\omega}(0,y') - p_{0,s+t}^{\omega}(0,y')\bigr|^{2}
    \Bigr]^{1/2}.
  \end{align*}
  Hence,
  \begin{align}
    \label{eq:hk:offdiag:est0}
    \mean\Bigl[
      \bigl| \nabla_{\!x}^1\, p_{0,t}^{\omega}(0,y')\bigr|^{2}
    \Bigr]^{1/2}
    \overset{\!\!\eqref{eq:gradient:lim:2}\!\!}{\;\leq\;}
    \limsup_{s \searrow 0}  
    \mean\Bigl[
      \bigl| p_{s,s+t}^{\omega}(x,y') - p_{0,s+t}^{\omega}(0,y')\bigr|^{2}\,
    \Bigr]^{1/2}.
  \end{align}
  Thus, by using again Minkowski's inequality, \eqref{eq:nash} and \eqref{eq:prob:N_y}, we further get that, for any $s > 0$,
  \begin{align*}
    &
    \mean\Bigl[
      \bigl| p_{s, s+t}^{\omega}(x, y') - p_{0, s+t}^{\omega}(0, y')\bigr|^{2}\,
    \Bigr]^{1/2}
    \nonumber\\[.5ex]
    &\mspace{36mu}\leq\;
    \mean\Bigl[
      \bigl| p_{s, s+t}^{\omega}(x, y') - p_{0, s+t}^{\omega}(0, y')\bigr|^{2}\,
      \indicator_{A_{0}(\omega) \cap A_{x}(\omega)}
    \Bigr]^{1/2}
    +\,
    c\, (1 \vee t)^{-d/2}\, |y'|^{-\beta/2}.
  \end{align*}
  By rewriting the first term as in \eqref{eq:gradient:lim:split:1} and using \eqref{eq:gradient:lim:1}, we find that 
  \begin{align}
    &\limsup_{s \searrow 0}
    \mean\Bigl[
      \bigl| p_{s,s+t}^{\omega}(x,y') - p_{0,s+t}^{\omega}(0,y')\bigr|^{2}\,
    \Bigr]^{1/2}
    \nonumber\\
    &\mspace{36mu}\leq\;
    \limsup_{s \searrow 0}
    \biggl(
      \frac{1}{K}  
      \mean\Bigl[
        s^{-1} p_{0, s}^{\omega}(0, x)\,
        \bigl| p_{s, s+t}^{\omega}(x,y') - p_{0, s+t}^{\omega}(0,y')\bigr|^{2}\,
        \indicator_{A_{0}(\omega) \cap A_{x}(\omega)}
      \Bigr]
    \biggr)^{\!\!1/2}
    \nonumber\\
    &\mspace{72mu}+\, 
    c\, t^{-d/2}\, |y'|^{-\beta/2},
    \label{eq:hk:offdiag:est1}  
  \end{align}
  where $K \ldef (1 \wedge \const[C]{lower:elliptic})/2$.  Moreover, we used in the second step that whenever $t \geq |y'| \geq 1$ it holds that $\sqrt{t} \vee |x-y'| \geq (\sqrt{t} \vee |y'|)/2$ for all $y' \in \mathbb{Z}^{d}$ and $(0,x) \in E_{d}$.  Note that $t^{-d/2} |y'|^{-\beta/2} \leq t^{-(d+1)/2} \bigl(|y'|/\sqrt{t}\bigr)^{-\alpha/2}$ for all $1 \leq |y'| \leq t \leq |y'|^{2}$ if $\beta \geq 1 + \alpha/2$.  Moreover, by combining the estimate \eqref{eq:gradient:chapman:final} with \eqref{eq:nash}, \eqref{eq:ass:hk:offdiag} and \eqref{eq:mean:entropy:est}, we obtain for any $s > 0$ and $|y'| \leq t \leq |y'|^{2}$ that
  \begin{align}
    &\mean\Bigl[
      p_{0, s}^{\omega}(0, x)\,
      \bigl|p_{s, s+t}^{\omega}(x, y') - p_{0, s+t}^{\omega}(0, y')\bigr|^{2}\,
      \indicator_{A_{0}(\omega) \cap A_{x}(\omega)}
    \Bigr]
    \nonumber\\[.5ex]
    &\mspace{32mu}\overset{\!\!\eqref{eq:nash}\!\!}{\;\leq\;}\,
    c\, t^{-d/2}\,  
    \mean\Bigl[
      \Mean_{0, 0}^{\omega}\Bigl[\Delta_{s, t/2}^{\omega}(0, X_{s})^{2}\Bigr]
      \bigl(p_{0, s+t}^{\omega}(0, y') + p_{s, s+t}^{\omega}(x, y')\bigr)\,
      \indicator_{A_{0}(\omega) \cap A_{x}(\omega)}
    \Bigr]
    \nonumber\\[.5ex]
    &\mspace{32mu}\overset{\!\!\eqref{eq:ass:hk:offdiag}\!\!}{\;\leq\;}
    c\, t^{-d}\, 
    \biggl(\frac{|y'|}{\sqrt{t}}\biggr)^{\!\!-\alpha}
    \mean\Bigl[
      \Mean_{0, 0}^{\omega}\Bigl[\Delta_{s, t/2}^{\omega}(0, X_{s})^{2}\Bigr]
    \Bigr]
    \nonumber\\[.5ex]
    &\mspace{32mu}\overset{\!\!\eqref{eq:mean:entropy:est}\!\!}{\;\leq\;}
    c\, s\, t^{-(d+1)}\,
    \biggl(\frac{|y'|}{\sqrt{t}}\biggr)^{\!\!-\alpha}. 
    \label{eq:hk:offdiag:est2}  
  \end{align}
  Thus, the assertion follows by combining \eqref{eq:hk:offdiag:est0}, \eqref{eq:hk:offdiag:est1} and \eqref{eq:hk:offdiag:est2}.
\end{proof}
\begin{prop}\label{prop:Green:function:estimates}
  Suppose that Assumption~\ref{ass:law} and \ref{ass:integrability} hold true.
  \begin{enumerate}[(i)]
    \item If $d \geq 3$ and, for any $p \geq 1$, Assumption~\ref{ass:offdiag} is satisfied for $\alpha > d-2$ and $\beta \geq p\alpha /2$, then there exists $C_{12} < \infty$ such that, for any $y, y' \in \mathbb{Z}^{d}$ with $y \ne y'$,
    \begin{align}\label{eq:green:p}
      \mean\Biggl[
        \biggl|
          \int_{0}^{\infty} p_{0, t}^{\omega}(y, y')\, \md t
        \biggr|^p
      \Biggr]^{1/p}
      \;\leq\;
      C_{12}\, |y - y'|^{2-d}.
    \end{align}
    \item If $d \geq 2$ and Assumption~\ref{ass:offdiag} is satisfied for $\alpha > 2(d-1)$ and $\beta > d + \alpha/2$, then there exists $C_{13} < \infty$ such that, for any $(0, x) \in E_{d}$ and $y, y' \in \mathbb{Z}^{d}$ with $y \ne y'$,
    \begin{align}\label{eq:gradient:green:1}
      \mean\Biggl[
        \biggl|
          \int_{0}^{\infty} 
            \nabla_{\!x}^{1}\, p_{0, t}^{\omega}(y, y')\, 
          \md t
        \biggr|^{2}
      \Biggr]^{1/2}
      \;\leq\;
      C_{13}\, |y - y'|^{1-d}.
    \end{align}
    \item If $d \geq 1$ and Assumption~\ref{ass:offdiag} is satisfied for $\alpha > 2d$ and $\beta > d + \alpha/2$, then there exists $C_{14} < \infty$ such that, for any $(0, x), (0, x') \in E_{d}$ and $y, y' \in \mathbb{Z}^{d}$ with $y \ne y'$,
    \begin{align}\label{eq:gradient:green:12}  
      \mean\Biggl[
        \biggl|
          \int_{0}^{\infty} 
            \nabla_{\!x}^1\, \nabla_{\!x'}^{2}\, p_{0, t}^{\omega}(y, y')\,
          \md t
        \biggr|
      \Biggr]
      \;\leq\;
      C_{14}\, |y - y'|^{-d}.
    \end{align}
  \end{enumerate}
\end{prop}
\begin{proof}
  First of all notice that by Lemma~\ref{lem:Nashclas} and Theorem~\ref{thm:gradient:est} the integrals appearing on the left-hand side of \eqref{eq:green:p}, \eqref{eq:gradient:green:1} and \eqref{eq:gradient:green:12}, respectively, are well defined.  Moreover, in view of \eqref{eq:shift:hk} and Assumption~\ref{ass:law}-(i) it suffices to show the assertion for $y = 0$.  Moreover, set $r \ldef |y'|$ and assume w.l.o.g.\ that $r \geq 4$.

  (i) By applying Minkowski's integral inequality, cf.\ \cite[Theorem~202]{HLP88}, together with the estimate \eqref{eq:hk:offdiag:p}, we obtain that, for any $p \in [1, \infty)$, $\alpha > 0$ and $\beta \geq p \alpha / 2$
  \begin{align*}
    \mean\Biggl[
     \biggl|\int_{0}^{\infty} p_{0, t}^{\omega}(0, y')\, \md t\biggr|^p
    \Biggr]^{1/p}
    \;\leq\;
    \int_{0}^{\infty}
      \mean\Bigl[\bigl| p^{\omega}_{0,t}(0,y') \bigr|^p\Bigr]^{1/p}
    \md t
    \overset{\!\!\eqref{eq:hk:offdiag:p}\!\!}{\;\leq\;}
    \const[C]{mean:p:heatkernel}\,
    \int_{0}^{\infty} k_{\alpha}(t,y')\, \md t.
  \end{align*}
  Thus, by estimating the resulting integral on the intervals $[0, r)$, $[r, r^{2})$ and $[r^{2}, \infty)$ separately, we get, for any $\alpha > d-2$ and $\beta \geq p\alpha/2$
  \begin{align*}
    &\int_{0}^{\infty} k_{\alpha}(t, y')\, \md t
    \\[.5ex]
    &\mspace{36mu}\leq\;
    \int_{0}^{r} r^{-d/2-\alpha/2}\, \md t
    +\,
    \int_{r}^{r^{2}} 
      t^{-d/2}\, \biggl(\frac{r}{\sqrt{t}}\biggr)^{\!\!-\alpha}\,
    \md t
    +\,
    \int_{r^{2}}^{\infty} t^{-d/2}\, \md t
    \;\leq\;
    c\, r^{-(d-2)},
  \end{align*}
  which concludes the proof of \eqref{eq:green:p}.

  (ii) By applying again Minkowski's integral inequality, we first get that
  \begin{align}
    \label{eq:gradient:1:split:1}
    &\mean\Biggl[
      \biggl| 
        \int_{0}^{\infty} \nabla^1_{\!x}\, p_{0, t}^{\omega}(0, y') 
      \biggr|^{2}
    \Biggr]^{1/2}\!\!
    \nonumber\\[.5ex]
    &\mspace{36mu}\leq\;
    \int_{0}^{r^{2}}
      \mean\Bigl[
        \bigl| \nabla^1_{\!x}\, p_{0,t}^{\omega}(0,y') \bigr|^{2}
      \Bigr]^{1/2}\,
    \md t
    \,+
    \int_{r^{2}}^{\infty}
      \mean\Bigl[
        \bigl| \nabla^1_{\!x}\, p_{0,t}^{\omega}(0,y') \bigr|^{2}
      \Bigr]^{1/2}\,
    \md t.
  \end{align}
  Thus, by Proposition~\ref{prop:gradient:offdiag}, we find that
  \begin{align*}
    \int_{r^{2}}^{\infty}
      \mean\Bigl[
        \bigl| \nabla^1_{\!x}\, p_{0,t}^{\omega}(0,y') \bigr|^{2}
      \Bigr]^{1/2}\,
    \md t
    \;\leq\;
    \int_{r^{2}}^{\infty} \const[C]{mean:gradient}\, t^{-(d+1)/2}\, \md t
    \;\leq\;
    c\, r^{-(d-1)}.
  \end{align*}
  Likewise, by splitting the domain of integration of the first integral on the right-hand side of \eqref{eq:gradient:1:split:1} into the contributions coming from the intervals $[0, r)$ and $[r, r^{2}]$ and using that $\alpha > 2(d-1)$, we obtain
  \begin{align*}
    &\int_{0}^{r^{2}}
      \mean\Bigl[
        \bigl| \nabla^1_{\!x}\, p_{0,t}^{\omega}(0,y') \bigr|^{2}
      \Bigr]^{1/2}\,
    \md t
    \nonumber\\[.5ex]
    &\mspace{36mu}\leq\;
    c\,
    \Biggl(
      \int_{0}^{r} (1 \vee t)^{-1/2}\, r^{-d/2-\alpha/4}\, \md t
      \,+
      \int_{r}^{r^{2}} 
        t^{-(d+1)/2} \biggl(\frac{r}{\sqrt{t}}\biggr)^{\!\!-\alpha/2}\,
      \md t
    \Biggr)
    \;\leq\;
    c\, r^{-(d-1)}. 
  \end{align*}
  Hence, by combining the estimates above, the assertion \eqref{eq:gradient:green:1} follows.
  
  (iii) First, an application of \eqref{eq:property:mixed:deriv} and the Cauchy-Schwarz inequality yields that
  \begin{align*}
    &\mean\Biggl[
      \biggl|
        \int_{0}^{\infty}
          \nabla^1_{\!-x}\, \nabla^{2}_{\!x'}\, p_{0, t}^{\omega}(0, y')\,
        \md t
      \biggr|
    \Biggr]
    \\[.5ex]
    &\mspace{36mu}\leq\;
    \int_{0}^{\infty}
      \sum_{z \in \mathbb{Z}^{d}}
      \mean\Bigl[
        \bigl| \nabla^1_{\!x}\, p_{0, t/2}^{\omega}(0,z) \bigr|^{2}
      \Bigr]^{1/2}
      \mean\Bigl[
        \bigl| \nabla^{2}_{\!x'}\, p_{t/2, t}^{\omega}(z,y') \bigr|^{2}
      \Bigr]^{1/2}\,
    \md t
    \\[.5ex]
    &\mspace{32mu}\overset{\!\!\eqref{eq:gradient:offdiag}\!\!}{\;\leq\;}
    \const[C]{mean:gradient}^{2}\,
    \int_{0}^{\infty}
      (1 \vee t)^{-1}
      \sum_{z \in \mathbb{Z}^{d}}
      k_{\alpha/2}(t/2,z)\, k_{\alpha/2}(t/2,y'-z)\,
    \md t 
    \\[.5ex]
    &\mspace{32mu}\overset{\!\!\eqref{eq:kernel:convolution}\!\!}{\;\leq\;}
    \const[C]{kernel:convolution} \const[C]{mean:gradient}^{2}\,
    \int_{0}^{\infty}
      (1 \vee t)^{-1}\, k_{\alpha/2}(t, y')\,
    \md t. 
  \end{align*}
  In order to estimate the remaining integral, we decompose the domain of integration into the contributions coming from the intervals $[0, r)$, $[r, r^{2})$ and $[r^{2}, \infty)$.  Since $\alpha > 2d$, we get that
  \begin{align*}
    &\int_{0}^{\infty}
      (1 \vee t)^{-1}\, k_{\alpha/2}(t, y')\,
    \md t
    \\[.5ex]
    &\mspace{36mu}\leq\;
    \int_{0}^{r} (1 \vee t)^{-1} r^{-d/2-\alpha/4}\, \md t
    \,+
    \int_{r}^{r^{2}} 
      t^{-d/2-1}\, \biggl(\frac{r}{\sqrt{t}}\biggr)^{\!\!-\alpha/2}\,
    \md t
    \,+
    \int_{r^{2}}^{\infty} t^{-d/2-1}\, \md t
    \\[1.5ex]
    &\mspace{36mu}\leq\;
    c\, r^{-d},
  \end{align*}
  which concludes the proof of \eqref{eq:gradient:green:12}. 
\end{proof}

\section{Annealed CLT and annealed local limit theorems}
\label{sec:alclt}
In this section, we discuss annealed CLT and establish annealed local limit theorems under mild conditions. 

Since properties of the environmental processes are crucial for the proof of the annealed CLT, we start with providing a proof of Lemma~\ref{lem:non-explosion}.
\begin{proof}[Proof of Lemma~\ref{lem:non-explosion}]
  (i) The assertion that, for $\prob$-a.e.\ $\omega$, both processes $X$ and $\tilde{X}$ do not explode $\Prob_{0, 0}^{\omega}$-a.s.\ follows directly from \cite[Lemma 4.1]{ACDS18}. In particular, for $\prob$-a.e.\ $\omega$, any $x \in \mathbb{Z}^{d}$ and $s \leq t$, we have that $1 = (P_{s, t}^{\omega} 1)(x)$ and $1 = (P_{-t, -s}^{\tilde{\omega}} 1)(x) = (P_{s, t}^{\omega} 1)^{*}(x)$ which implies \eqref{eq:p:sum}.

  (ii) Let $(\mathcal{F}_{t}^{X} : t \geq 0)$ be the natural filtration induced by $X$. Then, for $\prob$-a.e.\ $\omega$, any bounded, measurable $\varphi\colon \Omega \to \mathbb{R}$ and $s \leq t$ we obtain that, $\Prob_{0, 0}^{\omega}$-a.s.,
  \begin{align*}
    \Mean_{0, 0}^{\omega}\bigl[\varphi(\tau_{t, X_{t}} \omega) \mid \mathcal{F}_{s}^{X}\bigr]
    &\;=\;
    \Mean_{s, X_{s}}^{\omega}\bigl[\varphi(\tau_{t, X_{t}} \omega)\bigr]
    \\[.5ex]
    &\;=\;
    \sum_{y \in \mathbb{Z}^{d}} p_{0, t-s}^{\tau_{s, X_{s}}\omega}(0, y)\, \varphi(\tau_{t, X_{s} + y} \omega)
    \;=\;
    (\mathcal{P}_{t-s} \varphi)(\tau_{s, X_{s}} \omega),
  \end{align*}
  where $(\mathcal{P}_{t} \varphi)(\omega) \ldef \sum_{y \in \mathbb{Z}^{d}} p_{0, t}^{\omega}(0, y) \varphi(\tau_{t, y} \omega)$. Likewise, we obtain that, $\Prob_{0, 0}^{\tilde{\omega}}$-a.s.,
  \begin{align*}
    \Mean_{0, 0}^{\tilde{\omega}}\bigl[\varphi(\tau_{-t, \tilde{X}_{t}} \omega) \mid \mathcal{F}_{s}^{\tilde{X}}\bigr]
    &\;=\;
    \sum_{y \in \mathbb{Z}^{d}} p_{-(t-s), 0}^{\tau_{-s, \tilde{X}_{s}} \omega}(y, 0)\,
    \varphi(\tau_{-t, \tilde{X}_{s} + y} \omega)
    \;=\;
    (\tilde{\mathcal{P}}_{t-s} \varphi)(\tau_{-s, \tilde{X}_{s}} \omega),
  \end{align*}
  where $(\tilde{\mathcal{P}}_{t}\varphi)(\omega) \ldef \sum_{y \in \mathbb{Z}^{d}} p_{-t, 0}^{\omega}(y, 0) \varphi(\tau_{-t, y} \omega)$.
  This proves the Markov property for both the environment process and also identifies the corresponding transition probability semi-groups.

  Further, by using that $\prob$ is stationary with respect to time-space shift, it follows that, for any bounded, measurable $\varphi\colon \Omega \to \mathbb{R}$ and $t \geq 0$,
  \begin{align*}
    \mean\bigl[\mathcal{P}_{t} \varphi\bigr]
    \;=\;
    \sum_{y \in \mathbb{Z}^{d}}
    \mean\bigl[p_{0, t}^{\omega}(0, y) \varphi(\tau_{t, y} \omega)\bigr]
    \overset{\eqref{eq:shift:hk}}{\;=\;}
    \sum_{y \in \mathbb{Z}^{d}}
    \mean\bigl[p_{-t, 0}^{\omega}(-y, 0) \varphi(\omega)\bigr]
    \overset{\eqref{eq:p:sum}}{\;=\;}
    \mean[\varphi],
  \end{align*}
  and, likewise, $\mean[\tilde{\mathcal{P}}_{t} \varphi] = \mean[\varphi]$. Thus, $\prob$ is the invariant distribution for both environmental processes.
  
  Finally, we are left with proving that $\prob$ is ergodic, that is, $\prob[A] \in \{0, 1\}$ for any $A \in \mathcal{F}$ such that $\mathcal{P}_{t} \mathbbm{1}_{A} = \mathbbm{1}_{A}$ $\prob$-a.s.\ for any $t \geq 0$. Let $A \in \mathcal{F}$ be such a measurable subset satisfying $\prob[A] > 0$ (otherwise that assertion is immediate). Then, for $\prob$-a.e.\ $\omega$ and $t \geq 0$,
  \begin{align*}
    0
    \;=\;
    \mathbbm{1}_{A^{\mathrm{c}}}(\omega) (\mathcal{P}_{t} \mathbbm{1}_{A})(\omega)
    &\;=\;
    \sum_{z \in \mathbb{Z}^{d}} \mathbbm{1}_{A^{\mathrm{c}}}(\omega)\,
    p_{0, t}^{\omega}(0, z)\, \mathbbm{1}_{A}(\tau_{t, z} \omega).
  \end{align*}
  By dropping all but one summand, we find that, for $\prob$-a.e.\ $\omega$,
  \begin{align}
    \label{eq:indicator:hk:indicator:positive:time}
    0
    \;=\;
    \mathbbm{1}_{A^{\mathrm{c}}}(\omega)\, p_{0, t}^{\omega}(0, y)\, \mathbbm{1}_{A}(\tau_{t, y} \omega),
    \qquad \forall\, t \geq 0,\, y \in \mathbb{Z}^{d}.
  \end{align}
  Further, by replacing $\omega$ with $\tau_{-t, -y} \omega$, interchanging the role of $A$ and $A^{\mathrm{c}}$ and using \eqref{eq:shift:hk}, we conclude that, for $\prob$-a.e.\ $\omega$,
  \begin{align}
    \label{eq:indicator:hk:indicator:negative:time}
    0
    \;=\;
    \mathbbm{1}_{A^{\mathrm{c}}}(\omega)\, p_{-t, 0}^{\omega}(0, y)\, \mathbbm{1}_{A}(\tau_{-t, y} \omega),
    \qquad \forall\, t \geq 0,\, y \in \mathbb{Z}^{d}.
  \end{align}
  \indent
  We claim that the transition kernel is irreducible $\prob$-a.s., that is, for $\prob$-a.e.\ $\omega$,
  \begin{align}
    \label{eq:irreducibility}
    p_{s, t}^{\omega}(x, y) \;>\; 0
    \qquad \forall\, -\infty < s < t < \infty,\, x, y \in \mathbb{Z}^{d}.
  \end{align}
  In view of the Chapman-Kolmogorov equation, it suffices to prove the claim for any $-\infty < s < t < \infty$ and $x, y \in \mathbb{Z}^{d}$ such that $x = y$ or $\{x, y\} \in E_{d}$. By using the integrated backward equation \eqref{eq:backward:integrated} with $f = \mathbbm{1}_{y}$, we obtain that, for $\prob$-a.e.\ $\omega$,
  \begin{align*}
    \begin{split}
      p_{s, t}(y, y)
      &\;=\;
      (P_{s, t}^{\omega} \mathbbm{1}_{y})(y)
      \;\geq\;
      \me^{-\int_{s}^{t} \mu_{u}^{\omega}(y)\, \md u}
      \\[.5ex]
      p_{s, t}^{\omega}(x, y)
      &\;=\;
      \bigl(P_{s, t}^{\omega} \mathbbm{1}_{y}\bigr)(x)
      \;\geq\;
      \me^{-\int_{s}^{t} \mu_{u}^{\omega}(x)\, \md u}\,
      \me^{-\int_{s}^{t} \mu_{u}^{\omega}(y)\, \md u}\,
      \int_{s}^{t} \omega_{r}(x, y)\, \md r
    \end{split}
  \end{align*}
  Since, by assumption, $\prob\bigl[\omega_{t}(e) > 0\bigr] = 1$ for any $t \in \mathbb{R}$ and $e \in E_{d}$, Fubini's theorem implies that
  \begin{align}
    \label{eq:local:integrated:positivity}
    \prob\biggl[
      \int_{I} \omega_{s}(e)\, \mathrm{d}s > 0
    \biggr]
    \;=\;
    1,
    \qquad \text{for any finite interval } I \subset \mathbb{R},\,
    e \in E_{d}.
  \end{align}
  Thus, from \eqref{eq:local:integrability} and \eqref{eq:local:integrated:positivity} it follows that for $\prob$-a.e.\ $\omega$, any $s, t \in \mathbb{R}$ with $s < t$ and $e \in E_{d}$,
  \begin{align*}
    \int_{s}^{t} \mu_{r}^{\omega}(x)\, \md r \;<\; \infty,
    \quad \forall\, x \in \mathbb{Z}^{d}
    \qquad \text{and} \qquad
    \int_{s}^{t} \omega_{r}(e)\, \md r \;>\; 0,
    \quad \forall\, e \in E_{d}.
  \end{align*}
  This completes the proof of the claim \eqref{eq:irreducibility}.

  By combining \eqref{eq:indicator:hk:indicator:positive:time} and \eqref{eq:indicator:hk:indicator:negative:time} with \eqref{eq:irreducibility}, $\mathbbm{1}_{A^{\mathrm{c}}}(\omega)\, \mathbbm{1}_{A}(\tau_{t, y} \omega) = 0$ for $\prob$-a.e.\ $\omega$ and any $t \in \mathbb{R}$, $y \in \mathbb{Z}^{d}$. Now, \cite[Lemma~1.2.1]{CFS82} implies that there exists $B \in \mathcal{F}$ with $\tau_{t, y}^{-1}(B) = B$ for all $(t, y) \in \mathbb{R} \times \mathbb{Z}^{d}$ such that $\prob[A \triangle B] = 0$. Recalling that $\prob[A] > 0$, it is evident that $\prob[B] \geq \prob[A] - \prob[A \triangle B] > 0$. Since $\prob$ is ergodic with respect to time-space shifts it follows that $\prob[B] = 1$. This implies that $\prob[A] \geq \prob[B] - \prob[A \triangle B] = 1$, and consequently, $\prob[A] = 1$.

  By following literally the same argument it also follows that $\prob$ is ergodic with respect to $(\tau_{t, \tilde{X}_{t}} \omega : t \geq 0)$.
\end{proof}
Next, we give a proof of Proposition~\ref{prop:mean:displacement}.  For any $t \in \mathbb{R}$, let $T_{t}\colon L^{2}(\Omega, \prob) \to L^{2}(\Omega, \prob)$ be the map defined by $T_{t} \varphi \ldef \varphi \circ \tau_{t,0}$.  By Assumption~\ref{ass:law}-(ii), $\{T_{t} : t \in \mathbb{R}\}$ is a strongly continuous contraction semi-group on $L^{2}(\Omega, \prob)$, cf.~\cite[Section 7.1]{JKO94}.  Its infinitesimal generator $D_{0}\colon \dom(D_{0}) \subset L^{2}(\Omega, \prob) \to L^{2}(\Omega,\prob)$ is defined by
\begin{align}
  D_{0} \varphi
  \;\ldef\;
  \lim_{t \to 0} \frac{1}{t}\, \bigl(T_{t} \varphi - \varphi\bigr),
\end{align}
where $\dom(D_{0})$ denotes the set of all $\varphi \in L^{2}(\Omega, \prob)$ such that the limit above exists.  Notice that $\dom(D_{0})$ is dense in $L^{2}(\Omega, \prob)$, cf.~\cite[Theorem~34.4]{La02}.  Moreover, by Assumption~\ref{ass:law}-(i) it holds that $\langle D_{0} \varphi, \psi \rangle_{L^{2}(\prob)} = -\langle \varphi, D_{0} \psi\rangle_{L^{2}(\prob)}$ and $\langle D_{0} \varphi, \varphi \rangle_{L^{2}(\prob)} = 0$ for all $\varphi, \psi \in \dom(D_{0})$.

Consider the linear operator $\mathcal{L}\colon \dom(\mathcal{L}) \to L^{2}(\prob)$,
%
\begin{align}
  (\mathcal{L} \varphi)(\omega)
  \;=\;
  (D_{0} \varphi)(\omega)
  \,+\, \sum_{y \sim 0} \omega_{0}(0,y)\, 
  \bigl( \varphi(\tau_{0,y} \omega) - \varphi(\omega) \bigr),
\end{align}
and denote by $\mathcal{L}^{*}$ the $L^{2}(\prob)$-adjoint operator of $\mathcal{L}$. Notice that $\mathcal{C} \ldef \dom(D_{0}) \cap \mathbb{L}^{\infty}(\Omega, \prob)$ is a common core of $\mathcal{L}$ and $\mathcal{L}^{*}$.  We define a semi-norm, $\|\cdot\|_{H^1(\prob)}$, on $\mathcal{C}$ by $\|\varphi\|_{H^1(\prob)}^{2} \ldef \langle \varphi, -\mathcal{L} \varphi \rangle_{L^{2}(\prob)}$.  Again, Assumption~\ref{ass:law}-(i) ensures that
\begin{align*}
  \|\varphi\|_{H^1(\prob)}^{2}
  \;=\;
  \frac{1}{2}
  \sum_{x \sim 0}
  \mean\Bigl[\omega_{0}(0,x)\, \bigl(\varphi(\tau_x \omega) - \varphi(\omega)\bigr)^{2} \Bigr],
  \qquad \varphi \in \mathcal{C}.
\end{align*}
Let $H^{1}(\prob)$ be the completion of $\mathcal{C}$ with respect to $\|\cdot\|_{H^1(\prob)}$ taken modulo suitable equivalence classes of elements of $\mathcal{C}$.  We write $H^{-1}(\prob)$ to denote the dual of $H^1(\prob)$.  For $\psi \in L^{2}(\prob)$, let
\begin{align*}
  \|\psi\|_{H^{-1}(\prob)}^{2}
  \;\ldef\;
  \sup_{\varphi \in \mathcal{C}}
  \Bigl( 
    2 \langle \varphi, \psi \rangle_{L^{2}(\prob)} - \|\varphi\|_{H^1(\prob)}^{2} 
  \Bigr)
  \;\in\;
  [0, \infty].
\end{align*}
One may identify the Hilbert space $H^{-1}(\prob)$ with the completion of equivalence classes of element of $\{\varphi \in H^{1}(\prob) : \|\varphi\|_{H^{-1}(\prob)} < \infty\}$ with respect to $\| \cdot \|_{H^{-1}(\prob)}$.  For more details we refer to \cite[Section~2.2]{KLO12}.
\smallskip

Given all these notions, we first prove Proposition~\ref{prop:mean:displacement} that the maximal mean displacement of the process up to time $T > 0$ is of order $\sqrt{T}$.  
%
%
\begin{proof}[Proof of Proposition~\ref{prop:mean:displacement}]
  For any $i \in \{1, \ldots, d\}$, set $\mathbb{R} \times \mathbb{Z}^{d} \ni (t,x) \mapsto f_{i}(t,x) = x^i$, where $x^i$ is the $i$th component of $x$.  Since for $\prob$-a.e.\ $\omega$ the stochastic process $X$ does not explode $\Prob$-a.s., it follows that, for almost every $\omega$,
  \begin{align}\label{eq:Dynkin:martingale}
    M_{t}^i
    \;\ldef\;
    f_{i}(t, X_{t}) - f_{i}(0, X_{0})
    -
    \int_{0}^{t}
      (\partial_{s} f_{i})(s, X_{s}) 
      + \bigl(\mathcal{L}_{s}^{\omega} f_{i}(s, \cdot)\bigr)(X_{s})\;
    \md s
  \end{align}
  is a $(\Prob_{0,0}^{\omega}, \{\mathcal{F}_{t} : t \geq 0\})$-martingale with respect to the natural augmented filtration generated by the process $X$.  Note that the additive functional appearing on the right-hand side of \eqref{eq:Dynkin:martingale} can be expressed as
  \begin{align*}
    \int_{0}^{t}
      (\partial_{s} f_{i})(s, X_{s}) 
      + \bigl(\mathcal{L}_{s}^{\omega} f_{i}(s, \cdot)\bigr)(X_{s})\;
    \md s
    \;=\;
    \int_{0}^{t}
      \delta^i(\omega) \circ \tau_{s, X_{s}}\;
    \md s,
  \end{align*}
  where $\delta^{i}(\omega) = \sum_{x \sim 0} \omega_{0}(0,x) x^{i}$ is the so-called local drift.  Hence, by using Doob's $L^{p}$-inequality, we obtain
  \begin{align*}
    &\mean\biggl[
      \Mean_{0, 0}^{\omega}\biggl[
        \sup_{0 \leq t \leq T} \bigl| X_{t}^{i} \bigr|
      \biggr]
    \biggr]
    \\[.5ex]
    &\mspace{36mu}\leq\;
    \mean\bigl[
      \Mean_{0,0}^{\omega}\bigl[ \langle M^i \rangle_{T} \bigr] 
    \bigr]^{1/2}
    \,+\,
    \mean\biggl[
      \Mean_{0, 0}^{\omega}\biggl[
        \sup_{0 \leq t \leq T}
        \biggl| \int_{0}^{t} \delta^{i}(\omega) \circ \tau_{s, X_{s}}\, \md s \biggr|
      \biggr]
    \biggr],
  \end{align*}
  where the quadratic variation processes $(\langle M^{i} \rangle : t \geq 0)$ is given by
  \begin{align*}
    \langle M^{i} \rangle_{T}
    \;=\;
    \int_{0}^{T}
    \sum_{x \sim 0} \omega_{s}(X_{s}, X_{s} + x) (x^{i})^{2}\,
    \md s
    \;=\;
    \int_{0}^{T}
    \sum_{x \sim 0} \omega_{0}(0, x) (x^{i})^{2} \circ \tau_{s, X_{s}}\,
    \md s.
  \end{align*}
  By Fubini's theorem and the stationarity of the process as seen from the particle, it follows that $\mean\bigl[\Mean_{0,0}^{\omega}\bigl[\langle M^i \rangle_{T}\bigr]\bigr] = T \sum_{x \sim 0}\mean\bigl[\omega_{0}(0,x) (x^i)^{2}\bigr]$.  In order to derive an upper bound on the additive functional, set $\delta_{K}^{i}(\omega) \ldef \sum_{x \sim 0} (\omega_{0}(0,x) \wedge K)\, x^i$.  Obviously, $\delta_{K}^{i} \in L^{2}(\prob)$ for any $K < \infty$.  Moreover, in view of Assumption~\ref{ass:law}(i), we get, for any $\varphi \in H^1(\prob)$,
  \begin{align*}
    \langle \delta_{K}^{i}, \varphi \rangle_{L^{2}(\prob)}
    &\;=\;
    \frac{1}{2}\, \sum_{x \sim 0}
    \mean\bigl[(\omega_{0}(0,x) \wedge K) (\varphi(\tau_x \omega) - \varphi(\omega)) x^i\bigr].
    \\
    &\;\leq\;
    \| \varphi \|_{H^{1}(\prob)}\,
    \biggl(
      \frac{1}{2} \sum_{x \sim 0} \mean\bigl[\omega_{0}(0,x) (x^i)^{2}\bigr]
    \biggr)^{\!\!1/2}.
  \end{align*}
  Hence, $\delta_{K}^{i} \in H^{-1}(\prob) \cap L^{2}(\prob)$ with $\|\delta_{K}^{i}\|_{H^{-1}(\prob)}^{2} \leq \frac{1}{2}\sum_{x \sim 0} \mean\bigl[\omega_{0}(0,x) (x^i)^{2} \bigr]$.  Thus, by \cite[Lemma~2.4]{KLO12}, we obtain that
  \begin{align*}
    \mean\biggl[
      \Mean_{0, 0}^{\omega}\biggl[
        \sup_{0 \leq t \leq T}
        \biggl( 
          \int_{0}^{t} \delta_{K}^{i}(\omega) \circ \tau_{s, X_{s}}\, \md s 
        \biggr)^{\!\!2}
      \biggr]
    \biggr]
    \;\leq\;
    24 T\, \|\delta_{K}^{i}\|_{H^{-1}(\prob)}^{2}. 
  \end{align*}
  Hence,
  \begin{align*}
    &\mean\biggl[
      \Mean_{0, 0}^{\omega}\biggl[
        \sup_{0 \leq t \leq T} \bigl|X_{t}^i\bigr|
      \biggr]
    \biggr]
    \\[.5ex]
    &\mspace{36mu}\leq\;
    c\, \sqrt{T}\,
    \biggl(
      \sum_{x \sim 0}\, \mean\bigl[\omega_{0}(0,x) (x^i)^{2}\bigr]
    \biggr)^{\!1/2}
    + T\, \mean\bigl[|\delta_{K}^{i}(\omega) - \delta^{i}(\omega)|\bigr]. 
  \end{align*}
  Since, $|\delta_{K}^{i}(\omega) - \delta^{i}(\omega)| \leq 2 |\delta^{i}(\omega)|$ for any $K > 0$, $\delta^{i} \in L^{1}(\prob)$, and $\delta^{i}_{K}(\omega) \to \delta^{i}(\omega)$ as $K \to \infty$ for $\prob$-a.e.~$\omega$, the assertion \eqref{eq:mean:displacement} follows from Lebesgue's dominated convergence theorem.  Finally, for any $t > 0$, the tightness of the family $\{X_{t}^{(n)} : n \in \mathbb{N}\}$ follows immediately from \eqref{eq:mean:displacement} by Markov's inequality.
\end{proof}
In order to prove the annealed CLT, we follow the common approach for proving the corresponding quenched result by introducing the so-called \emph{harmonic} coordinates, $\Phi\colon \Omega \times \mathbb{R} \times \mathbb{Z}^{d} \to \mathbb{R}^{d}$ and the \emph{corrector}, $\chi\colon \Omega \times \mathbb{R} \times \mathbb{Z}^{d} \to \mathbb{R}^{d}$.  The latter solves component-wise the time-inhomogeneous Poisson equation
\begin{align*}
  \partial_{t} u + \mathcal{L}_{t}^{\omega} u \;=\; \mathcal{L}_{t}^{\omega} \Pi,
\end{align*}
where $\Pi$ denotes the identity map on $\mathbb{Z}^{d}$. As a consequence, the harmonic coordinates defined by
\begin{align}\label{eq:corrector:def}
  \Phi(\omega, t, x)
  \;=\;
  x - \chi(\omega, t, x)
\end{align}
is component-wise a time-space harmonic function.  By inspecting the proofs of \cite[Theorem~2.5 and Corollary~2.7]{ACDS18}, the existence of $\Phi$ and $\chi$ are ensured if Assumption~\ref{ass:law} and \ref{ass:integrability} holds true.  Moreover, in view of \cite[Proposition~3.3]{ACDS18}, we also know that, for any $v \in \mathbb{R}^{d}$, $\delta > 0$, $t \in (0, \infty)$ and $K \in (0, \infty)$,
\begin{align}\label{eq:sublinearity:average}
  \lim_{n \to \infty} \frac{1}{n^{d+2}}
  \int_{0}^{t n^{2}}\mspace{-12mu}
    \sum_{x \in B(0, Kn)}\mspace{-6mu} \indicator_{|v \cdot \chi(\omega, s, x)| > \delta n}\,
  \md s
  \;=\;
  0
  \qquad \prob\text{-a.s.\ and in } L^{1}(\prob).
\end{align}
\begin{proof}[Proof of Proposition~\ref{prop:annealed:CLT}]
  In order to lighten notations we set $\Phi_{t} \ldef \Phi(\omega, t, X_{t})$ and $\chi_{t} \ldef \chi(\omega, t, X_{t})$ for any $t \geq 0$.  Obviously, \eqref{eq:corrector:def} implies that $X_{t} = \Phi_{t} + \chi_{t}$ for any $t \geq 0$.  By \cite[Corollary~4.2]{ACDS18}, we know that, for any $v \in \mathbb{R}^{d}$, the stochastic process $(v \cdot \Phi_{t} : t \geq 0)$ is a $\Prob_{0, 0}^{\omega}$-martingale for $\prob$-a.e.\ $\omega$ with $v \cdot \Phi_{0} = 0$ whose quadratic variation process is given by
  \begin{align*}
    \langle v \cdot \Phi \rangle_{t}
    \;=\;
    \int_{0}^{t}
      \sum_{y}
      \omega_{0}(0, y) \bigl(v \cdot \Phi(\omega, 0, y)\bigr)^{2} \circ \tau_{s, X_{s}}\,
    \md s.
  \end{align*}
  In particular, by using the Cauchy-Schwarz inequality and the stationary of the process as seen from the particle, we obtain that, for any $0 \leq s \leq t < \infty$,
  \begin{align}\label{eq:harmonic:process:incr}
    \mean\bigl[
      \Mean_{0, 0}^{\omega}\bigl[|v \cdot \Phi_{t} - v \cdot \Phi_{s}|\bigr]
    \bigr]^{2}
    &\;\leq\;
    \mean\bigl[
      \Mean_{0, 0}^{\omega}\bigl[
        \langle v \cdot \Phi \rangle_{t} - \langle v \cdot \Phi \rangle_{s}
      \bigr]
    \bigr]
    \;=\;
    \bigl(v \cdot \Sigma^{2} v\bigr) (t-s),  
  \end{align}
  where
  \begin{align*}
    v \cdot \Sigma^{2}v
    \;\ldef\;
    \mean\biggl[
      \sum\nolimits_{y \in \mathbb{Z}^{d}} \omega_{0} \bigl(v \cdot \Phi(\omega, 0, y)\bigr)^{2}
    \biggr].
  \end{align*}
  Furthermore, from \cite[Proposition~4.4]{ACDS18}, we know that the diffusively rescaled martingale, $(\Phi_{t}^{(n)} : t \geq 0)$, defined as $\Phi_{t}^{(n)} \ldef n^{-1} \Phi_{t n^{2}}$ for any $t \geq 0$, converges for $\prob$-a.e.\ $\omega$ under $\Prob_{0, 0}^{\omega}$ to a Brownian motion $(W_{t} : t \geq 0)$ on $\mathbb{R}^{d}$ with a deterministic and non-degenerate covariance matrix $\Sigma^{2}$. In particular, by employing Lebesgue's dominated convergence theorem, for any $k \in \mathbb{N}$ and $0 \leq t_{1} < t_{2} < \cdots < t_{k} < \infty$, we deduce that the distribution of $(\Phi_{t_{1}}^{(n)}, \ldots, \Phi_{t_{k}}^{(n)})$ under $\prob_{0, 0}^{*}$ converges weakly to the distribution of $(W_{t_{1}}, \ldots, W_{t_{k}})$.  Therefore, by Slutzki's theorem, it suffices to show that, for any $\delta > 0$
  \begin{align}\label{eq:convergence:corrector:process:max}
    \lim_{n \to \infty}
    \mean\biggl[
      \Prob_{0, 0}^{\omega}\biggl[
        \sup_{1 \leq i \leq k} \bigl| \chi_{t_{i}}^{(n)} \bigr| > \delta
      \biggr]
    \biggr]
    \;=\;
    0
  \end{align}
  in order to conclude that the distribution of $(X_{t_{1}}^{(n)}, \ldots, X_{t_{k}}^{(n)})$ under $\prob_{0, 0}^{*}$ converges to $(W_{t_{1}}, \ldots, W_{t_{k}})$, thereby completing the proof of Proposition~\ref{prop:annealed:CLT}.

  Thus, it remains to show \eqref{eq:convergence:corrector:process:max}. By using a union bound and applying the Cram\'{e}r-Wold device, see e.g.\ \cite[Theorem~3.10.6]{Du19}, \eqref{eq:convergence:corrector:process:max} will follow once we have shown that, for any $t > 0$ and $v \in \mathbb{R}^{d}$
  \begin{align}\label{eq:convergence:corrector:process}
    \lim_{n \to \infty}
    \mean\Bigl[
      \Prob_{0, 0}^{\omega}\Bigl[
        \Bigl| n^{-1}v \cdot \chi_{t n^{2}} \Bigr| > \delta
      \Bigr]
    \Bigr]
    \;=\;
    0
  \end{align}
  For this purpose, fix some $\varepsilon > 0$.  Further, set $K \ldef 3 \const[C]{mean:displacement} \sqrt{t}/\varepsilon$ and choose some $\tau \in (0, t \wedge (\delta \varepsilon)^{2}/(16 \const[C]{heatkernel:ub}^{2}))$.  Then,
  \begin{align}\label{eq:corrector:process:split1}
    &\mean\Bigl[
      \Prob_{0, 0}^{\omega}\Bigl[
        \bigl|n^{-1} v \cdot \chi_{t n^{2}}\bigr| > \delta
      \Bigr]
    \Bigr]
    \nonumber\\
    &\mspace{36mu}\leq\;
    \mean\biggl[
      \Prob_{0, 0}^{\omega}\biggl[
        \max_{0 \leq s \leq tn^{2}} |X_{s}| > Kn
      \biggr]
    \biggr]
    \,+\,
    \mean\Bigl[
      \Prob_{0, 0}^{\omega}\Bigl[
        \bigl|n^{-1} v \cdot\chi_{t n^{2}}\bigr| > \delta, T_{Kn} > tn^{2}
      \Bigr]
    \Bigr]
    \nonumber\\
    &\mspace{32mu}%
    \overset{\mspace{-9mu}\eqref{eq:mean:displacement}\mspace{-9mu}}{\;\leq\;}
    \frac{\const[C]{mean:displacement} \sqrt{t}}{K}
    \,+\,
    \mean\Bigl[
      \Prob_{0, 0}^{\omega}\Bigl[
        \bigl|n^{-1} v \cdot\chi_{t n^{2}}\bigr| > \delta, T_{Kn} > tn^{2}
      \Bigr]
    \Bigr],
  \end{align}
  where we recall that $T_{K n} = \inf\{t \geq 0 : X_{t} \not\in B(0, Kn)\}$ denotes the first exit time from the ball $B(0, Kn)$.  It is evident that $\tau \in (0, t)$. Therefore, we can bound the second term in \eqref{eq:corrector:process:split1} by
  \begin{align}\label{eq:corrector:process:split2}
    &\mean\Bigl[
      \Prob_{0, 0}^{\omega}\Bigl[
        \bigl|n^{-1} v \cdot \chi_{tn^{2}} \bigr| > \delta, T_{Kn} > tn^{2}
      \Bigr]
    \Bigr]
    \nonumber\\
    &\mspace{36mu}\leq\;
    \frac{2}{\delta \tau n}  
    \int_{t - \tau}^{t}
      \mean\Bigl[
        \Mean_{0, 0}^{\omega}\Bigl[
          |v \cdot \chi_{t n^{2}} - v \cdot \chi_{s n^{2}}|
        \Bigr]
      \Bigr]\,
    \md s
    \nonumber\\
    &\mspace{72mu}+\,
    \frac{1}{\tau}
    \int_{t - \tau}^{t}
      \mean\Bigl[
        \Prob_{0, 0}^{\omega}\Bigl[
          \bigl| v \cdot \chi_{s n^{2}}\bigr| > \delta n/2, T_{Kn} > sn^{2}
        \Bigr]
      \Bigr]\,
    \md s.
  \end{align}
  Since, for any $0 \leq s \leq t < \infty$, an application of the Markov property, \eqref{eq:shift:hk} and the stationarity of the process as seen from the particle yields that
  \begin{align*}
    \mean\Bigl[
      \Mean_{0, 0}^{\omega}\Bigl[ \bigl|v \cdot X_{t} - v \cdot X_{s}\bigr| \Bigr]
    \Bigr]
    \;=\;
    \mean\Bigl[
      \Mean_{0, 0}^{\omega}\Bigl[ \bigl|v \cdot X_{t-s}\bigr| \Bigr]
    \Bigr]  
    \overset{\mspace{-9mu}\eqref{eq:mean:displacement}\mspace{-9mu}}{\;\leq\;}
    \const[C]{mean:displacement} \sqrt{t-s}  
  \end{align*}
  we get that
  \begin{align*}
    &\mean\Bigl[
      \Mean_{0, 0}^{\omega}\Bigl[ \bigl|v \cdot \chi_{t} - v \cdot \chi_{s}\bigr| \Bigr]
    \Bigr]
    \nonumber\\
    &\mspace{36mu}\leq\;
    \mean\Bigl[
      \Mean_{0, 0}^{\omega}\Bigl[ \bigl|v \cdot X_{t} - v \cdot X_{s}\bigr| \Bigr]
    \Bigr]
    +
    \mean\Bigl[
      \Mean_{0, 0}^{\omega}\Bigl[ \bigl|v \cdot \Phi_{t} - v \cdot \Phi_{s}\bigr| \Bigr]
    \Bigr]
    \overset{\mspace{-9mu}\eqref{eq:harmonic:process:incr}\mspace{-9mu}}{\;\leq\;}  
    c_{1} \sqrt{t-s},
  \end{align*}
  where $c_{1} \ldef \const[C]{mean:displacement} + \sqrt{v \cdot \Sigma^{2}v}$.  Therefore, the first term in \eqref{eq:corrector:process:split2} is bounded from above by $4 c_{1} \sqrt{\tau}/(3 \delta)$.  On the other hand, by using Lemma~\ref{lem:Nashclas}, we obtain that, for any $s \in [t-\tau, t]$,
  \begin{align*}
    &\mean\Bigl[
      \Prob_{0, 0}^{\omega}\Bigl[
        \bigl| v \cdot \chi_{s n^{2}}\bigr| > \delta n/2, T_{Kn} > sn^{2}
      \Bigr]
    \Bigr]
    \\
    &\mspace{32mu}\overset{\mspace{-9mu}\eqref{eq:nash}\mspace{-9mu}}{\;\leq\;}
    \frac{\const[C]{heatkernel:ub}}{(t - \tau)^{d/2} n^{d}} \sum_{x \in B(0, Kn)}\mspace{-6mu}
    \prob\Bigl[\bigl|v \cdot \chi(\omega, s n^{2}, x)\bigr| > \delta n/2\Bigr].  
  \end{align*}
   Notice that, in view of \eqref{eq:sublinearity:average}, there exists $n_{0} \equiv n_{0}(\delta, \varepsilon, t, \tau) \in \mathbb{N}$ such that
  \begin{align*}
    \frac{1}{n^{d+2}}
    \int_{0}^{t n^{2}}\mspace{-12mu}
      \sum_{x \in B(0, Kn)} \mspace{-6mu}
      \prob\Bigl[ \bigl|v \cdot \chi(\omega, s, x) \bigr| > \delta n/2 \Bigr]\,
    \md s
    \;\leq\;
    \frac{\varepsilon \tau (t - \tau)^{d/2}}{3 \const[C]{heatkernel:ub}}
    \qquad \forall\, n \geq n_{0}.
  \end{align*}
  Thus, by combining the above estimates, we finally obtain that, for any $n \geq n_{0}$,
  \begin{align*}
    \mean\Bigl[
      \Prob_{0, 0}^{\omega}\Bigl[
        \bigl|n^{-1} v \cdot \chi_{t n^{2}}\bigr| > \delta
      \Bigr]
    \Bigr]
    \;\leq\;
    \frac{\const[C]{mean:displacement} \sqrt{t}}{K} \,+\, \frac{4c_{1} \sqrt{\tau}}{3 \delta}
    \,+\, \frac{\varepsilon}{3}
    \;\leq\;
    \varepsilon.
  \end{align*}
  Since $\varepsilon > 0$ was arbitrarily chosen, the claim \eqref{eq:convergence:corrector:process} follows.
\end{proof}
\begin{proof}[Proof of Theorem~\ref{thm:annealed:LCLT}]
  (i) Let us start with proving the assertion \eqref{eq:alclt:time:constant}. From \cite[Theorem~3.1]{ACS21}, it immediately follows that, for any $0 < T_{0} < T_{1} < \infty$ and $K \subset \mathbb{R}^{d}$ compact,
  \begin{align}\label{eq:alclt:local}
    \lim_{n \to \infty} \sup_{y \in K} \sup_{t \in [T_{0}, T_{1}]}
    \bigl|
      n^{d} \bar{p}_{0, t n^{2}}(0, [y n]) - k_{t}^{\Sigma}(y)
    \bigr|
    \;=\;
    0
  \end{align}
  once we have verified the conditions $(\mathcal{G}.1)$ -- $(\mathcal{G}.4)$ as specified in Section~3 of \cite{ACS21}. Since the underlying graph we are considering is the Euclidean lattice $\mathbb{Z}^{d}$, $(\mathcal{G}.1)$ and $(\mathcal{G}.2)$ are obviously satisfied, whereas $(\mathcal{G}.3)$ follows from the annealed CLT as shown in Proposition~\ref{prop:annealed:CLT}.  To verify condition $(\mathcal{G}.4)$, we rely on the additional Assumption~\ref{ass:time:derivative}. Indeed, by using the estimates  \eqref{eq:gradient:on-diag} together with the Cauchy-Schwarz inequality and \eqref{eq:time:derivative}, we obtain that, for any $\delta \in (0, 1)$, $t \geq 4 \delta^{2} \vee T_{0}$ and $y \in \mathbb{R}^{d}$,
  \begin{align*}
    n^{d} \bigl| \bar{p}_{0, s n^{2}}(0, z) - \bar{p}_{0, s' n^{2}}(0, z') \bigr|
    &\;\leq\;
    \const[C]{gradient:point}\, \frac{|z - z'|}{n s^{(d+1)/2}}
    + \const[C]{time:elliptic}\, \frac{s'-s}{s^{d/2 + 1}}
    \\[.5ex]
    &\;\leq\;
    \biggl(
      \frac{2\, \const[C]{gradient:point}\, \delta / \sqrt{T_{0}}}{(1 - \delta^{2} / t)^{(d+1)/2}}
      +
      \frac{\const[C]{time:elliptic}\, \delta^{2} / T_{0}}{(1 - \delta^{2} / t)^{d/2 + 1}}
    \biggr)\, t^{-d/2}
    \\[.5ex]
    &\;\leq\;
    c\, \biggl(
      \frac{\delta}{\sqrt{T_{0}}} + \frac{\delta^{2}}{T_{0}} 
    \biggr)\, t^{-d/2}
  \end{align*}
  for every $z, z' \in B([y n], \delta n)$ and $t - \delta^{2} < s \leq s' \leq t$ which implies $(\mathcal{G}.4)$. Since
  \begin{align*}
    \lim_{t \to \infty} \sup_{y \in K} 
    \bigl|
      n^{d+1} \bar{p}_{0, t n^{2}}(0, [y n])
    \bigr|
    \;=\;
    0
    \qquad \text{and} \qquad
    \lim_{t \to \infty} \sup_{y \in K}
    \bigl| k_{t}^{\Sigma}(y)\bigr|
    \;=\;
    0,
  \end{align*}
  \eqref{eq:alclt:local} can be extended to hold uniformly on $[T_{0}, \infty)$ for any $T_{0} > 0$, thereby completing the proof of \eqref{eq:alclt:time:constant}.

  Next, we address the assertion~\eqref{eq:alclt}. For any $t > 0$ and $\delta \in (0, 1)$, the estimate given in \eqref{eq:gradient:on-diag}, combined with the Cauchy-Schwarz inequality, yields
  \begin{align*}
    \sup_{z, z' \in B([y n], \delta n)} n^{d} \bigl| \bar{p}_{0, t n^{2}}(0, z) - \bar{p}_{0, t n^{2}}(0, z') \bigr| 
    \;\leq\;
    2\, \const[C]{gradient:point}\, \delta\, t^{-(d+1)/2}.
  \end{align*}
  Thus, by employing this estimate in place of $(\mathcal{G}.4)$ and following the arguments presented in the proof of \cite[Theorem~3.1]{ACS21}, we conclude that for any $t > 0$ and any compact set $K \subset \mathbb{R}^{d}$
  \begin{align*}
    \lim_{n \to \infty} \sup_{y \in K}
    \bigl| n^{d} \bar{p}_{0, t n^{2}}(0, [y n]) - k_{t}^{\Sigma}(y) \bigr|
    \;=\;
    0.
  \end{align*}

  (ii) \textsc{Step 1.} Fix some unit vector $e_{i} \sim 0$, $i \in \{1, \ldots, d\}$.  We will first show that, for any $t > 0$ and $y \in \mathbb{R}^{d}$,
  \begin{align}\label{eq:gradient:lclt:pointwise}
    \lim_{n \to \infty}
    \bigl|
      n^{d+1} \nabla_{e_{i}}^{2} \bar{p}_{0, t n^{2}}(0, [y n])
      \,-\,
      \bigl(\partial_{i} k_{t}^{\Sigma}\bigr)(y)
    \bigr|
    \;=\;
    0.
  \end{align} 
  For any $\delta > 0$ let $f_{\delta} \in C^{2}(\mathbb{R}^{d}, \mathbb{R})$ be a non-negative function such that $\int_{\mathbb{R}^{d}} f_{\delta}(z)\, \md z = 1$ and $\supp f_{\delta} \subset C(y, \delta) \ldef y + (-\delta, \delta)^{d}$.  Then, for any $n \in \mathbb{N}$, by using both a summation and an integration by parts, we rewrite the left-hand side of \eqref{eq:gradient:lclt:pointwise} as
  \begin{align*}
    \bigl| 
      n^{d+1} \nabla_{e_{i}}^{2} \bar{p}_{0, t n^{2}}(0, [y n])
      \,-\,
      \bigl(\partial_{i} k_{t}^{\Sigma}\bigr)(y)
    \bigr|
    \;=\;
    \biggl| \sum_{i=1}^4 J_{i}(\delta, n) + J_{5}(\delta) \biggr|
  \end{align*}
  where
  \begin{align*}
    J_{1}(\delta, n)
    &\;\ldef\;
    \biggl(1 - \frac{1}{n^{d}} \sum_{z \in \mathbb{Z}^{d}} f_{\delta}(z/n)\biggr)\,
    n^{d+1} \nabla_{e_{i}}^{2} \bar{p}_{0, tn^{2}}(0, [y n]),
    \\
    J_{2}(\delta, n)
    &\;\ldef\;
    \sum_{z \in \mathbb{Z}^{d}} f_{\delta}(z/n)\,
    \Bigl(
      n \nabla_{e_{i}}^{2} \bar{p}_{0, tn^{2}}(0, [y n])
      -
      n \nabla_{e_{i}}^{2} \bar{p}_{0, tn^{2}}(0, z)
    \Bigr),
    \\
    J_{3}(\delta, n)
    &\;\ldef\;
    \sum_{z \in \mathbb{Z}^{d}}
    \Bigl(
      n \nabla_{-e_{i}} f_{\delta}(z/n) - \bigl(-\partial_{i} f_{\delta}\bigr)(z/n)
    \Bigr)\,
    \bar{p}_{0, t n^{2}}(0, z),
    \\
    J_{4}(\delta, n)
    &\;\ldef\;
    \mean\Bigl[
      \Mean_{0, 0}^{\omega}\Bigl[
        (-\partial_{i} f_{\delta})\bigl(X_{t}^{(n)}\bigr)
      \Bigr]
    \Bigr]
    \,-\,
    \int_{\mathbb{R}^{d}} (-\partial_{i} f_{\delta})(z)\, k_{t}^{\Sigma}(z)\, \md z,
    \\
    J_{5}(\delta)
    &\;\ldef\;
    \int_{\mathbb{R}^{d}} 
      f_{\delta}(z)\, 
      \Bigl( 
        \bigl(\partial_{i} k_{t}^{\Sigma}\bigr)(z)
        -
        \bigl(\partial_{i} k_{t}^{\Sigma}\bigr)(y)
      \Bigr)\,
    \md z.
  \end{align*}
  Thus, in order to prove \eqref{eq:gradient:lclt:pointwise} it suffices to show that
  \begin{align*}
    \limsup_{\delta \searrow 0} \limsup_{n \to \infty} |J_{i}(\delta, n)| 
    \;=\; 
    0
    \quad \forall\, i \in \{1, \ldots, 4\}
    \qquad \text{and} \qquad
    \limsup_{\delta \searrow 0} |J_{5}(\delta)| \;=\; 0.
  \end{align*}
  First, from Theorem~\ref{thm:gradient:est}-(i) is follows that
  \begin{align}\label{eq:gradient:rescaled}
    \bigl|
      n^{d+1} \nabla_{e_{i}}^{2} \bar{p}_{0, t n^{2}}(0,[y n])
    \bigr| 
    \;\leq\; 
    \const[C]{gradient:point}\, t^{-(d+1)/2},
    \qquad \forall\, n \in \mathbb{N}.
  \end{align}
  Therefore, the convergence of $|J_{1}(\delta, n)|$ to zero $n$ tends to infinity is an immediate consequence of the approximation of the integral $1 = \int_{\mathbb{R}^{d}} f_{\delta}(z)\, \md z$ by a $d$-dimensional Riemann sum.  In order to handle the summand involving $|J_{2}|$, notice that an application of Theorem~\ref{thm:gradient:est}-(iii) yields
  \begin{align}\label{eq:gradient:22:uniform}
    \max_{(0, x), (0, x') \in E_{d}} \sup_{z, z' \in \mathbb{Z}^{d}}
    \bigl|
      n^{d+2} \nabla_{x}^{2} \nabla_{x'}^{2} \bar{p}_{0, t n^{2}}(z, z')
    \bigr| 
    \;\leq\; 
    \const[C]{gradient:sum}\, t^{-(d+2)/2}
    \qquad \forall\, n \in \mathbb{N}.
  \end{align}
  This estimate implies that, for any $z \in n C(y, \delta) \cap \mathbb{Z}^{d}$,
  \begin{align}\label{eq:gradient:22:difference}
    n^{d+1}
    \bigl|
      \nabla_{e_{i}}^{2} \bar{p}_{0, tn^{2}}(0, [y n])
      -
      \nabla_{e_{i}}^{2} \bar{p}_{0, tn^{2}}(0, z)
    \bigr|
    \overset{\eqref{eq:gradient:22:uniform}}{\;\leq\;}
    \const[C]{gradient:sum}\, t^{-(d+2)/2}\, \biggl|\frac{[y n] - z}{n}\biggr|_{\infty}.
  \end{align}
  Hence,
  \begin{align*}
    |J_{2}(\delta, n)|
    \;\leq\;
    \const[C]{gradient:sum}\, t^{-(d+2)/2}\, \biggl(\delta + \frac{2}{n}\biggr)\,
    \frac{1}{n^{d}} \sum_{z \in \mathbb{Z}^{d}} f_{\delta}(z).
  \end{align*}
  This estimate clearly implies that $\limsup_{\delta \downarrow 0} \limsup_{n \to \infty} |J_{2}(\delta, n)| = 0$.  Next, we consider the term $|J_{3}|$.  Since $f_{\delta}$ is assumed to be twice continuously differentiable, we get that $n \bigl| n \nabla_{-e_{i}} f_{\delta}(z/n) - \bigl(-\partial_{i} f_{\delta}\bigr)(z/n) \bigr| \leq \sup_{x \in \mathbb{R}^{d}} |(\partial_{ii}^{2} f_{\delta})(x)| < \infty$ for any $z \in \mathbb{R}^{d}$.  Since $\sum_{z \in \mathbb{Z}^{d}} \bar{p}_{0, tn^{2}}(0, z) = 1$, we therefore obtain that, for any $\delta > 0$,
  \begin{align*}
    |J_{3}(\delta, n)|
    \;\leq\;
    \frac{1}{n}\, \sup_{x \in \mathbb{R}^{d}} |(\partial_{ii}^{2} f_{\delta})(x)|
    \underset{n \to \infty}{\;\longrightarrow\;}
    0.
  \end{align*}
  Further, from the specifying properties of the function $f_{\delta}$ it follows that $-\partial_{i} f_{\delta}$ is bounded and continuous.  Hence, the convergence $\lim_{n \to \infty}|J_{4}(\delta, n)| = 0$ for any $\delta > 0$ is an immediate consequence of the annealed CLT as stated in Proposition~\ref{prop:annealed:CLT}.  Finally, since $x \mapsto k_{t}^{\Sigma}(x)$ is a smooth function that vanishes at infinity, we obtain that
  \begin{align*}
    \delta^{-1}\, 
    \bigl|
      \bigl( \partial_{i} k_{t}^{\Sigma}\bigr)(z) 
      - 
      \bigl(\partial_{i} k_{t}^{\Sigma}\bigr)(y)
    \bigr| 
    \;\leq\; 
    \sup_{x \in \mathbb{R}^{d}} 
    \bigl|\bigl(\partial_{ii}^{2} k_{t}^{\Sigma}\bigr)(x)\bigr| 
    \;<\; 
    \infty
    \qquad \forall\, z \in C(y, \delta),
  \end{align*}
  which implies that $\lim_{\delta \downarrow 0} |J_{5}(\delta)| = 0$.  This completes the proof of \eqref{eq:gradient:lclt:pointwise}.

  \textsc{Step 2.} To extend the pointwise convergence stated in \eqref{eq:gradient:lclt:pointwise} to uniform convergence on a compact subset $K \subset \mathbb{R}^{d}$, we use a covering argument akin to that in \cite[Proof of Theorem~3.1 -- \emph{Step~2}]{ACS21}.  For $\eta \in (0, 1)$, we define $K_{\eta} \ldef K \cap \eta \mathbb{Z}^{d}$. It is evident that $K_{\eta} \neq \emptyset$ for sufficiently small $\eta$. Utilising the fact that the number of elements in $K_{\eta}$ is finite, \emph{Step~1} implies that for any $\varepsilon > 0$ there exists $\const[N]{N:aclt:gradient}(\varepsilon, \eta) \in \mathbb{N}$ such that, for any $t > 0$ and $i \in \{1, \ldots, d\}$,
  \begin{align*}
    \sup_{y \in K_{\eta}}
    \bigl|
      n^{d+1} \nabla_{e_{i}}^{2} \bar{p}_{0, t n^{2}}(0, [y n])
      \,-\,
      \bigl(\partial_{i} k_{t}^{\Sigma}\bigr)(y)
    \bigr|
    \;\leq\;
    \varepsilon,
    \qquad \forall\, n \geq \const[N]{N:aclt:gradient}(\varepsilon, \eta).
  \end{align*}
  Furthermore, for any $y \in K \setminus K_{\eta}$, let $y_{0}$ be a closest point in $K_{\eta}$. Then, we have $|y - y_{0}|_{\infty} \leq \eta$. In particular, it follows that $|(\partial_{i} k_{t}^{\Sigma})(y) - (\partial_{i} k_{t}^{\Sigma})(y_{0})| \leq c \eta t^{-(d+2)/2}$, and by using \eqref{eq:gradient:22:difference}, we obtain
  \begin{align*}
    n^{d+1} \bigl|
      \nabla_{e_{i}}^{2} \bar{p}_{0, t n^{2}}(0, [y n]) - \nabla_{e_{i}}^{2} \bar{p}_{0, t n^{2}}(0, [y_{0} n])
    \bigr|
    \;\leq\;
    c\, \eta\,t^{-(d+2)/2}.
  \end{align*}
  Thus, by choosing $\eta \ldef \varepsilon t^{(d+2)/2} / (2 c)$, we obtain that for any $n \geq \const[N]{N:aclt:gradient}(\varepsilon, \eta)$,
  \begin{align*}
    \sup_{y \in K}\,
    \bigl|
      n^{d+1} \nabla_{e_{i}}^{2} \bar{p}_{0, t n^{2}}(0, [y n])
      - \bigl(\partial_{i} k_{t}^{\Sigma}\bigr)(y)
    \bigr|
    \;\leq\;
    \varepsilon + 2 c\, \eta\, t^{-(d+2)/2}
    \;=\;
    2 \varepsilon,
  \end{align*}
  which concludes the proof of \eqref{eq:alclt:gradient}.
  %
\end{proof}
Using the local CLT we can deduce near-diagonal lower bounds for the annealed heat kernel and its first derivative.  Let us stress that the usual arguments to obtain on-diagonal lower bound (see for instance \cite[Proposition 4.3.4]{Kum14}) cannot be used in the time-dependent case since the proof uses crucially the symmetry of the heat kernel.  In particular, this shows that our upper bound is correct, at least near the diagonal. Notice that a local CLT and corresponding lower bound for the annealed second derivative is still open.
\begin{cor}
  \label{cor:hk:lbs}
  Suppose that Assumption~\ref{ass:law}, \ref{ass:integrability} and \ref{ass:time:derivative} hold true.  Then, there exists $n_{1} \in \mathbb{N}$ and $\const[C]{annealed:heatkernel:lb} > 0$ such that for all $t \geq n_{1}$ and $y, y' \in \mathbb{Z}^{d}$ with $|y'-y| \leq \sqrt{t}$
  \begin{align}\label{eq:hk:lb}
    \bar{p}_{0,t}(y, y') 
    \;\geq\;
    \const[C]{annealed:heatkernel:lb}\, t^{-d/2}.
  \end{align} 
\end{cor}
%
\begin{proof}
  Again, by \eqref{eq:shift:hk} and Assumption~\ref{ass:law}-(i) it suffices to prove the above assertions for $y = 0$.  The proof of Corollary~\ref{cor:hk:lbs}-(i) follows from Theorem~\ref{thm:annealed:LCLT} by exactly the arguments as given in \cite[Lemma~5.3]{ADS16}.  
\end{proof}  
\begin{cor}\label{cor:lclt:L2}
  Suppose that Assumption~\ref{ass:law} and \ref{ass:integrability} hold.  Additionally, assume that 
  \begin{align}\label{eq:ass:time-independent}
    \prob\bigl[\omega_{t}(e) = \omega_{0}(e)\; \forall\, t \in \mathbb{R} \bigr] 
    \;=\; 
    1,
    \qquad \forall\, e \in E_{d},
  \end{align}
  that is, the conductances are time-independent.  Then, for any $t > 0$ and $y \in \mathbb{R}^{d}$,
  \begin{align}\label{eq:lclt:L2}
    \lim_{n \to \infty}
    \mean\Bigl[
      \bigl| n^{d} p^{\omega}_{0, tn^{2}}(0, [y n]) - k_{t}^{\Sigma}(y)\bigr|^{2}
    \Bigr]
    \;=\;
    0.
  \end{align}
\end{cor}
\begin{proof}
  The proof comprises three steps.

  \textit{Step 1.}  We will first show that, for any $t > 0$,
  \begin{align}\label{eq:lclt:var:sum}
    \lim_{n \to \infty} 
      \frac{1}{n^{d}} 
      \sum_{z \in \mathbb{Z}^{d}}
      \var\Bigl[
        n^{d} p_{0, t n^{2}}^{\omega}(0, z)
      \Bigr]
    \;=\;
    0.
  \end{align}
  So, fix some $t > 0$.  Then, by applying the Chapman-Kolmogorov equation, it follows that, for any $n \in \mathbb{N}$,
  \begin{align*}
    \bar{p}_{0, 2t n^{2}}(0, 0)
    \;=\;
    \sum_{z \in \mathbb{Z}^{d}}
    \biggl(
      \bar{p}_{0, t n^{2}}(0, z)\, \bar{p}_{t n^{2}, 2t n^{2}}(z,0)
      \,+\,
      \cov\Bigl[ p_{0, tn^{2}}^{\omega}(0, z); p_{t n^{2}, 2 t n^{2}}^{\omega}(z,0) \Bigr]
    \biggr).
  \end{align*}
  By exploiting the stationarity of the law $\prob$ with respect to time-space shift implied by Assumption~\ref{ass:law}-(i) and using \eqref{eq:shift:hk}, we get that $\bar{p}_{t n^{2}, 2t n^{2}}(z, 0) = \bar{p}_{0, t n^{2}}(0, -z)$.  Further, by additionally using Lemma~\ref{lem:time_inverse}, the above covariance can be rewritten as
  \begin{align*}
    &\cov\Bigl[
      p_{0, t n^{2}}^{\omega}(0, z); p_{t n^{2}, 2t n^{2}}^{\omega}(z,0)
    \Bigr]
    \\[.5ex]
    &\mspace{32mu}\overset{\mspace{-9mu}\eqref{eq:shift:hk}\mspace{-9mu}}{\;=\;}
    \cov\Bigl[
      p_{-t n^{2}, 0}^{\omega}(-z, 0); p_{0, t n^{2}}^{\omega}(0, -z)
    \Bigr]
    \overset{\mspace{-9mu}\eqref{eq:hk:time_inverse}\mspace{-9mu}}{\;=\;}
    \cov\Bigl[
      p_{0, t n^{2}}^{\tilde{\omega}}(0, -z); p_{0, t n^{2}}^{\omega}(0, -z)
    \Bigr].
  \end{align*}
  Notice that, under the additional Assumption~\eqref{eq:ass:time-independent}, we have that
  \begin{align}\label{eq:time-shift:pt}
    \cov\Bigl[
      p_{0, t n^{2}}^{\tilde{\omega}}(0, -z); p_{0, t n^{2}}^{\omega}(0, -z)
    \Bigr]
    \;=\;
    \var\Bigl[ p_{0, t n^{2}}^{\omega}(0, z)\Bigr]
    \qquad \forall\, z \in \mathbb{Z}^{d}.
  \end{align}
  Thus, for any $n \in \mathbb{N}$ and $k \in \mathbb{N}$, the left-hand side of \eqref{eq:lclt:var:sum} can be bounded from above by
  \begin{align*}
    &
      \frac{1}{n^{d}} \sum_{z \in \mathbb{Z}^{d}}
      \var\Bigl[n^{d} p_{0, t n^{2}}^{\omega}(0, z)\Bigr]
    \;\leq\;
    I_{1}(n) + I_{2}(n,k) + I_{3}(n, k),
  \end{align*}
  where
  \begin{align*}
    I_{1}(n)
    &\;\ldef\;
    \bigl| 
      n^{d} \bar{p}_{0, 2t n^{2}}(0, 0) - k_{2t}^{\Sigma}(0)
    \bigr|
    \\[1ex]
    I_{2}(n, k)
    &\;\ldef\;
    \sum_{z \in \mathbb{Z}^{d} \setminus B(0, k n)}
    n^{d} \bar{p}_{0, t n^{2}}(0, z)\, \bar{p}_{0, t n^{2}}(0, -z)
    \\
    I_{3}(n, k)
    &\;\ldef\;
    \frac{1}{n^{d}} \sum_{z \in B(0, kn)} 
    \bigl|
      n^{d} \bar{p}_{0, t n^{2}}(0, z)\, n^{d} \bar{p}_{0, t n^{2}}(0, -z)
      -
      k_{2t}^{\Sigma}(0)
    \bigr|
  \end{align*}
  Thus, \eqref{eq:lclt:var:sum} follows once we have shown that the terms $I_{1}$, $I_{2}$ and $I_{3}$ vanishes when taking first the limit $n \to \infty$ and then $k \to \infty$.
  
  Clearly, Theorem~\ref{thm:annealed:LCLT}-(i) implies that $\lim_{n \to \infty} I_{1}(n) = 0$.  Moreover, by using Lemma~\ref{lem:Nashclas}, the Markov inequality and Proposition~\ref{prop:mean:displacement},
  \begin{align*}
    I_{2}(n, k)
    \;\leq\;
    c\, t^{-d/2}\, 
    \mean\bigl[
      \Prob_{0, 0}^{\omega}\bigl[ |X_{t n^{2}}| \geq k n\bigr]
    \bigr]
    \;\leq\;
    c\, \frac{t^{-(d-1)/2}}{k},
  \end{align*}
  which implies that $\limsup_{k \to \infty} \limsup_{n \to \infty} I_{2}(n, k) = 0$.  Thus, it remains to consider the term $I_{3}$.  But, by an elementary computation that involves multiple applications of Theorem~\ref{thm:annealed:LCLT}-(i), we obtain that
  \begin{align*}
    \limsup_{n \to \infty} I_{3}(n, k)
    &\;\leq\;
    \Biggl|
      \int_{\mathbb{R}^{d}} 
        \indicator_{\{|x| \leq k\}}\, k_{t}^{\Sigma}(x)\, k_{t}^{\Sigma}(0 - x)\, 
      \md x
      \,-\,
      k_{2t}^{\Sigma}(0)
    \Biggr|.
  \end{align*}
  Therefore, by using again the Chapman-Kolmogorov equation as well as the monotone convergence theorem which yields
  \begin{align*}
    \lim_{k \to \infty} 
    \int_{\mathbb{R}^{d}} 
      \indicator_{\{|x| > k\}}\, k_{t}^{\Sigma}(x)\, k_{t}^{\Sigma}(0 - x)\, 
    \md x 
    \;=\; 
    0,
  \end{align*}
  we conclude that $\limsup_{k \to \infty} \limsup_{n \to \infty} I_{3}(n, k) = 0$.
  
  \textit{Step 2.}  Fix some $t > 0$ and $y \in \mathbb{R}^{d}$.  For any $\delta > 0$, $n \in \mathbb{N}$ and $z \in B([y n], \delta n)$ let $\gamma_z$ be a shortest path in the graph metric between $[y n]$ and $z$, that is, $\gamma_z = (y_{0}, \ldots, y_k)$ with $y_{0} = [y n]$, $y_k = z$ and $\{y_{i-1}, y_{i}\} \in E_{d}$ for all $i \in \{1, \ldots, k\}$ and $k \leq c_{1} \delta n$ for some $c_{1} \in [1, \infty)$ that is independent of $\delta$ and $n$.  Then, by applying the Minkowski inequality, we get
  \begin{align*}
    &\var\Bigl[ n^{d}\, p_{0, t n^{2}}^{\omega}(0, [y n]) \Bigr]^{1/2}
    \\
    &\mspace{36mu}\leq\;
    \var\Bigl[ n^{d}\, p_{0, t n^{2}}^{\omega}(0, z) \Bigr]^{1/2}
    +\,
    \sum_{i=1}^k 
    \max_{(0,x) \in E_{d}}
    \var\Bigl[ 
      n^{d}\, \nabla_x^{2}\, p_{0, t n^{2}}^{\omega}(y_{i-1}, y_{i}) 
    \Bigr]^{1/2}.
  \end{align*}
  By using Theorem~\ref{thm:gradient:est}-(i), we get that, for any $(0, x) \in E_{d}$,
  \begin{align*}
    \var\Bigl[ 
      n^{d}\, \nabla_x^{2} p_{0, t n^{2}}^{\omega}(y_{i-1}, y_{i}) 
    \Bigr]^{1/2}
    \;\leq\;
    2 \mean\Bigl[
      \bigl|
        n^{d}\, \nabla_x^{2}\, p_{0, t n^{2}}^{\omega}(y_{i-1}, y_{i}) 
      \bigr|^{2}
    \Bigr]^{1/2}
    \overset{\eqref{eq:gradient:on-diag}}{\;\leq\;}
    \frac{2 \const[C]{gradient:point}}{n\, t^{(d+1)/2}}.
  \end{align*}
  By combining the estimates above, we obtain that, for any $\delta > 0$,
  \begin{align*}
    \var\Bigl[ n^{d}\, p_{0, t n^{2}}^{\omega}(0, [y n]) \Bigr]
    \;\leq\;
    \frac{2}{|B([yn], \delta n)|}\,
    \sum_{z \in B([y n], \delta n)}
    \var\Bigl[ n^{d}\, p_{0, t n^{2}}^{\omega}(0, z) \Bigr]
    \,+\,
    \frac{2c \delta^{2}}{t^{d+1}}.
  \end{align*}
  Thus, by applying \eqref{eq:lclt:var:sum}, we conclude the right-hand side of the above estimate vanishes when we first let $n \to \infty$ and then $\delta \downarrow 0$.  This proves that
  \begin{align}\label{eq:lclt:variance}
    \lim_{n \to \infty} \var\Bigl[n^{d} p^{\omega}_{0, tn^{2}}(0, [y n])\Bigr]
    \;=\;
    0.
  \end{align}

  \textit{Step 3.}  Finally, for any $t > 0$ and $y \in \mathbb{R}$, we rewrite the left-hand side of \eqref{eq:lclt:L2} as
  \begin{align}\label{eq:lclt:L2:split}
    &\mean\Bigl[
      \bigl| n^{d} p^{\omega}_{0, tn^{2}}(0, [y n]) - k_{t}^{\Sigma}(y)\bigr|^{2}
    \Bigr]
    \nonumber\\[.5ex]
    &\mspace{36mu}=\;
    \var\Bigl[n^{d} p^{\omega}_{0, tn^{2}}(0, [y n])\Bigr]
    \,+\,
    \bigl|n^{d}\, \bar{p}_{0, t n^{2}}(0, [yn]) - k_{t}^{\Sigma}(y)\bigr|^{2}.
  \end{align}
  Now, by applying Theorem~\ref{thm:annealed:LCLT}-(i) and \eqref{eq:lclt:variance}, we immediately obtain that the right-hand side of \eqref{eq:lclt:L2:split} vanishes as $n$ tends to infinity. This concludes the proof.
\end{proof}
\begin{remark}  
  In order to deduce \eqref{eq:time-shift:pt}, it is enough to assume
  \begin{align}\label{eq:ass:time-symmetric}
    \prob\bigl[
      \omega_{t}(e) = \omega_{-t}(e)\; \forall\, t \in \mathbb{R}
    \bigr]
    \;=\;
    1,
    \qquad \forall\, e \in E_{d}.
  \end{align}
  However, Assumption~\ref{ass:law}-(i) and \eqref{eq:ass:time-symmetric} imply 
  \eqref{eq:ass:time-independent}, so Corollary~\ref{cor:lclt:L2} applies to the case of time-independent conductances only. 
\end{remark}
\begin{cor}
  Suppose that Assumption~\ref{ass:law} and \ref{ass:integrability} hold true.
  \begin{enumerate}[(i)]
    \item If $d \geq 3$ and Assumption~\ref{ass:offdiag} is satisfied for some $\alpha > d$ and $\beta \geq \alpha/2$, then for any compact subset $K \subset \mathbb{R}^{d} \setminus \{0\}$
    \begin{align}\label{eq:alclt:green}
      \lim_{n \to \infty} \sup_{y \in K}\,
      \biggl|
        n^{d-2}\, \int_{0}^{\infty} \bar{p}_{0, t}(0, [y n])\, \md t 
        - 
        \int_{0}^{\infty} k_{t}^{\Sigma}(y)\, \md t
      \biggr|
      \;=\;
      0.
    \end{align}
    In particular, there exists $\const[C]{annealed:greenkernel:lb} > 0$ such that, for any $y, y' \in \mathbb{Z}^{d}$ with $y \ne y'$,
    \begin{align}\label{eq:green:lb}
      \int_{0}^{\infty} \bar{p}_{0, t}(y, y')\, \md t
      \;\geq\;
      \const[C]{annealed:greenkernel:lb}\, |y-y'|^{-d+2}.
    \end{align}
    \item If $d \geq 2$ and Assumption~\ref{ass:offdiag} is satisfied for some $\alpha \geq 2(d-1)$ and $\beta > d + \alpha/2$, then for any compact subset $K \subset \mathbb{R}^{d} \setminus \{0\}$ and $i \in \{1, \ldots, d\}$
    \begin{align}\label{eq:alclt:green:gradient}
      \lim_{n \to \infty} \sup_{y \in K}\,
      \biggl|
        n^{d-1}\, 
        \int_{0}^{\infty} 
          \nabla_{\!e_{i}}^1\, \bar{p}_{0, t}(0, [y n])\, 
        \md t 
        - 
        \int_{0}^{\infty} \bigl(\partial_{i} k_{t}^{\Sigma}\bigr)(y)\, \md t
      \biggr|
      \;=\;
      0.
    \end{align}
    Moreover, there exists $\const[C]{annealed:gradient:green:lb} \in (0, \infty)$ such that, for any $(0, x) \in E_{d}$ and $y, y' \in \mathbb{Z}^{d}$ with $y \ne y'$,
    \begin{align}
      \biggl|
        \int_{0}^{\infty} \nabla_{\!x}^1\, \bar{p}_{0, t}(y, y')\, \md t
        \biggr|
      \;\geq\;
      \const[C]{annealed:gradient:green:lb}\, |y-y'|^{-d+1}.
    \end{align}
  \end{enumerate}
\end{cor}
\begin{proof}
  (i) Fix some compact subset $K \subset \mathbb{R}^{d} \setminus \{0\}$.  To lighten notation, we set 
  \begin{align*}
    \bar{g}(y, y')
    \;\ldef\;
    \int_{0}^{\infty} \bar{p}_{0, t}(y, y')\, \md t
    \qquad \text{and} \qquad
    g^{\Sigma}(y)
    \;\ldef\;
    \int_{0}^{\infty} k_{t}^{\Sigma}(y)\, \md t.
  \end{align*}
  Fix some $y \in  K$.  Then, for any $\delta > 0$ and $n \in \mathbb{N}$, the left-hand side of \eqref{eq:alclt:green} is bounded from above by
  \begin{align}\label{eq:alclt:green:est}
    \sup_{y \in K}\,
    \bigl|
      n^{d-2}\, \bar{g}(0, [y n]) 
      - g^{\Sigma}(y)
    \bigr|
    \;\leq\;
    I_{1}(\delta, n) + I_{2}(\delta, n) + I_{3}(\delta),
  \end{align}
  where
  \begin{align*}
    I_{1}(\delta, n)
    &\;\ldef\;
    \int_{\delta}^{\infty}
      \sup_{y \in K}\,
      \bigl|
        n^{d} \bar{p}_{0, tn^{2}}(0, [y n]) - k_{t}^{\Sigma}(y)
      \bigr|\,
    \md t,
    \\[.5ex]
    I_{2}(\delta, n)
    &\;\ldef\;
    \sup_{y \in K}\, 
    n^{d-2} \int_{0}^{\delta n^{2}} \bar{p}_{0, t}(0, [y n])\, \md t,
    \\[.5ex]
    I_{3}(\delta)
    &\;\ldef\;
    \sup_{y \in K}\,
    \int_{0}^{\delta} k_{t}^{\Sigma}(y)\, \md t.
  \end{align*}
  Thus, it suffices to show that the right-hand side of \eqref{eq:alclt:green:est} vanishes when we first take the limit $n \to \infty$ and then $\delta \downarrow 0$.  
  
  First, by Theorem~\ref{thm:annealed:LCLT}, the function $f_{n}(t) \ldef \sup_{y \in K} |n^{d} \bar{p}_{0, t n^{2}}(0, [y n]) - k_{t}(y)|$ converges pointwise to zero as $n$ tends to infinity.  (Note that for this fact, we fix time, so Assumption~\ref{ass:time:derivative} is not needed.)  Moreover, in view of Lemma~\ref{lem:Nashclas} and the properties of $k^{\Sigma}$, it holds that there exists $c_{1} < \infty$ such that $f_{n}(t) \leq c_{1} t^{-d/2}$ for any $n \in \mathbb{N}$.  Since for any $d \geq 3$ and $\delta > 0$ the function $t \mapsto t^{-d/2}$ is integrable over $[\delta, \infty)$, an application of Lebesgue's dominated convergence theorem yields that $\lim_{n \to \infty} I_{1}(\delta, n) = 0$. 
  
  Next, consider the term $I_{2}$.  Since $K$ is a compact subset of $\mathbb{R}^{d} \setminus \{0\}$ and the map $y \mapsto |y|$ is continuous, there exist $0 < k_{*} \leq k_{**} < \infty$ such that $\min_{y \in K} |y| \geq k_{*}$ and $\max_{y \in K} |y| \leq k_{**}$.  Then, for any fix $\delta > 0$ there exists $n(\delta, k_{*}, k_{**}) \in \mathbb{N}$ such that for any $n \geq n(\delta, k_{*}, k_{**})$ we have that $|[y n] / n| \geq k_{*} / 2$ and $|[y n]| \leq \delta n^{2}$ for any $y \in K$.  In particular, by applying Lemma~\ref{lem:hk:offdiag:p} with $p=1$, $\alpha > d$ and $\beta \geq \alpha/2$ we obtain that, for any $y \in K$,
  \begin{align*}
    &n^{d-2} \int_{0}^{\delta n^{2}} \bar{p}_{0, t}(0, [y n])\, \md t
    \\[.5ex]
    &\mspace{32mu}\overset{\mspace{-9mu}\eqref{eq:hk:offdiag:p}\mspace{-9mu}}{\;\leq\;}
    \const[C]{mean:p:heatkernel}\, n^{d-2} 
    \biggl(
      \int_{0}^{|[y n]|} k_{\alpha}(t, [y n])\, \md t
      +
      \int_{|[y n]|}^{\delta n^{2}} k_{\alpha}(t, [y n])\, \md t
    \biggr)
    \\[.5ex]
    &\mspace{36mu}\leq\;
    c\,
    \biggl(
      n^{d/2-\alpha/2-1}\, k_{*}^{-d/2-\alpha/2+1}
      \,+\,
      \delta^{-d/2+\alpha/2+1}\, k_{*}^{-\alpha}
    \biggr).
  \end{align*}
  Hence,
  \begin{align*}
    \limsup_{\delta \searrow 0} \limsup_{n \to \infty} I_{2}(\delta, n)
    \;=\;
    0.
  \end{align*}
  Finally, let us address the term $I_{3}$.  Since $\sup_{t \geq 0} \sup_{y \in K} k_{t}^{\Sigma}(y) \leq c(k_{*}) < \infty$, we get that $I_{3}(\delta) \leq \delta c(k_{*})$, and $\limsup_{\delta \downarrow 0} I_{3}(\delta)$ follows immediately.  This completes the proof of \eqref{eq:alclt:green}.

  It remains to show the claimed lower bound for the annealed Green's kernel.  In view of \eqref{eq:shift:hk} and Assumption~\ref{ass:law}-(i) it suffices to prove \eqref{eq:green:lb} for $y = 0$.  First, note that there exists $c_{1} \in (0, \infty)$ such that $g^{\Sigma}(x) \geq 2c_{1} |x|^{2-d}$ for all $x \in \mathbb{R}^{d} \setminus \{0\}$.  Moreover, \eqref{eq:alclt:green} implies that there exists $n_{0} \in \mathbb{N}$ such that, for any $n \geq n_{0}$,
  \begin{align*}
    \sup_{1/2 \leq |x| \leq 2} 
    \bigl| n^{d-2} \bar{g}(0, [x n]) - g^{\Sigma}(x)\bigr| 
    \;\leq\; 
    c_{1}
  \end{align*}
  Hence, for any $y' \in \mathbb{Z}^{d}$ such that $n \ldef |y'| \geq n_{0}$, we obtain that
  \begin{align*}
    \bar{g}(0, y')
    \;\geq\;
    n^{2-d}\,
    \Bigl(
      n^{d-2} g^{\Sigma}(y' / n)
      -
      \sup_{1/2 \leq |x| \leq 2} 
      \bigl|
        n^{d-2} \bar{g}(0, [x n]) - g^{\Sigma}(x)
      \bigr|
    \Bigr)
    \;\geq\;
    c_{1} |y'|^{2-d}.
  \end{align*}
  Thus, by setting $\const[C]{annealed:greenkernel:lb} \ldef \min\{|y'|^{d-2} \bar{g}(0, y') : y' \in \mathbb{Z}^{d} \text{ s.th. } |y'| < n_{0}\} \wedge c_{1}$, the assertion \eqref{eq:green:lb} follows.

  (ii) The proof of (ii) relies on the same arguments as in (i) except for the fact that in the proof of \eqref{eq:alclt:green:gradient} we apply Theorem~\ref{thm:gradient:est}-(i) instead of Lemma~\ref{lem:Nashclas} and Proposition~\ref{prop:gradient:offdiag} instead of Lemma~\ref{lem:hk:offdiag:p}.  Moreover, we use that there exists $c(k_{*}, k_{**}) < \infty$ such that $\sup_{t \geq 0} \sup_{y \in K} \bigl(\partial_{i} k_{t}^{\Sigma}\bigr)(y) \leq c(k_{*}, k_{**})$ for any $i \in \{1, \ldots, d\}$.
\end{proof}

\appendix
\section{Forward and backward equation for the semigroup}
\label{appendix:semigroup}
Let us briefly recall the construction of the time-inhomogeneous Markov process $X$ starting at time $s \in \mathbb{R}$ in $x \in \mathbb{Z}^{d}$, cf.~\cite[Section~4]{ACDS18}.  Let $(E_{n} : n \in \mathbb{N})$ be a sequence of independent $\mathop{\mathrm{Exp}}(1)$-distributed random variables.  Further, set $\pi_{t}^{\omega}(x, y) \ldef \omega_{t}(x, y) / \mu_{t}^{\omega}(x) \indicator_{\{(x, y) \in E_{d}\}}$, where $\mu_{t}^{\omega}(x) \ldef \sum_{y : (x, y) \in E_{d}} \omega_{t}(x, y)$ for any $t \in \mathbb{R}$, $x \in \mathbb{Z}^{d}$.  We specify both the sequence of jump times, $(J_{n} : n \in \mathbb{N}_{0})$ and positions, $(Y_{n} : n \in \mathbb{N}_{0})$, inductively.  For this purpose, set $J_{0} = s$ and $Y_{0} = x$.  Suppose that, for any $n \geq 1$, we have already constructed the random variables $(J_{0}, Y_{0}, \ldots, J_{n-1}, Y_{n-1})$.  Then, $J_{n}$ is given by
\begin{align*}
  J_{n}
  \;=\;
  J_{n-1} \,+\,
  \inf\biggl\{
    t \geq 0
    \;:\;
    \int_{J_{n-1}}^{J_{n-1} + t} \mu_{u}^{\omega}(Y_{n-1})\, \md u \geq E_{n}
  \biggr\},
\end{align*}
and at the jump time $J_{n}$ the distribution of $Y_{n}$ is given by $\pi_{J_{n}}^{\omega}(Y_{n-1}, \cdot)$.  Since, under Assumption~\ref{ass:law}, $\sup_{n \in \mathbb{N}_{0}} J_{n} = \infty$ the Markov process $X$ is given by
\begin{align*}
  X_{t} \;=\; Y_{n} \quad \text{on} \quad [J_{n}, J_{n+1})
  \qquad \forall\, n \in \mathbb{N}_{0}.
\end{align*}
Note that, under $\Prob_{s, x}^{\omega}$, $J_{0} = s$ and $Y_{0} = x$ almost surely, the conditional law of $J_{n}$ given $(J_{0}, Y_{0}, \ldots, J_{n-1}, Y_{n-1})$ (also called survival distribution with time-dependent hazard rate $\mu_{t}(Y_{n-1})$) is
\begin{align*}
  \mu_{t}^{\omega}(Y_{n-1})\, \me^{-\int_{J_{n-1}}^{t} \mu_{u}^{\omega}(Y_{n-1})\, \md u}\,
  \indicator_{\{t > J_{n-1}\}}\, \md t,
\end{align*}
and the conditional law of $Y_{n}$ given $(J_{0}, Y_{0}, \ldots, J_{n-1}, Y_{n-1}, Y_{n})$ is $\pi_{J_{n}}^{\omega}(Y_{n-1}, \cdot)$.

Note that for the Markov process $X$ as constructed above the strong Markov property holds true.  Thus, an application of the strong Markov property yields that $(P_{s,t}^{\omega} : t \geq s)$ satisfies the integrated backward equation, that is, for $\prob$-a.e.\ $\omega$,
\begin{align}
  \label{eq:backward:integrated}
  (P_{s,t}^{\omega} f)(x)
  \;=\;
  \me^{-\int_{s}^{t} \mu_{u}^{\omega}(x)\, \md u} f(x)
  \,+\,
  \int_{s}^{t}
    \me^{-\int_{s}^{r} \mu_{u}^{\omega}(x)\, \md u}\!
    \sum_{y:(x,y) \in E_{d}}\mspace{-9mu} \omega_{r}(x,y)\, (P_{r,t}^{\omega} f)(y)\,
  \md r
\end{align}
for any $f \in \ell^{\infty}(\mathbb{Z}^{d})$, $-\infty < s < t < \infty$ and $x \in \mathbb{Z}^{d}$.  
\begin{prop}\label{prop:backward}
  For $\prob$-a.e.\ $\omega$, every $x,y \in \mathbb{Z}^{d}$ and $f \in \ell^{\infty}(\mathbb{Z}^{d})$ the following hold.
  \begin{enumerate}[(i)]
  \item
    For every $t \in \mathbb{R}$, the map $(-\infty, t) \ni s \mapsto p_{s,t}^{\omega}(x,y)$ is absolutely continuous and hence differentiable for almost every $s \in (-\infty, t)$.  In particular, $\lim_{s \uparrow t} p_{s,t}^{\omega}(x,y) = p_{t,t}^{\omega}(x,y) = \indicator_{\{y\}}(x)$.

  \item
    For every $s \in \mathbb{R}$, the map $(s,\infty) \ni t \mapsto p_{s, t}^{\omega}(x, y)$ is continuous.  Hence, it holds that $\lim_{t \downarrow s} p_{s, t}^{\omega}(x, y) = p_{s, s}^{\omega}(x,y) = \indicator_{\{y\}}(x)$.

  \item
    \textnormal{(Backward equation)} 
    It holds that for every $t \in \mathbb{R}$,
    \begin{align} \label{eq:diff:backward}
      - \partial_{s}\, (P_{s,t}^{\omega} f)(x)
      \;=\;
      \bigl(\mathcal{L}_{s}^{\omega} (P_{s,t}^{\omega} f)\bigr)(x),
      \qquad \text{for a.e. } s \in (-\infty, t).
    \end{align}
    In particular, the function $(s, x) \mapsto u(s, x) \ldef p^{\omega}_{s, 0}(x, 0)$ solves
    \begin{align}
      \partial_{s} u(s, x) \;=\; -(\mathcal{L}_{s}^{\omega} u(s, \cdot))(x),
      \qquad \forall\, x \in \mathbb{Z}^{d} \text{ and a.e. } s \in (-\infty, 0).
    \end{align}
    Moreover, for any function $g$ with finite support and $t \in \mathbb{R}$,
    \begin{align}\label{eq:beqweak}
      \partial_{s} \scpr{g}{P_{s,t}^{\omega}f}{\mathbb{Z}^{d}}
      \;=\;
      \mathcal{E}_{s}^{\omega}\bigl(g, P_{s,t}^{\omega}f\bigr),
      \qquad \text{for a.e. } s \in (-\infty, t).
    \end{align}
  \end{enumerate}
\end{prop}
\begin{proof}
  (i) By \eqref{eq:local:integrability} the map $t \mapsto \mu_{t}^{\omega}(x)$ is $\prob$-a.s.\ locally integrable for every $x \in \mathbb{Z}^{d}$.  Hence, the absolute continuity of the Lebesgue integral implies that, for every $\varepsilon > 0$, there exists $\delta \equiv \delta(x) > 0$ such that
  \begin{align*}
    \int_{D} \mu^{\omega}_{u}(x)\, \md u
    \;<\;
    \varepsilon
    \qquad \forall\, D \in \mathcal{B}(\mathbb{R})
    \text{ with Lebesgue measure less than } \delta.
  \end{align*}
  Thus, it remains to observe that, for each $f \in \ell^{\infty}(\mathbb{Z}^{d})$ and $t \in \mathbb{R}$, using the integrated backward equation as well as the Cauchy-Schwarz inequality,
  \begin{align*}
    &\sum_{i=1}^{n}
    \bigl|(P_{s_{i},t}^{\omega} f)(x) - (P_{r_{i},t}^{\omega}f)(x)\bigr|
    \\[0ex]
    &\mspace{36mu}\leq\;
    2 \|f\|_{\infty}\,
    \sum_{i=1}^{n} \Bigl(1 - \me^{-\int_{r_{i}}^{s_{i}} \mu_{u}^{\omega}(x)\, \md u} \Bigr)
    \;\leq\;
    2 \|f\|_{\infty}\, \biggl(\int_{D} \mu_{u}^{\omega}(x)\, \md u \biggr)
    \;<\;
    2\varepsilon \|f\|_{\infty}
  \end{align*}
  for any union $D = \bigcup_{i=1}^{n} (r_{i}, s_{i})$ of pairwise disjoint intervals $(r_{i}, s_{i}) \subset (-\infty, t]$ of total length less than $\delta$.
  
  (ii) This can be deduced from the absolute continuity of the Lebesgue integral.  Indeed, $\prob$-a.s., for any $f \in \ell^{\infty}(\mathbb{Z}^{d})$ and $r,t \in (s, \infty)$ with $r < t$,
  \begin{align*}
    \bigl|(P_{s,r}^{\omega} f)(x) - (P_{s,t}^{\omega}f)(x)\bigr|
    &\;=\;
    \bigl| \bigl(P_{s,r}^{\omega} (f - P_{r,t}^{\omega} f) \bigr)(x) \bigr|
    \\[.5ex]
    &\;\leq\;
    2 \|f\|_{\infty}\,
    \sum_{y \in \mathbb{Z}^{d}} p_{s,r}^{\omega}(x,y)
    \Bigl(1 - \me^{-\int_{r}^{t} \mu_{u}^{\omega}(y)\, \md u}\Bigr).
  \end{align*}
  Thus, by applying Lebesgue's dominated convergence theorem, it follows that
  \begin{align*}
    \lim_{t \searrow r} |(P_{s,r}^{\omega} f)(x) - (P_{s,t}^{\omega} f)(x)| \;=\; 0.
  \end{align*}
  
  (iii) Note that the right-hand side of \eqref{eq:backward:integrated} can be rewritten as
  \begin{align*}
    \me^{-\int_{s}^{t} \mu_{u}^{\omega}(x)\, \md u}\, f(x)
    \,+\,
    \me^{\int_{0}^{s} \mu_{u}^{\omega}(x)\, \md u}\,
    \int_{s}^{t}
    \me^{-\int_{0}^{r} \mu_{u}^{\omega}(x)\, \md u}
    \sum_{y:(x,y) \in E_{d}}\mspace{-9mu} \omega_{r}(x,y) (P_{r,t}^{\omega} f)(y)\,
    \md r.
  \end{align*}
  Since $t \mapsto \mu_{t}^{\omega}(x)$ is locally integrable, the differential form of the backward equation in weak sense follows from \cite[Theorem~6.3.6]{Co13} together with an application of the chain and product rule.
\end{proof}
\begin{lemma} \label{lem:time_inverse}
  Define $\tilde{\omega}_{t}(e) \ldef \omega_{-t}(e)$ for any $t \in \mathbb{R}$ and $e \in E_{d}$.  Then,
  \begin{align}\label{eq:hk:time_inverse} 
    p_{s,t}^{\omega}(x, y)
    \;=\;
    p_{-t,-s}^{\tilde{\omega}}(y, x),
    \qquad \forall\, x, y \in \mathbb{Z}^{d}, \quad s \leq t.
  \end{align}
\end{lemma}
\begin{proof}
  Write $B_{n} \ldef B(0, n)$ and $T_{n} \ldef \inf\{t \geq s: |X_{t}| > n \}$ with $\inf \emptyset \ldef \infty$.  We denote by $p^{\omega, B_{n}}_{s, t}(x, y) \ldef \Prob_{s, x}^{\omega}\!\big[X_{t} = y, \, t < T_{n} \big]$ the heat kernel associated with the process $X$ killed upon exiting $B_{n}$, and we write $(P_{s, t}^{\omega, B_{n}}: t \geq s)$ for the corresponding transition semigroup.  Recall that the associated time-dependent generator, still denoted by $\mathcal{L}^{\omega}_{t}$, is acting on functions with Dirichlet boundary condition.  By similar arguments as in Proposition~\ref{prop:backward} one can establish a backward equation for $(P_{s, t}^{\omega, B_{n}} : t \geq s)$, which gives, for any $s \leq t$,
  \begin{align*}
    &\partial_{u} \scpr{P_{-u, -s}^{\tilde{\omega}, B_{n}} g}{P_{u, t}^{\omega, B_{n}} f}{B_{n}}
    \\[1ex]
    &\mspace{36mu}=\;
    \scpr{\mathcal{L}_{-u}^{\tilde{\omega}} P_{-u, -s}^{\tilde{\omega}, B_{n}} g}{P_{u, t}^{\omega, B_{n}} f}{B_{n}} 
    -
    \scpr{ P_{-u, -s}^{\tilde{\omega}, B_{n}} g }{\mathcal{L}_{u}^{\omega} P_{u, t}^{\omega, B_{n}} f}{B_{n}}
    \;=\;
    0,
  \end{align*}
  where we used in the last step that $\mathcal{L}_{-u}^{\tilde{\omega}} = \mathcal{L}_{u}^{\omega}$.  By integration over $[s, t]$ we get
  \begin{align*}
    \scpr{P_{-t, -s}^{\tilde{\omega}, B_{n}} g}{f}{B_{n}}
    - \scpr{g}{P_{s, t}^{\omega, B_{n}} f}{B_{n}}
    \;=\;
    0,
  \end{align*}
  and by choosing $f = \indicator_{\{y\}}$ and $g = \indicator_{\{x\}}$ we obtain $p_{-t, -s}^{\tilde{\omega}, B_{n}}(y, x) = p_{s, t}^{\omega, B_{n}}(x, y)$.  Finally, since, by a similar reasoning as in the proof of Lemma~\ref{lem:Nashclas}, $\lim_{n \to \infty} p_{s, t}^{\omega, B_{n}}(x, y) = p_{s, t}^{\omega}(x, y)$ for all $x, y \in \mathbb{Z}^{d}$, $t \geq s$ and $\omega \in \Omega$, the result follows by taking the limit $n \to \infty$.
\end{proof}
\begin{prop}[Forward equation] \label{prop:forward}
  For $\prob$-a.e.\ $\omega$, every $x \in \mathbb{Z}^{d}$, $s \in \mathbb{R}$ and finitely supported $f\colon \mathbb{Z}^{d} \to \mathbb{R}$, the map $t \mapsto (P_{s, t}^{\omega}f)(x)$ is differentiable at almost every $t \in (s, \infty)$ and
  \begin{align}\label{eq:diff:forward}
    \partial_{t} (P_{s, t}^{\omega} f)(x)
    \;=\;
    \bigl(P_{s, t}^{\omega} (\mathcal{L}_{t}^{\omega} f)\bigr)(x),
    \qquad \text{for a.e. } t \in (s, \infty).
  \end{align}
  %
  In particular, the function $(t, x) \mapsto u(t, x) \ldef p^{\omega}_{0, t}(0, x)$ solves
  \begin{align}
    \partial_{t} u(t, x) \;=\; (\mathcal{L}_{t}^{\omega} u(t, \cdot))(x),
    \qquad \forall\, x \in \mathbb{Z}^{d} \text{ and a.e. } t \in (0,\infty).
  \end{align}
\end{prop}
\begin{proof}
  This follows from the backward equation  in Proposition~\ref{prop:backward} and Lemma~\ref{lem:time_inverse}.  Indeed, let $\tilde{\omega}$ be defined as in Lemma~\ref{lem:time_inverse}, then we have for any $f,g \in \ell^{2}(\mathbb{Z}^{d})$, 
  \begin{align*}
    \scpr{P_{s, t}^{\omega} f}{g}{\mathbb{Z}^{d}}
    \;=\;
    \scpr{f}{(P_{s, t}^{\omega})^{*}g}{\mathbb{Z}^{d}}
    \overset{\!\!\eqref{eq:hk:time_inverse}\!\!}{\;=\;}
    \scpr{f}{P_{-t, -s}^{\tilde{\omega}}g}{\mathbb{Z}^{d}}
  \end{align*}
  for any $s \leq t$.  Thus, for any finitely supported $f\colon \mathbb{Z}^{d} \to \mathbb{R}$ and $g = \indicator_{\{x\}}$, we get
  \begin{align*}
    \partial_{t} (P_{s, t}^{\omega} f)(x)
    &\;=\;
    \lim_{h \to 0} \frac{1}{h}
    \scpr{P_{s, t+h}^{\omega} f - P_{s, t}^{\omega} f}{g}{\mathbb{Z}^{d}}
    \\
    &\;=\;
    \lim_{h \to 0} \frac{1}{h}
    \scpr{f}{P_{-(t+h), -s}^{\tilde{\omega}} g - P_{-t, -s}^{\tilde{\omega}} g}{\mathbb{Z}^{d}}
    \\[1ex]
    &\;=\;
    - \scpr{f}{\partial_{r} (P_{r, -s}^{\tilde{\omega}} g)\big|_{r = -t}}{\mathbb{Z}^{d}}.
  \end{align*}
  Hence, by applying the differential backward equation to $P_{r, -s}^{\tilde{\omega}}$, we get
  \begin{align*}
    -\scpr{f}{\partial_{r} (P_{r, -s}^{\tilde{\omega}} g)\big|_{r=-t}}{\mathbb{Z}^{d}}
    &\;=\;
    \scpr{f}{\mathcal{L}_{-t}^{\tilde{\omega}}(P_{-t, -s}^{\tilde{\omega}} g)}{\mathbb{Z}^{d}}
    \\[.5ex]
    &\overset{\!\!\eqref{eq:hk:time_inverse}\!\!}{\;=\;}
    \scpr{f}{\mathcal{L}_{t}^{\omega} ((P_{s, t}^{\omega})^{*}g)}{\mathbb{Z}^{d}}
    \;=\;
    P_{s, t}^{\omega}(\mathcal{L}_{t}^{\omega} f)(x),
  \end{align*}
  which yields \eqref{eq:diff:forward}.  Finally, consider the function $u(t, x) \ldef p_{0, t}^{\omega}(0,x)$.  Then, by applying \eqref{eq:diff:forward}, we find that
  \begin{align*}
    \partial_{t} u(t, x)
    \;=\;
    \partial_{t} (P_{0, t}^{\omega} \indicator_{\{x\}})(0)
    &\;=\;
    P_{0, t}^{\omega}(\mathcal{L}_{t}^{\omega} \indicator_{\{x\}})(0)
    \\
    &\;=\;
    \scpr{u(t, \cdot)}{\mathcal{L}_{t}^{\omega} \indicator_{\{x\}}}{\mathbb{Z}^{d}}
    \;=\;
    \big(\mathcal{L}_{t}^{\omega} u(t, \cdot)\big)(x),
  \end{align*}
  which concludes the proof.
\end{proof}
\begin{cor} \label{prop:forward:backward:adjoint}
  For $\prob$-a.e.\ $\omega$, and every $x \in \mathbb{Z}^{d}$ the following holds.
  \begin{enumerate}[(i)]
  \item \textnormal{(Forward equation)} For every $s \in \mathbb{R}$ and $f \in \ell^{\infty}(\mathbb{Z}^{d})$
    \begin{align}\label{eq:diff:forward:adjoint}
      \partial_{t} \bigl((P_{s, t}^{\omega})^{*} f\bigr)(x)
      \;=\;
      \bigl(\mathcal{L}_{t}^{\omega} (P_{s, t}^{\omega})^{*} f\bigr)(x),
      \qquad \text{for a.e. } t \in (s, \infty).
    \end{align}
  \item \textnormal{(Backward equation)} For every $t \in \mathbb{R}$ and $f\colon \mathbb{Z}^{d} \to \mathbb{R}$ finitely supported
    \begin{align}\label{eq:diff:backward:adjoint}
      -\partial_{s} \bigl((P_{s, t}^{\omega})^{*}f \bigr)(x)
      \;=\;
      \bigl((P_{s, t}^{\omega})^{*}(\mathcal{L}_{s}^{\omega} f)\bigr)(x),
      \qquad \text{for a.e. } s \in (0, t).
    \end{align}
  \end{enumerate}
\end{cor}
\begin{proof}
  (i) For $\prob$-a.e.\ $\omega$, every $x \in \mathbb{Z}^{d}$, $s \in \mathbb{R}$ and $f \in \ell^{\infty}(\mathbb{Z}^{d})$ we have that
  \begin{align*}
    \partial_{t} \bigl((P_{s, t}^{\omega})^{*} f\bigr)(x)
    \overset{\!\!\eqref{eq:hk:time_inverse}\!\!}{\;=\;}
    \partial_{t} \bigl(P_{-t, -s}^{\tilde{\omega}} f\bigr)(x)
    \overset{\!\!\eqref{eq:diff:backward}\!\!}{\;=\;}
    \bigl(\mathcal{L}_{-t}^{\tilde{\omega}} (P_{-t, -s}^{\tilde{\omega}} f)\bigr)(x)
    \overset{\!\!\eqref{eq:hk:time_inverse}\!\!}{\;=\;}
    \bigl(\mathcal{L}_{t}^{\omega} (P_{s, t}^{\omega})^{*} f\bigr)(x)
  \end{align*}
  for a.e.\ $t \in (s, \infty)$.

  (ii) Recall that $(P_{s, t}^{\omega})^{*}$ denotes the adjoint of $P_{s, t}^{\omega}$ on $\ell^{2}(\mathbb{Z}^{d})$. Hence, $\prob$-a.s.\ for every $x \in \mathbb{Z}^{d}$, $t \in \mathbb{R}$ and finitely supported $f\colon \mathbb{Z}^{d} \to \mathbb{R}$,
  \begin{align*}
    \partial_{s} \bigl((P_{s, t}^{\omega})^{*}f\bigr)(x)
    &\;=\;
    \partial_{s} \scpr{P_{s, t}^{\omega} \indicator_{\{x\}}}{f}{\mathbb{Z}^{d}}
    \\[.5ex]
    &\overset{\!\!\eqref{eq:diff:forward}\!\!}{\;=\;}
    - \scpr{\mathcal{L}_{s}^{\omega}(P_{s, t}^{\omega} \indicator_{\{x\}})}{f}{\mathbb{Z}^{d}}
    \;=\;
    -\bigl((P_{s, t}^{\omega})^{*}(\mathcal{L}_{s}^{\omega}f)\bigr)(x)
  \end{align*}
  for a.e.\ $s \in (-\infty,t)$.
\end{proof}

\section{Technical estimates}\label{appendix:technical}
\begin{proof}[Proof of Lemma~\ref{lem:kernel:convolution}]
  Fix some $t > 0$ and $y \in \mathbb{Z}^{d}$, and set $\mathbb{H} \ldef \{z \in \mathbb{Z}^{d} : |z| \leq |y-z|\}$.  First, we claim that
  \begin{align}\label{eq:kernel:est}
    k_{\alpha}(t/2, y - z)
    \;\leq\;
    2^{d/2 + \alpha/2}\, k_{\alpha}(t, y)
    \qquad 
    \forall\, z \in \mathbb{H}.
  \end{align}
  Indeed, in case that either $t/2 < |y|/2$ or $t/2 \geq |y-z|$, the estimate \eqref{eq:kernel:est} follows directly from the definition of the kernel $k_{\alpha}$, cf.\ \eqref{def:kernel}, together with the fact that $|y-z| \geq |y|/2$ for any $z \in \mathbb{H}$.  On the other hand, if $|y|/2 \leq t/2 < |y-z|$ we get that
  \begin{align*}
    k_{\alpha}(t/2, y - z)
    &\;=\;
    |y-z|^{-d/2 - \alpha/2}
    \\[.5ex]
    &\;\leq\;
    k_{\alpha}(t, y)\,
    \bigl(1 \vee t\bigr)^{d/2} \bigl(1 \vee |y|\bigr)^{\alpha/2}\, |y-z|^{-d/2 - \alpha/2}
    \;\leq\;
    2^{d/2 + \alpha/2}\, k_{\alpha}(t, y), 
  \end{align*}
  where we used in the last step that $0 < t/2 < |y-z|$ implies that $|y-z| \geq 1$.  
  
  Let us now address the proof of \eqref{eq:kernel:convolution}.  By symmetry, we obtain that
  \begin{align*}
    \sum_{z \in \mathbb{Z}^{d}} k_{\alpha}(t/2, z)\, k_{\alpha}(t/2, y - z)
    \;\leq\;
    2\,
    \sum_{z \in \mathbb{H}} k_{\alpha}(t/2, z)\, k_{\alpha}(t/2, y - z).
  \end{align*}
  Thus, in view of \eqref{eq:kernel:est}, the assertion of Lemma~\ref{lem:kernel:convolution} will follow, once we have shown that for any $\alpha > d$ there exists $c < \infty$ such that, for any $t > 0$,
  \begin{align}
    \label{eq:kernel:norm}
    \sum_{z \in \mathbb{Z}^{d}} k_{\alpha}(t, z)
    \;=\;
    \sum_{|z| \leq t}
    \bigl(1 \vee t\bigr)^{-d/2}\, 
    \biggl(1 \vee \frac{|z|}{\sqrt{t}}\biggr)^{\!\!-\alpha}
    +\,
    \sum_{|z| > t} |z|^{-d/2 - \alpha/2}
    \;\leq\;
    c.
  \end{align}
  Note that for any $t < 1$ the first sum consists of a single element, namely $z=0$, and is therefore bounded from above by 1. On the other hand, for any $t \geq 1$, we can rewrite the first sum on the right-hand side of \eqref{eq:kernel:norm} as
  \begin{align*}
    &\sum_{|z| \leq \sqrt{t}} t^{-d/2} +\, 
    \sum_{\sqrt{t} < |z| \leq t} 
    t^{-d/2}\, \biggl(\frac{|z|}{\sqrt{t}}\biggr)^{\!\!-\alpha}
    \\
    &\mspace{36mu}\leq\;
    t^{-d/2} \bigl|B\bigl(0, \sqrt{t}\bigr)\bigr| 
    \,+\,
    t^{-d/2} 
    \sum_{k=1}^{\infty} 
    k^{-\alpha}\,
    \bigl|
      B\bigl(0, (k + 1)\sqrt{t}\bigr) \setminus B\bigl(0, k \sqrt{t}\bigr)
    \bigr|
    \\
    &\mspace{36mu}\leq\;
    c\, \biggl(1 \,+\, \sum_{k=1}^{\infty} k^{d-\alpha-1} \biggr).
  \end{align*}
  Since $\alpha > d$, the remaining sum is bounded from above by a constant that only depends on $d$ and $\alpha$.  Let us now consider the second sum on the right-hand side of \eqref{eq:kernel:norm}. Since $\alpha > d$ it follows that
  \begin{align*}
    \sum_{|z| > t} |z|^{-d/2-\alpha/2}
    \;\leq\;
    c\, \sum_{k=1}^{\infty} k^{d/2 - \alpha/2 - 1}
    \;\leq\;
    c
    \;<\;
    \infty,
  \end{align*}
  which completes the proof of \eqref{eq:kernel:norm}.
\end{proof}

\subsubsection*{Acknowledgment}
This work was initiated while Jean-Dominique Deuschel was visi\-ting RIMS in Kyoto.  The research of Takashi Kumagai is supported by JSPS KAKENHI Grant Number JP22H00099 and by the Alexander von Humboldt Foundation.  

\bibliographystyle{plain}
\bibliography{literature}

\end{document}